\documentclass[11pt,twoside]{article}

\usepackage[english]{babel}

\usepackage{amsmath}
\usepackage{amssymb}

\usepackage{mathrsfs}
\usepackage{amsthm}

\usepackage{indentfirst}
\usepackage{color}
\usepackage{txfonts}

\usepackage{anysize}

\usepackage[colorlinks=true,
linkcolor=blue,
citecolor=red,
urlcolor=magenta,
backref=page
]{hyperref}

\allowdisplaybreaks

\newtheorem{theorem}{Theorem}[section]
\newtheorem{lemma}[theorem]{Lemma}

\newtheorem{proposition}[theorem]{Proposition}

\theoremstyle{definition}
\newtheorem{remark}[theorem]{Remark}
\newtheorem{definition}[theorem]{Definition}

\numberwithin{equation}{section}

\begin{document}

\title{\bf\Large Sharp Fractional
Capacitary-Distance Hardy Inequalities\footnotetext{\hspace{-0.35cm} 2020 {\it
Mathematics Subject Classification}. Primary 26D10; Secondary 31C45, 35J92, 46E35.
\endgraf
{\it Key words and phrases}. Hardy inequality, capacitary distance, nonlinear resolvent,
Sobolev capacity, tree capacity.
\endgraf
This project is partially supported by the National Natural Science
Foundation of China (Grant Nos.~12431006, 12371093, and 12501118),
the Beijing Natural Science Foundation (Grant No.~1262011),
the Fundamental Research Funds for the Central Universities
(Grant No.~2253200028), and NSERC of Canada (Grant No.~202979).
R.~L.~F. acknowledges partial support from the US National Science
Foundation (Grant No.~DMS-1954995) and the German Research Foundation
through EXC-2111-390814868, TRR~352--Project-ID~470903074,
and FR~2664/3-1.
}}
\author{Yiqun Chen, Rupert L. Frank, Jie Xiao, Dachun Yang,\\Wen Yuan and
Yangyang Zhang
}
\maketitle

\vspace{-0.8cm}

\begin{center}
\begin{minipage}{13cm}
{\small {\bf Abstract:}\quad
Let $d\in\mathbb N$, $p\in(1,\infty)$, $p':=\frac{p}{p-1}$, and
$s\in(0,\min\{1+\frac1p,\frac dp\})$.
For an open set $\Omega\subset\mathbb{R}^d$, let $F:=\mathbb{R}^d\setminus\Omega$.
For any $\alpha\in(0,1]$ and $x\in\Omega$, let
\begin{align*}
d_{\alpha,s,p}(x)
:=\inf\left\{r\in(0,\infty):
\frac{\operatorname{cap}_{s,p}(\overline{F\cap B(x,r)})}
{\operatorname{cap}_{s,p}(\overline{B(x,r)})}\ge\alpha\right\},
\end{align*}
where $\operatorname{cap}_{s,p}$ is the capacity associated with the
fractional Sobolev energy $[\,\cdot\,]_{s,p}^p$.
In this article, we prove that there exists a positive constant $C$,
depending only on $d$, $s$, and $p$, such that, for any
$\alpha\in(0,1]$ and $u\in C_{\rm c}^\infty(\Omega)$,
\begin{align*}
\int_\Omega\frac{|u(x)|^p}{d_{\alpha,s,p}(x)^{sp}}\,dx
\le C\alpha^{-p'}[u]_{s,p}^p.
\end{align*}
We also show that both the exponent $p'$ in $\alpha^{-p'}$ as
$\alpha\to0^+$ over the class of all open sets
and the range $s\in(0,1+\frac1p)$ within $sp<d$ are sharp.
Our proofs use variational minimization, equilibrium measures on trees,
and weighted Hardy inequalities.
}
\end{minipage}
\end{center}

\vspace{0.15cm}

\tableofcontents

\vspace{0.2cm}

\section{Introduction}

Hardy-type inequalities originate in the classical one-dimensional
inequalities of Hardy \cite{Hardy1920,Hardy1925} and have become
important tools in harmonic analysis, potential theory, and partial
differential equations. Let $d\in\mathbb N$ and $p\in(1,\infty)$.
A basic question is to determine the open sets $\Omega\subset\mathbb{R}^d$
for which there exists a positive constant $C$ such that, for any
$u\in C_{\rm c}^\infty(\Omega)$ (the set of all infinitely differentiable functions
with compact support in $\Omega$),
\begin{align}\label{eq:intro-classical-p-hardy}
\int_\Omega\frac{|u(x)|^p}{\operatorname{dist}(x,\partial\Omega)^p}\,dx
\le C\int_\Omega|\nabla u(x)|^p\,dx.
\end{align}
Ne\v{c}as \cite{Necas1962} studied such inequalities on Lipschitz
domains. Subsequent work of Ancona \cite{AnconaHardy},
Lewis \cite{LewisUniformlyFat}, and Wannebo \cite{Wannebo1990}
showed that capacitary thickness of the complement can replace
boundary regularity. Capacitary characterizations further explain
the role of the geometry of $\Omega$; see
\cite[pp.~733--735]{MazyaSobolev} and
\cite[Theorem~2.1]{KinnunenKorteHardy}.
Fractional analogues based on difference quotients and nonlocal
energies have also been studied extensively; see, for example,
\cite{Dyda2004,FrankSeiringer2008,imv25}.
For fractional Hardy inequalities with homogeneous weights,
see, for example, \cite{Roy2026}.

The Euclidean distance in \eqref{eq:intro-classical-p-hardy} records
how close a point is to the complement, but does not measure the
capacitary size of that complement. In
\cite[Problem~8, Section~3.8]{MazProblems}, Maz'ya asked, when $p=2$,
whether replacing the Euclidean distance by a capacitary distance makes
the Hardy inequality valid on every domain.
Chen et al. \cite{CXYZ} answered this question
for Newtonian capacity by establishing the corresponding estimate
on every open set and proving the optimal dependence $\alpha^{-2}$
on the capacity density $\alpha$.
The aim of the present article is to extend this capacitary-distance
principle \cite[Theorem 1.3]{CXYZ} to local, fractional, and higher-order $p$-energies, and to
determine the sharp exponent of the density parameter $\alpha$
and the admissible range of order $s$.

For any $s\in(0,1)$, any $p\in(1,\infty)$, and any measurable function
$u$, let
\begin{align*}
[u]_{s,p}^p
:=\iint_{\mathbb{R}^d\times\mathbb{R}^d}
\frac{|u(x)-u(y)|^p}{|x-y|^{d+sp}}\,dx\,dy.
\end{align*}
For any $u\in\mathcal S(\mathbb{R}^d)$ and $k\in\mathbb Z_+$, let
$D^ku:=(\partial^\beta u)_{\beta\in\mathbb Z_+^d,\,|\beta|=k}$,
where $\mathcal S(\mathbb{R}^d)$ denotes the
\emph{Schwartz space}. For any $s\in(0,\infty)$, the
\emph{notation $\lfloor s\rfloor$} denotes the
largest integer not greater than $s$.
For any $s\in[1,\infty)$, $p\in(1,\infty)$, and
$u\in\mathcal S(\mathbb{R}^d)$, let $m:=\lfloor s\rfloor$ and
\begin{align*}
[u]_{s,p}^p:=
\begin{cases}
\|D^m u\|_{L^p(\mathbb{R}^d)}^p,&s\in\mathbb N,\\[4pt]
\displaystyle\iint_{\mathbb{R}^d\times\mathbb{R}^d}
\frac{|D^m u(x)-D^m u(y)|^p}{|x-y|^{d+(s-m)p}}\,dx\,dy,
&s\notin\mathbb N.
\end{cases}
\end{align*}
For the corresponding higher-order fractional Sobolev spaces, see,
for instance, \cite[Subsection~2.2.2]{Triebel}.
Throughout this article, functions in $C_{\rm c}^\infty(\Omega)$
are always understood to be extended by zero to $\mathbb R^d$.
For any $s\in(0,\infty)$, any $p\in(1,\infty)$, and any compact set
$K\subset\mathbb{R}^d$, we consider the following \emph{absolute capacity}
\begin{align}\label{ac}
\operatorname{cap}_{s,p}(K)
:=\inf\left\{[\varphi]_{s,p}^p:
\varphi\in C_{\rm c}^\infty(\mathbb{R}^d),\ \varphi\ge1\text{ on an open neighbourhood of }K\right\}.
\end{align}
Let $\Omega\subset\mathbb{R}^d$ be open and $F:=\mathbb{R}^d\setminus\Omega$.
For any $s\in(0,\infty)$, $p\in(1,\frac ds)$,
$\alpha\in(0,\infty)$, and $x\in\Omega$, let
\begin{align*}
d_{\alpha,s,p}(x)
&:=\inf\left\{r\in(0,\infty):
\frac{\operatorname{cap}_{s,p}(\overline{F\cap B(x,r)})}
{\operatorname{cap}_{s,p}(\overline{B(x,r)})}\ge\alpha\right\}.
\end{align*}
We use the conventions that $\inf\varnothing:=\infty$ and $\infty^{-a}:=0$ for any
$a\in(0,\infty)$.
For any $p\in(1,\infty)$,
$s\in(0,\min\{1+\frac1p,\frac dp\})$, and $\alpha\in(0,\infty)$, let
\begin{align*}
\mathcal H_{d,s,p}(\alpha)
:=\sup_{\genfrac{}{}{0pt}{}{\Omega\subset\mathbb{R}^d\ \mathrm{open}}{0\ne u\in C_{\rm c}^\infty(\Omega)}}
\frac{1}{[u]_{s,p}^p}\int_\Omega
\frac{|u(x)|^p}{d_{\alpha,s,p}(x)^{sp}}\,dx.
\end{align*}
In what follows, let $p':=\frac{p}{p-1}$ be the
\emph{conjugate index} of $p\in(1,\infty]$.
Our first main result is as follows.

\begin{theorem}\label{thm:generalp}
Let $p\in(1,\infty)$ and
$s\in(0,\min\{1+\frac1p,\frac dp\})$.
Then the following statements hold.
\begin{itemize}
\item[{\rm(i)}]
There exists a positive constant $C$ such that, for any open set
$\Omega\subset\mathbb{R}^d$, any $\alpha\in(0,1]$, and any
$u\in C_{\rm c}^\infty(\Omega)$,
\begin{align}\label{eq:generalp-main}
\int_\Omega\frac{|u(x)|^p}{d_{\alpha,s,p}(x)^{sp}}\,dx
\le C\alpha^{-p'}[u]_{s,p}^p.
\end{align}

\item[{\rm(ii)}]
There exist a positive constant $c$ and $\alpha_*\in(0,1)$
such that, for any $\alpha\in(0,\alpha_*]$,
\begin{align*}
\mathcal H_{d,s,p}(\alpha)\ge c\alpha^{-p'}.
\end{align*}
\end{itemize}
\end{theorem}

In the above theorem, the restriction $s<1+\frac1p$ is sharp within the subcritical range
$sp<d$. More precisely, the following theorem gives bounded-domain
counterexamples whenever the complementary range is nonempty.

\begin{theorem}\label{thm:counterexample}
Let $p\in(1,d-1)$ and
$s\in[1+\frac1p,\frac dp)$.
Then there exist a bounded open set $\Omega\subset\mathbb{R}^d$,
$\alpha_0\in(0,1]$, and a sequence
$\{u_j\}_{j\in\mathbb N}\subset C_{\rm c}^\infty(\Omega)$ such that
\begin{align*}
\sup_{j\in\mathbb N}[u_j]_{s,p}<\infty\quad and\quad
\lim_{j\to\infty}\int_\Omega
\frac{|u_j(x)|^p}{d_{\alpha_0,s,p}(x)^{sp}}\,dx=\infty.
\end{align*}
\end{theorem}

It is known that the absolute homogeneous capacity is
degenerate in the case $sp\ge d$;
see, for instance, \cite[Lemmas~4.1 and~5.1]{BrascoGomezCastroVazquez2021}.
To remove this restriction, we need to work with relative capacity.
For related work on relative capacities, see, for example,
\cite{AdamsXiao2025}.
For any $s\in(0,\infty)$, any $p\in(1,\infty)$, any open set
$U\subset\mathbb R^d$, and any compact set $K\subset U$,
we consider the following \emph{relative capacity}
\begin{align*}
\operatorname{cap}_{s,p}(K;U)
:=\inf\left\{[\varphi]_{s,p}^p:
\varphi\in C_{\rm c}^\infty(U),\ \varphi\ge1
\text{ on an open neighbourhood of }K\right\}.
\end{align*}
For any open set $\Omega\subset\mathbb R^d$, any
$\alpha\in(0,\infty)$, and any $x\in\Omega$, let
\begin{align}\label{eq:relative-distance}
\widetilde d_{\alpha,s,p}(x):=\inf\left\{r\in(0,\infty):
\frac{\operatorname{cap}_{s,p}(\overline{F\cap B(x,r)};B(x,2r))}
{\operatorname{cap}_{s,p}(\overline{B(x,r)};B(x,2r))}\ge\alpha\right\}.
\end{align}
For any $p\in(1,\infty)$, $s\in(0,\infty)$, and
$\alpha\in(0,\infty)$, let
\begin{align*}
\widetilde{\mathcal H}_{d,s,p}(\alpha)
:=\sup_{\genfrac{}{}{0pt}{}{\Omega\subset\mathbb R^d\ \mathrm{open}}{0\ne u\in C_{\rm c}^\infty(\Omega)}}
\frac{1}{[u]_{s,p}^p}\int_\Omega
\frac{|u(x)|^p}{\widetilde d_{\alpha,s,p}(x)^{sp}}\,dx.
\end{align*}
For $p\in(1,\infty)$, let
$\mathcal A_{d,p}:=(0,1+\frac1p)\cup\bigcup_{k\in\mathbb N}(k+\frac dp,k+1+\frac1p)$.
\begin{theorem}\label{thm:relative-capacity}
Let $p\in(1,\infty)$ and
$s\in\mathcal A_{d,p}$.
Then the following assertions hold.
\begin{itemize}
\item[{\rm(i)}]
There exists a positive constant $C$
such that, for any open set
$\Omega\subset\mathbb R^d$, any $\alpha\in(0,1]$, and any
$u\in C_{\rm c}^\infty(\Omega)$,
\begin{align}\label{eq:relative-main}
\int_\Omega\frac{|u(x)|^p}{\widetilde d_{\alpha,s,p}(x)^{sp}}\,dx
\le C[u]_{s,p}^p
\begin{cases}
\alpha^{-p'}& if\ s\in(0,\frac dp],\\
1 &if\ s\in(\frac dp,\infty).
\end{cases}
\end{align}
\item[{\rm(ii)}]
There exist a positive constant $c$ and $\alpha_*\in(0,1)$,
depending only on $d$, $s$, and $p$,
such that, for any $\alpha\in(0,\alpha_*]$,
\begin{align*}
\widetilde{\mathcal H}_{d,s,p}(\alpha)\ge c
\begin{cases}
\alpha^{-p'} &if\ s\in(0,\frac dp],\\
1 &if\ s\in(\frac dp,\infty).
\end{cases}
\end{align*}
\end{itemize}
\end{theorem}

The following theorem shows the sharpness of the admissible range
$s\in\mathcal A_{d,p}$

\begin{theorem}\label{thm:relative-classification}
Let $p\in(1,\infty)$ and $s\in(0,\infty)\setminus\mathcal A_{d,p}$.
Then there exist a bounded open set $\Omega\subset\mathbb R^d$,
$\alpha_0\in(0,1)$, and a sequence
$\{u_j\}_{j\in\mathbb N}\subset C_{\rm c}^\infty(\Omega)$ such that
\begin{align*}
\sup_{j\in\mathbb N}[u_j]_{s,p}<\infty
\quad\text{and}\quad
\lim_{j\to\infty}\int_\Omega
\frac{|u_j(x)|^p}{\widetilde d_{\alpha_0,s,p}(x)^{sp}}\,dx=\infty.
\end{align*}
\end{theorem}

\begin{remark} We have the following observations.
\begin{itemize}
\item[{\rm(i)}]
The dependence on $\alpha$ in Theorems~\ref{thm:generalp}(i)
and~\ref{thm:relative-capacity}(i) is sharp as $\alpha\to0^+$.
The assumption $sp<d$ in Theorems~\ref{thm:generalp} and \ref{thm:counterexample}
is needed because, when $sp\ge d$,
$\operatorname{cap}_{s,p}(K)=0$
for every compact set $K\subset\mathbb R^d$.
Within $sp<d$, Theorem~\ref{thm:counterexample} shows that
the restriction $s<1+\frac1p$ is also sharp.

\item[{\rm(ii)}]
By replacing absolute capacity with relative capacity,
Theorem~\ref{thm:relative-capacity} removes the restriction $sp<d$
in Theorem~\ref{thm:generalp}.
Moreover, when $sp<d$, Theorems~\ref{thm:generalp}
and \ref{thm:relative-capacity} in this case are equivalent up to a fixed
rescaling of $\alpha$; see \eqref{eq:relative-best-comparison} below.

\item[{\rm(iii)}]
If $p\le d-1$, then $\mathcal A_{d,p}=(0,1+\frac1p)$.
If $p>d-1$, then $\mathcal A_{d,p}$ also contains the nonempty
intervals $(k+\frac dp,k+1+\frac1p)$ for any $k\in\mathbb N$.
\end{itemize}
\end{remark}

We next briefly describe the main ideas used in this article.

For $s=1$, Lemma~\ref{lem:plateau} reduces the proof to estimating
the weighted measure of $\{f=1\}$ for a function $f$ taking values in $[0,1]$.
We then construct comparison functions $v_t$ at different scales by
minimizing an energy functional \eqref{ef} that balances the gradient energy against
the deviation from $f$. Using this minimizing property,
we find that the gradient energy
outside $\{f=1\}$ controls the overall decrease of $v_t$ from 1
on this set.
By varying $t$, we choose for each point a scale adapted to its
capacitary distance, at which the capacity condition forces
$v_t$ to decrease by at least a constant multiple of
$\alpha^{\frac1{p-1}}$.
Combining this lower bound with the energy estimate, we obtain
the Hardy inequality with the sharp factor $\alpha^{-p'}$.

For $s\neq1$, we first select dyadic cubes associated with the scale of
the capacitary distance, and then decompose the function into high-
and low-frequency terms on
each selected cube. Young's inequality and the Besov--Sobolev
characterization control the high-frequency terms.
For the low-frequency terms, we use affine approximations and
control the cubes on which the linear part dominates by
Lemma~\ref{lem:affine-gradient}.
On the cubes where the constant part dominates, we reduce the
estimate to the low-frequency values on the zero set.
Since $u=0$ on this zero set,
it follows that these values can be written
as sums of high-frequency blocks.
Using weights constructed from a tree-capacity minimizer and
its associated measure, together with the weighted Hardy inequality
in Lemma~\ref{lem:rayhardy} along nested cubes, we obtain the
desired estimate in this case.

To prove the sharpness of power $p'$ of $\alpha^{-p'}$, we use a compact regular
set of Hausdorff dimension $d-sp$ and place small disjoint balls at separated
points of that set. The resulting obstacles have capacity density
bounded below by a logarithmic factor. We construct smooth test functions that decrease from one to zero
through a sequence of smooth steps towards these obstacles,
yielding the sharp dependence of the Hardy constant on the
capacity parameter.
For the counterexamples in Theorem~\ref{thm:counterexample}, we instead
use a regular set of Hausdorff dimension $d-p(s-1)$ contained in a hyperplane.
Its capacity is quantitatively large at small scales.
Smooth cutoffs of a function vanishing linearly on the hyperplane
have uniformly bounded energy, while their Hardy integrals
diverge logarithmically.

To remove the restriction $sp<d$, we replace absolute capacity by
relative capacity and determine the complete range of admissible
order $s$. For the additional admissible orders, we prove the relative
capacitary-distance Hardy inequality
by iterating a weighted first-order Hardy inequality and reducing
to the lower-order case.
For orders outside the admissible range, we construct counterexamples
by replacing the linear function in the preceding construction
with a polynomial.

The remainder of this article is organized as follows.
In Section~\ref{sec:local}, we establish the capacitary barrier and
resolvent comparison, and prove Theorem~\ref{thm:generalp}(i) for $s=1$.
Section~\ref{sec:higher} contains the dyadic preliminaries, the affine
estimate, and the proof of Theorem~\ref{thm:generalp}(i) for $s\ne1$.
In Section~\ref{sec:alpha-sharpness}, we construct the critical-dimensional
obstacles and smooth test functions, and prove
Theorem~\ref{thm:generalp}(ii).
Section~\ref{sec:counterexamples} is devoted to the proof of
Theorem~\ref{thm:counterexample}.
Finally, in Section~\ref{sec:relative-capacity}, we compare absolute and
relative capacities, give the extension to $sp\ge d$, and prove
Theorems~\ref{thm:relative-capacity} and~\ref{thm:relative-classification}.

We end this introduction by making some notational conventions.
Let $\mathbb N:=\{1,2,\ldots\}$ and
$\mathbb Z_+:=\{0,1,2,\ldots\}$.
The symbols $\mathbb R$, $\mathbb Z$, and $\mathbb Q$ denote the sets
of real numbers, integers, and rational numbers, respectively.
For any $x,y\in\mathbb{R}^d$, let $x\cdot y$ and $|x|$ denote the Euclidean
inner product and norm, respectively. For any $x\in\mathbb{R}^d$ and
$r\in(0,\infty)$, let $B(x,r):=\{y\in\mathbb{R}^d:|y-x|<r\}$.
For any Lebesgue measurable set $A\subset\mathbb{R}^d$, let $|A|$ denote its
Lebesgue measure and $\mathbf1_A$ its characteristic function.
For any $x\in\mathbb{R}^d$ and $A\subset\mathbb{R}^d$, let
$\operatorname{dist}(x,A):=\inf_{y\in A}|x-y|$, with $\operatorname{dist}(x,\varnothing):=\infty$.
For any $a\in\mathbb{R}$, let $a_+:=\max\{a,0\}$.
For any open set $O\subset\mathbb{R}^d$, let
\begin{align*}
\operatorname{cap}_{s,p}(O)
:=\sup\left\{\operatorname{cap}_{s,p}(K):K\subset O,\ K\text{ is compact}\right\}.
\end{align*}
For any set $A\subset\mathbb{R}^d$, its \emph{outer capacity} is defined by
\begin{align*}
\operatorname{cap}_{s,p}(A)
:=\inf\left\{\operatorname{cap}_{s,p}(O):A\subset O\subset\mathbb{R}^d,
\ O\text{ is open}\right\}.
\end{align*}
Let $\mathcal S'(\mathbb{R}^d)$ denote the \emph{space of tempered distributions}.
For any $f\in\mathcal S(\mathbb{R}^d)$,
its Fourier transform $\widehat f$ is defined by setting,
for any $\xi\in\mathbb{R}^d$,
$$\widehat f(\xi):=\frac 1{(2\pi)^{\frac d2}}\int_{\mathbb{R}^d}
e^{-ix\cdot\xi}f(x)\,dx.$$
For any dyadic cube $Q$, let $x_Q$, $\ell(Q)$, $\widehat Q$, and
$j(Q):=-\log_2\ell(Q)\in\mathbb Z$ denote its center, edge length,
parent, and generation, respectively. For any $\lambda\in(0,\infty)$,
let $\lambda Q:=x_Q+\lambda(Q-x_Q)$.
For any open set $U\subset\mathbb{R}^d$, let
$$W^{1,p}(U):=\left\{v\in L^p(U):\nabla v\in L^p\left(U;\mathbb{R}^d\right)\right\}$$
and $W^{1,p}_0(U):=\overline{C_{\rm c}^\infty(U)}^{\,W^{1,p}(U)}$,
where $\nabla v$ denotes the weak gradient and the closure is taken
with respect to the norm
$$\|v\|_{W^{1,p}(U)}:=\left[\|v\|_{L^p(U)}^p+
\|\nabla v\|_{L^p(U;\mathbb{R}^d)}^p\right]^{\frac1p}$$
for any $v\in W^{1,p}(U)$.
For $p\in(1,\infty)$, a property holds \emph{$p$-quasi-everywhere}
($p$-q.e.) if its exceptional set $N\subset\mathbb R^d$ satisfies
$$\operatorname{cap}_{1,p}(N\cap B(x,r);B(x,2r))=0$$
for any $x\in\mathbb R^d$ and $r\in(0,\infty)$;
see, for instance, \cite[Section~3.1, p.~143]{KilpelainenMaly1994}.
For $p\in(1,d)$, this is equivalent to
$\operatorname{cap}_{1,p}(N)=0$.
Pointwise conditions on $W^{1,p}$ functions are understood in terms
of their $p$-quasicontinuous representatives, which are uniquely
determined up to sets of $p$-capacity zero; see, for instance,
\cite[Theorem~4.14]{HeinonenKilpelainenMartio}.
We shall also use the standard Sobolev characterization of
$1$-capacity. If $K\subset U$ is compact and $U\subset\mathbb R^d$
is open, then
\begin{align*}
\operatorname{cap}_{1,p}(K;U)
=\inf\left\{
\int_U|\nabla v|^p\,dx:
 v\in W^{1,p}_0(U),\
 \widetilde v\ge1\quad p\text{-q.e. on }K
\right\},
\end{align*}
where $\widetilde v$ denotes the $p$-quasicontinuous representative
of $v$. When $p<d$, the corresponding absolute characterization is
\begin{align*}
\operatorname{cap}_{1,p}(K)
=\inf\left\{
\int_{\mathbb R^d}|\nabla v|^p\,dx:
 v\in W^{1,p}(\mathbb R^d),\ v\text{ has compact support},\
 \widetilde v\ge1\quad p\text{-q.e. on }K
\right\}.
\end{align*}
These formulations are equivalent to the smooth definitions above
by approximation and quasicontinuity.
The symbols $C$ and $c$ denote positive constants that may change
from line to line. We write $f\lesssim g$ when $f\le Cg$,
$f\gtrsim g$ when $g\lesssim f$, and $f\sim g$ when both inequalities hold.
Finally, in all proofs we
consistently retain the notation introduced in the
original theorem (or related statement).

\section{Proof of Theorem \ref{thm:generalp}(i) in the Case $s=1$}
\label{sec:local}

In this section, we prove Theorem \ref{thm:generalp}(i) in the case $s=1$.
To this end, we need some technical lemmas.

\begin{lemma}\label{lem:plateau}
Let $\nu$ be a nonnegative measure on $\Omega$, and let $A\in(0,\infty)$.
Assume that, for any $f\in W^{1,p}_0(\Omega)$ satisfying $0\le f\le1$,
\begin{align}\label{eq:plateau-hyp-local}
\nu(\{x\in\Omega:\,f(x)=1\})\le A[f]_{1,p}^p.
\end{align}
Then, for any $u\in C_{\rm c}^\infty(\Omega)$,
\begin{align*}
\int_\Omega |u(x)|^p\,d\nu(x)\le4^p A[u]_{1,p}^p.
\end{align*}
\end{lemma}

\begin{proof}
Let $u\in C_{\rm c}^\infty(\Omega)$ and
$f_k:=2^{-k}\min\{(|u|-2^k)_+,2^k\}$ for any $k\in\mathbb{Z}$. Using the Sobolev chain rule and
$\nabla|u|=0$ almost everywhere on
$\{x\in\Omega:|u(x)|=t\}$ for any $t\in[0,\infty)$, we obtain
\begin{align}\label{eq:local-truncation-energy}
\sum_{k\in\mathbb{Z}}2^{kp}[f_k]_{1,p}^p
=\sum_{k\in\mathbb{Z}}
\int_{\{x\in\Omega:\,2^k<|u(x)|<2^{k+1}\}}
|\nabla|u|(x)|^p\,dx=\int_\Omega|\nabla|u|(x)|^p\,dx\le[u]_{1,p}^p.
\end{align}
By the definition of $f_k$, \eqref{eq:plateau-hyp-local}, the
homogeneity of $[\,\cdot\,]_{1,p}^p$, and
\eqref{eq:local-truncation-energy}, we conclude that
\begin{align*}
\int_\Omega |u(x)|^p\,d\nu(x)
&\le\sum_{k\in\mathbb{Z}}2^{p(k+2)}
\nu\left(\left\{x\in\Omega:|u(x)|\ge2^{k+1}\right\}\right)\\
&=\sum_{k\in\mathbb{Z}}2^{p(k+2)}\nu(\{x\in\Omega:f_k(x)=1\})
\le A\sum_{k\in\mathbb{Z}}2^{p(k+2)}[f_k]_{1,p}^p\\
&=4^p A\sum_{k\in\mathbb{Z}}2^{kp}[f_k]_{1,p}^p
\le4^p A[u]_{1,p}^p.
\end{align*}
This completes the proof of Lemma~\ref{lem:plateau}.
\end{proof}

For $sp<d$, let
$\kappa_{d,s,p}:=\operatorname{cap}_{s,p}(\overline{B(0,1)}).$
By translation invariance and scaling,
$\operatorname{cap}_{s,p}(\overline{B(x,r)})
=\kappa_{d,s,p}r^{d-sp}$
for every $x\in\mathbb R^d$ and $r\in (0,\infty)$.
The following lemma shows the measurability of the capacitary distance.

\begin{lemma}\label{lem:distance-measurable}
Let $s\in(0,\infty)$, $p\in(1,\frac ds)$, and
$\Omega\subset\mathbb{R}^d$ be open. Then, for any
$\alpha\in(0,\infty)$, $d_{\alpha,s,p}$ is Lebesgue measurable.
\end{lemma}

\begin{proof}
We follow the proof of \cite[Lemma~2.7]{CXYZ} with the metric
therein replaced by
$\rho(\phi,\psi):=\|\phi-\psi\|_{L^\infty(\mathbb R^d)}
+[\phi-\psi]_{s,p}$
and the energies therein replaced by $e_j:=[\phi_j]_{s,p}^p$.
Using the derivative and difference-quotient representations of
$[\,\cdot\,]_{s,p}$ and the separability of
$L^p$ spaces, we obtain the separability of
$C_{\rm c}^\infty(\mathbb R^d)$ with respect to $\rho$.
By continuity and approximation, for any compact set $K\subset\mathbb R^d$,
we may replace the condition $\phi\ge1$ on an open neighbourhood of $K$
by $\phi\ge1$ on $K$ in the definition \eqref{ac} of $\operatorname{cap}_{s,p}(K)$.
Therefore, the same countable-test-function and projection
argument as that used in the proof of \cite[Lemma~2.7]{CXYZ} gives the Lebesgue
measurability of $d_{\alpha,s,p}$.
This completes the proof of Lemma~\ref{lem:distance-measurable}.
\end{proof}

For any $t\in(0,\infty)$, $g\in L^p(\mathbb{R}^d)$, and
$v\in W^{1,p}(\mathbb{R}^d)$, let
\begin{align}\label{ef}
\mathcal J_t(v;g):=\int_{\mathbb{R}^d}
\left[|\nabla v(x)|^p+\frac{|v(x)-g(x)|^p}{t}\right]\,dx.
\end{align}

\begin{lemma}\label{lem:local-barrier}
Let $p\in(1,\infty)$, $z\in\mathbb R^d$, $R\in(0,\infty)$,
$E\subset\overline{B(z,R)}$ be compact, and $D=B(z,8R)$.
Then there exists $h\in W^{1,p}_0(D)$ satisfying $0\le h\le1$
and $h=1$ $p$-quasi-everywhere on $E$ such that
$\mathcal J_{R^p}(h;0)=\inf \mathcal J_{R^p}(v;0)$,
where the infimum is taken over all $v\in W^{1,p}_0(D)$ satisfying
$v\ge1$ $p$-quasi-everywhere on $E$.
Moreover, there exists a positive constant $C$, depending only on
$d$ and $p$, such that
\begin{align}\label{eq:barrier-energy-oscillation}
\mathcal J_{R^p}(h;0)
\le CR^{d-p}
\left[\operatorname*{ess\,inf}_{y\in B(z,4R)}h(y)\right]^{p-1}.
\end{align}
\end{lemma}

\begin{proof}
Let $$\mathcal K:=\left\{v\in W^{1,p}_0(D):v\ge1
\text{ $p$-quasi-everywhere on }E\right\}.$$
Applying \cite[Corollary~3.23]{Brezis2011} to
$W^{1,p}_0(D)$ with $A:=\mathcal K$ and
$\varphi:=\mathcal J_{R^p}(\,\cdot\,;0)|_{\mathcal K}$,
we obtain a unique $h\in\mathcal K$ satisfying
\begin{align}\label{eq:barrier-new-minimum}
\mathcal J_{R^p}(h;0)=\inf_{v\in\mathcal K}\mathcal J_{R^p}(v;0).
\end{align}
For any $t\in\mathbb{R}$, let $T(t):=\min\{1,\max\{0,t\}\}$.
By the Sobolev chain rule and the definition of $\mathcal K$,
we have $T(h)\in\mathcal K$ and
\begin{align*}
\mathcal J_{R^p}(T(h);0)
&=\int_{\{x\in D:\,0<h(x)<1\}}|\nabla h(x)|^p\,dx
+R^{-p}\int_D|T(h(x))|^p\,dx\le \mathcal J_{R^p}(h;0).
\end{align*}
Using this and the uniqueness of the minimizer, we conclude that
$h=T(h)$ in $W^{1,p}_0(D)$, $0\le h\le1$ a.e. in $D$, and
$h=1$ $p$-quasi-everywhere on $E$.

By $h\in\mathcal K$, we find that, for any
$t\in[0,\infty)$ and $\phi\in C_{\rm c}^\infty(D)$ satisfying $\phi\ge0$,
$h+t\phi\in\mathcal K$.
From this and \eqref{eq:barrier-new-minimum}, we deduce that, for any such $\phi$,
\begin{align}\label{eq:barrier-new-supersolution}
0\le\left.\frac1p\frac{d}{dt}\mathcal J_{R^p}(h+t\phi;0)\right|_{t=0}
=\int_D|\nabla h|^{p-2}\nabla h\cdot\nabla\phi\,dx
+R^{-p}\int_Dh^{p-1}\phi\,dx.
\end{align}
By the definition of $\mathcal K$, we have $h+t\phi\in\mathcal K$
for any $t\in\mathbb{R}$ and $\phi\in W^{1,p}_0(D)$ satisfying
$\phi=0$ $p$-quasi-everywhere on $E$.
Using this and \eqref{eq:barrier-new-minimum}, we obtain, for any such $\phi$,
\begin{align}\label{eq:barrier-new-stationarity}
\int_D|\nabla h|^{p-2}\nabla h\cdot\nabla\phi\,dx
+R^{-p}\int_Dh^{p-1}\phi\,dx
=\left.\frac1p\frac{d}{dt}\mathcal J_{R^p}(h+t\phi;0)\right|_{t=0}=0.
\end{align}
From this with $\phi\in C_{\rm c}^\infty(D\setminus E)$, we infer that
$-\Delta_p h+R^{-p}h^{p-1}=0$ weakly in $D\setminus E$, where
$\Delta_p h:=\operatorname{div}(|\nabla h|^{p-2}\nabla h)$ denotes
the \emph{$p$-Laplacian} of $h$.
Under the rescaling $w(y):=h(z+Ry)$, this equation becomes
$-\Delta_p w+w^{p-1}=0.$
This equation satisfies the structural hypotheses of the weak Harnack
and Harnack theorems in \cite[Theorems 1.1 and 1.2]{Trudinger},
with constants depending only on $d$ and $p$. Hence all Harnack
constants below are independent of $R$, $z$, and $E$.
Let $\eta\in C_{\rm c}^\infty(B(z,3R))$ satisfy $0\le\eta\le1$,
$\eta(x)=1$ for any $x\in B(z,2R)$, and
$|\nabla\eta(x)|\lesssim R^{-1}$ for any $x\in\mathbb{R}^d$.
Using $\eta-h=0$ $p$-quasi-everywhere on $E$ and using
\eqref{eq:barrier-new-stationarity} with $\phi:=\eta-h$, we obtain
\begin{align}\label{eq:barrier-new-identity}
\mathcal J_{R^p}(h;0)=\int_D|\nabla h|^{p-2}\nabla h\cdot\nabla\eta\,dx
+R^{-p}\int_Dh^{p-1}\eta\,dx.
\end{align}
Let $m:=\operatorname*{ess\,inf}_{x\in B(z,4R)}h(x)$.
From \eqref{eq:barrier-new-supersolution} and the weak Harnack inequality
(see, for instance, \cite[Theorem~1.2]{Trudinger}),
it follows that
\begin{align}\label{eq:barrier-new-moment}
\int_{B(z,4R)}h(x)^{p-1}\,dx\lesssim R^d m^{p-1}.
\end{align}
Using $E\subset\overline{B(z,R)}$,
we find that, for any $y\in\overline{B(z,\frac72R)}\setminus B(z,\frac32R)$,
$$B\left(y,\frac R4\right)\subset B(z,4R)\setminus\overline{B(z,R)}
\subset D\setminus E.$$
By this, \eqref{eq:barrier-new-stationarity}, the Harnack inequality
(see, for instance, \cite[Theorem~1.1]{Trudinger}),
and \eqref{eq:barrier-new-moment},
we conclude that, for any $y\in\overline{B(z,\frac72R)}\setminus B(z,\frac32R)$,
\begin{align*}
\left[\operatorname*{ess\,sup}_{x\in B(y,\frac R8)}h(x)\right]^{p-1} &\lesssim\left[\operatorname*{ess\,inf}_{x\in B(y,\frac R8)}h(x)\right]^{p-1} \le\frac1{|B(0,\frac R8)|}\int_{B(z,4R)}h(x)^{p-1}\,dx \lesssim m^{p-1}.
\end{align*}
From this and
$\overline{B(z,\frac72R)}\setminus B(z,\frac32R)
\subset\bigcup_{i=1}^M B(y_i,\frac R8)$
for some $M\in\mathbb N$ and
$\{y_i\}_{i=1}^M\subset\overline{B(z,\frac72R)}\setminus B(z,\frac32R)$,
we infer that
\begin{align}\label{eq:barrier-new-annulus}
\operatorname*{ess\,sup}_{x\in\overline{B(z,\frac72R)}\setminus B(z,\frac32R)}
h(x)\lesssim m.
\end{align}
Let $\zeta\in C_{\rm c}^\infty(B(z,\frac72R)\setminus
\overline{B(z,\frac32R)})$ satisfy $0\le\zeta\le1$,
$\zeta(x)=1$ for any $x\in\overline{B(z,3R)}\setminus B(z,2R)$,
and $|\nabla\zeta(x)|\lesssim R^{-1}$ for any $x\in\mathbb{R}^d$.
Using $h\zeta^p\in W^{1,p}_0(D)$, $h\zeta^p=0$ on $E$,
and \eqref{eq:barrier-new-stationarity} with $\phi:=h\zeta^p$,
we obtain
\begin{align*}
\int_D\zeta^p|\nabla h|^p\,dx
+R^{-p}\int_Dh^p\zeta^p\,dx
=-p\int_Dh\zeta^{p-1}|\nabla h|^{p-2}
\nabla h\cdot\nabla\zeta\,dx.
\end{align*}
By this and H\"older's inequality, we conclude that
\begin{align*}
\int_D\zeta^p|\nabla h|^p\,dx
\le p\left(\int_D\zeta^p|\nabla h|^p\,dx\right)^{\frac{p-1}{p}}
\left(\int_Dh^p|\nabla\zeta|^p\,dx\right)^{\frac1p},
\end{align*}
which, together with the properties of $\zeta$ and \eqref{eq:barrier-new-annulus},
further implies that
\begin{align*}
\int_{B(z,3R)\setminus B(z,2R)}|\nabla h(x)|^p\,dx &\le\int_D\zeta(x)^p|\nabla h(x)|^p\,dx \le p^p\int_Dh(x)^p|\nabla\zeta(x)|^p\,dx \lesssim R^{d-p}m^p.
\end{align*}
Using this, \eqref{eq:barrier-new-identity}, $0\le\eta\le1$,
$\operatorname{supp}\nabla\eta\subset
\overline{B(z,3R)}\setminus B(z,2R)$, H\"older's inequality,
$|\nabla\eta|\lesssim R^{-1}$,
and \eqref{eq:barrier-new-moment}, we find that
\begin{align*}
\mathcal J_{R^p}(h;0)
&\le\left[\int_{B(z,3R)\setminus B(z,2R)}|\nabla h(x)|^p\,dx
\right]^{\frac{p-1}{p}}
\|\nabla\eta\|_{L^p(D)}
+R^{-p}\int_{B(z,4R)}h(x)^{p-1}\,dx\\
&\lesssim\left(R^{d-p}m^p\right)^{\frac{p-1}{p}}
R^{\frac{d-p}{p}}+R^{d-p}m^{p-1}
\lesssim R^{d-p}m^{p-1}.
\end{align*}
This then completes the proof of Lemma~\ref{lem:local-barrier}.
\end{proof}

\begin{lemma}\label{lem:local-comparison}
Let $p,$ $z,$ $R,$ $E,$ $D$, and $h$ be as in
Lemma~\ref{lem:local-barrier}, and let $f\in W^{1,p}_0(\Omega)$ satisfy
$0\le f\le1$.
Assume that $E\subset F$.
Then, for any $t\in(0,\infty)$, there exists a unique
$v_t\in W^{1,p}_0(\Omega)$ satisfying $0\le v_t\le1$
$p$-quasi-everywhere on $\Omega$ such that
$$\mathcal J_t(v_t;f)=\inf_{v\in W^{1,p}_0(\Omega)}\mathcal J_t(v;f).$$
Moreover, for any $t\in[R^p,\infty)$ and almost every
$x\in\Omega\cap B(z,8R)$,
$v_t(x)\le1-h(x).$
\end{lemma}

\begin{proof}
Applying \cite[Corollary~3.23]{Brezis2011}
to $W^{1,p}_0(\Omega)$ with
$A:=W^{1,p}_0(\Omega)$ and $\varphi:=\mathcal J_t(\,\cdot\,;f)$,
and using the strict convexity of $\mathcal J_t(\,\cdot\,;f)$, we conclude that, for any
$t\in(0,\infty)$, there exists a unique $v_t\in W^{1,p}_0(\Omega)$ satisfying
\begin{align}\label{eq:local-resolvent-minimum}
\mathcal J_t(v_t;f)=\inf_{v\in W^{1,p}_0(\Omega)}\mathcal J_t(v;f).
\end{align}
For any $t\in\mathbb{R}$, let $T(t):=\min\{1,\max\{0,t\}\}$.
By the Sobolev chain rule, $T(0)=0$, and $0\le f\le1$,
we have $T(v_t)\in W^{1,p}_0(\Omega)$ and
\begin{align*}
\mathcal J_t(T(v_t);f)
&=\int_{\{x\in\Omega:\,0<v_t(x)<1\}}|\nabla v_t(x)|^p\,dx
+t^{-1}\int_\Omega|f(x)-T(v_t(x))|^p\,dx
\le\mathcal J_t(v_t;f).
\end{align*}
Using this, \eqref{eq:local-resolvent-minimum}, and the uniqueness
of the minimizer, we conclude that $v_t=T(v_t)$ in
$W^{1,p}_0(\Omega)$ and $0\le v_t\le1$
$p$-quasi-everywhere on $\Omega$.
From \eqref{eq:local-resolvent-minimum} and
$v_t+\varepsilon\phi\in W^{1,p}_0(\Omega)$ for any
$\varepsilon\in\mathbb{R}$ and $\phi\in W^{1,p}_0(\Omega)$,
it follows that, for any such $\phi$,
\begin{align}\label{eq:nonlinear-resolvent}
0 &=\left.\frac1p\frac{d}{d\varepsilon} \mathcal J_t(v_t+\varepsilon\phi;f)\right|_{\varepsilon=0} =\int_\Omega|\nabla v_t|^{p-2}\nabla v_t\cdot\nabla\phi\,dx -t^{-1}\int_\Omega|f-v_t|^{p-2}(f-v_t)\phi\,dx.
\end{align}

For any $t\in[R^p,\infty)$, let
$\phi:=\max\{v_t+h-1,0\}$ on $\mathbb{R}^d$.
Using $v_t\in W^{1,p}_0(\Omega)$,
$h\in W^{1,p}_0(D)$, and $v_t,h\in[0,1]$, we obtain
$\phi\in W^{1,p}_0(\Omega\cap D)$ and
$\phi=0$ $p$-quasi-everywhere on $E$.
By \eqref{eq:barrier-new-stationarity},
$t\ge R^p$,
and $f\le1$, we conclude that
\begin{align*}
-\int_{\Omega\cap D}|\nabla h|^{p-2}\nabla h\cdot\nabla\phi\,dx
&=R^{-p}\int_{\Omega\cap D}h^{p-1}\phi\,dx
\ge t^{-1}\int_{\Omega\cap D}h^{p-1}\phi\,dx\\
&\ge t^{-1}\int_{\Omega\cap D}|f-1+h|^{p-2}(f-1+h)\phi\,dx.
\end{align*}
Using this, the fact that $\nabla\phi=\mathbf1_{\{v_t+h>1\}}(\nabla v_t+\nabla h)$
a.e. in $\Omega\cap D$, and \eqref{eq:nonlinear-resolvent} with the same $\phi$,
we obtain
\begin{align}\label{eq:local-comparison-difference}
0&\le\int_{\{x\in\Omega\cap D:\,v_t(x)+h(x)>1\}}
\left(|\nabla v_t|^{p-2}\nabla v_t+
|\nabla h|^{p-2}\nabla h\right)
\cdot(\nabla v_t+\nabla h)\,dx\notag\\
&\le t^{-1}\int_{\Omega\cap D}
|f-v_t|^{p-2}(f-v_t)\phi\,dx
+\int_{\Omega\cap D}|\nabla h|^{p-2}
\nabla h\cdot\nabla\phi\,dx\notag\\
&\le t^{-1}\int_{\Omega\cap D}
\left[|f-v_t|^{p-2}(f-v_t)-
|f-1+h|^{p-2}(f-1+h)\right]\phi\,dx.
\end{align}
From the definition of $\phi$, we infer that
$f-v_t<f-1+h$ a.e. on
$$\{x\in\Omega\cap D:\phi(x)>0\}=\{x\in\Omega\cap D:\phi(x)\neq0\},$$
which, combined with \eqref{eq:local-comparison-difference},
further implies that
$\phi=0$ a.e. in $\Omega\cap D$.
By the definition of $\phi$, we obtain the desired conclusion, which
completes the proof of Lemma~\ref{lem:local-comparison}.
\end{proof}

\begin{lemma}\label{lem:resolvent-sum}
Let $p\in(1,\infty)$, $\Omega\subset\mathbb R^d$ be open, and
$f\in W^{1,p}_0(\Omega)$ satisfy $0\le f\le1$.
For any $t\in(0,\infty)$, let $v_t$ satisfy
\eqref{eq:local-resolvent-minimum}.
For any $x\in\Omega$, let
$$S_f(x):=\sum_{j\in\mathbb Z}\frac{|f(x)-v_{2^j}(x)|^p}{2^j}.$$
Then
\begin{align}\label{eq:resolvent-sum}
\int_\Omega S_f(x)\,dx
=\sum_{j\in\mathbb Z}\frac{\|f-v_{2^j}\|_{L^p(\Omega)}^p}{2^j}
\le2[f]_{1,p}^p.
\end{align}
\end{lemma}

\begin{proof}
For any $j\in\mathbb{Z}$, let $m_j:=\mathcal J_{2^j}(v_{2^j};f)$.
By \eqref{eq:local-resolvent-minimum} and $\mathcal J_{2^j}(f;f)=[f]_{1,p}^p$,
we have $m_j\le[f]_{1,p}^p$.
Using \eqref{eq:local-resolvent-minimum} with $t:=2^{j+1}$
and $v:=v_{2^j}$ and the definition of $\mathcal J_t(\,\cdot\,;f)$,
we obtain, for any $j\in\mathbb{Z}$,
\begin{align*}
m_{j+1} &\le\mathcal J_{2^{j+1}}(v_{2^j};f) =[v_{2^j}]_{1,p}^p+\frac{\|f-v_{2^j}\|_{L^p(\Omega)}^p}{2^{j+1}} =m_j-\frac{\|f-v_{2^j}\|_{L^p(\Omega)}^p}{2^{j+1}}.
\end{align*}
From this and $m_j\le[f]_{1,p}^p$ for any $j\in\mathbb{Z}$,
we deduce that, for any $a,b\in\mathbb{Z}$ satisfying $a\le b$,
\begin{align*}
\sum_{j=a}^b \frac{\|f-v_{2^j}\|_{L^p(\Omega)}^p}{2^j}
&\le2\sum_{j=a}^b(m_j-m_{j+1})
=2(m_a-m_{b+1})\le2[f]_{1,p}^p.
\end{align*}
Using this and letting $a\to-\infty$ and $b\to\infty$,
we complete the proof of Lemma~\ref{lem:resolvent-sum}.
\end{proof}

\begin{proposition}\label{prop:local-plateau}
Let $p\in(1,d)$, and let $\Omega\subset\mathbb R^d$ be open.
Then there exists a positive constant $C$ such that,
for any $\alpha\in(0,1]$ and $f\in W^{1,p}_0(\Omega)$ with
$0\le f\le1$,
\begin{align*}
\int_{\{x\in\Omega:\,f(x)=1\}}\frac{1}{d_{\alpha,1,p}(x)^p}\,dx
\le C\alpha^{-p'}[f]_{1,p}^p.
\end{align*}
\end{proposition}

\begin{proof}
Let $\alpha\in(0,1]$ and $f\in W^{1,p}_0(\Omega)$ satisfy
$0\le f\le1$.
For any $j\in\mathbb{Z}$ and $z\in\mathbb Q^d$, let
$$E_{j,z}:=F\cap\overline{B(z,2^{\frac{j}{p}})}.$$
Let $S_f$ be the same as in Lemma~\ref{lem:resolvent-sum}.
For any $(j,z)\in\mathbb Z\times\mathbb Q^d$, let $h_{j,z}$
be the minimizer in Lemma~\ref{lem:local-barrier} with
$E:=E_{j,z}$ and $R:=2^{\frac jp}$.
By Lemmas~\ref{lem:local-barrier} and \ref{lem:local-comparison}
and the countability of $\mathbb Z\times\mathbb Q^d$,
we conclude that there exists a Lebesgue-null set $N_0\subset\Omega$
such that, for any $(j,z)\in\mathbb Z\times\mathbb Q^d$ and
$x\in[\Omega\cap B(z,4\cdot2^{\frac jp})]\setminus N_0$,
$$1-v_{2^j}(x)\ge h_{j,z}(x)
\gtrsim\left[\frac{\mathcal J_{2^j}(h_{j,z};0)}
{2^{\frac{j(d-p)}p}}\right]^{\frac1{p-1}}.$$

Let $x\in\Omega\setminus N_0$ satisfy $f(x)=1$ and
$d_{\alpha,1,p}(x)<\infty$,
$r_x\in[d_{\alpha,1,p}(x),2d_{\alpha,1,p}(x))$
be admissible in the definition of $d_{\alpha,1,p}(x)$,
$j_x\in\mathbb Z$ satisfy
$4d_{\alpha,1,p}(x)\le2^{\frac{j_x}p}
<2^{2+\frac1p}d_{\alpha,1,p}(x)$, and
$z\in\mathbb Q^d$ satisfy $|z-x|<d_{\alpha,1,p}(x)$.
From $r_x+|z-x|<3d_{\alpha,1,p}(x)<2^{\frac{j_x}p}$,
we infer that $\overline{F\cap B(x,r_x)}\subset E_{j_x,z}$.
By the Sobolev characterization of $\operatorname{cap}_{1,p}$,
we have
$$\operatorname{cap}_{1,p}(E_{j_x,z})\le[h_{j_x,z}]_{1,p}^p
\le\mathcal J_{2^{j_x}}(h_{j_x,z};0).$$
Using this, Lemmas~\ref{lem:local-barrier} and
\ref{lem:local-comparison}, the admissibility of $r_x$,
the monotonicity of $\operatorname{cap}_{1,p}$, and $r_x\sim2^{\frac{j_x}p}$, we obtain
\begin{align*}
1-v_{2^{j_x}}(x)
&\gtrsim\left[\frac{\mathcal J_{2^{j_x}}(h_{j_x,z};0)}
{2^{\frac{j_x(d-p)}p}}\right]^{\frac1{p-1}}
\ge\left[\frac{\operatorname{cap}_{1,p}(E_{j_x,z})}
{2^{\frac{j_x(d-p)}p}}\right]^{\frac1{p-1}}\\
&\ge\left[\frac{\operatorname{cap}_{1,p}
(\overline{F\cap B(x,r_x)})}
{2^{\frac{j_x(d-p)}p}}\right]^{\frac1{p-1}}
\gtrsim\alpha^{\frac1{p-1}}.
\end{align*}
From this and $f(x)=1$, we deduce that
\begin{align}\label{eq:local-plateau-new-deficit}
S_f(x)
&\ge \frac{|f(x)-v_{2^{j_x}}(x)|^p}{2^{j_x}}
\gtrsim\frac{\alpha^{p'}}{2^{j_x}}
\gtrsim\alpha^{p'}\frac{1}{d_{\alpha,1,p}(x)^p}.
\end{align}
By Lemma~\ref{lem:resolvent-sum}, we have
$\int_\Omega S_f(x)\,dx\le2[f]_{1,p}^p.$
Using this and \eqref{eq:local-plateau-new-deficit},
we obtain
\begin{align*}
\int_{\{x\in\Omega:\,f(x)=1\}}\frac{1}{d_{\alpha,1,p}(x)^p}\,dx
&\lesssim\alpha^{-p'}\int_\Omega S_f(x)\,dx
\le2\alpha^{-p'}[f]_{1,p}^p.
\end{align*}
This completes the proof of Proposition~\ref{prop:local-plateau}.
\end{proof}

Now, we can prove Theorem~\ref{thm:generalp}(i) in the case $s=1$.

\begin{proof}[Proof of Theorem~\ref{thm:generalp}(i) for $s=1$]
Let $\alpha\in(0,1]$.
From Lemma~\ref{lem:distance-measurable}, we infer that
$d\nu(x):=\frac{1}{d_{\alpha,1,p}(x)^p}\,dx$ defines a nonnegative measure
on $\Omega$.
Using Proposition~\ref{prop:local-plateau} and Lemma~\ref{lem:plateau},
we obtain \eqref{eq:generalp-main} in the case $s=1$.
This completes the proof of Theorem~\ref{thm:generalp}(i) in the case $s=1$.
\end{proof}

\section{Proof of Theorem~\ref{thm:generalp}(i)
in the Case \texorpdfstring{$s\ne1$}{s is not equal to 1}}\label{sec:higher}

In this section, we prove Theorem~\ref{thm:generalp}(i) in the case $s\ne1$.
Specifically, in Subsection~\ref{s3.1} we recall the Besov--Sobolev characterization,
compare Euclidean capacities and tree capacities,
and establish a weighted Hardy inequality.
In Subsection~\ref{s3.2}, we estimate the affine approximations on cubes
where the linear part dominates.
Finally, in Subsection~\ref{s3.3} we combine these estimates with a high- and
low-frequency decomposition to prove a dyadic Hardy inequality,
and then select cubes at the scale of the capacitary distance
to complete the proof.

\subsection{Dyadic Preliminaries}\label{s3.1}

Let $\varphi\in C_{\rm c}^\infty(\mathbb{R}^d)$ be real-valued and radial such that
$\operatorname{supp}\varphi\subset B(0,2)$,
$0\le\varphi(\xi)\le1$ for any $\xi\in\mathbb{R}^d$, and
$\varphi(\xi)=1$ for any $\xi\in B(0,1)$.
For any $j\in\mathbb{Z}$ and $f\in\mathcal S'(\mathbb{R}^d)$, let
$$\Delta_j f:=\mathcal F^{-1}\left(\left[\varphi\left(2^{-j-1}\cdot\right)
-\varphi\left(2^{-j}\cdot\right)\right]\widehat f\right),\quad
S_j f:=\mathcal F^{-1}\left(\varphi\left(2^{-j}\cdot\right)\widehat f\right),$$
and $H_j:=I-S_j$, where $I$ is the identity operator.
By the properties of $\varphi$, we have,
for any $j\in\mathbb{Z}$ and $f\in\mathcal S(\mathbb{R}^d)$,
$f=\sum_{k\in\mathbb{Z}}\Delta_k f$ and
$S_j f=\sum_{k\in\mathbb{Z}:\,k<j}\Delta_k f$ in $\mathcal S'(\mathbb{R}^d)$.
Based on the above preparations, we now recall the homogeneous Besov
quasi-seminorm; see, for example,
\cite{Triebel,Triebel2013,YuanSickelYang2010}.

\begin{definition}
Let $\sigma\in\mathbb{R}$ and $p,q\in(0,\infty]$.
For any $f\in\mathcal S'(\mathbb{R}^d)$, the \emph{homogeneous Besov quasi-seminorm}
is defined by
\begin{align*}
\|f\|_{\dot B^\sigma_{p,q}(\mathbb{R}^d)}
:=\left\|\left\{2^{j\sigma}\|\Delta_j f\|_{L^p(\mathbb{R}^d)}\right\}_{j\in\mathbb{Z}}
\right\|_{\ell^q(\mathbb{Z})}.
\end{align*}
\end{definition}

Using the Besov--Sobolev characterization (see, for instance,
\cite[Theorem~6.3.1]{BerghLofstrom}), we obtain, for any
$s\in(0,2)\setminus\{1\}$, $p\in(1,\infty)$, and
$u\in\mathcal S(\mathbb{R}^d)$,
\begin{align}
[u]_{s,p}^p
&\sim\sum_{j\in\mathbb Z}2^{jsp}\|\Delta_j u\|_{L^p(\mathbb{R}^d)}^p.
\label{eq:besovLP}
\end{align}

\begin{lemma}\label{lem:sampling}
Let $p\in(1,\infty)$, $C_0,C_1\in(0,\infty)$, and $\mathscr D$
be a dyadic grid. Then there exists a positive constant $C$ such that,
for any $j\in\mathbb Z$ and $f\in L^p(\mathbb{R}^d)$ satisfying
$\operatorname{supp}\widehat f\subset B(0,C_0 2^j)$,
\begin{align*}
\sum_{\genfrac{}{}{0pt}{}{Q\in\mathscr D}{\ell(Q)=2^{-j}}}
|Q|\,\sup_{x\in C_1Q}|f(x)|^p
\le C\|f\|_{L^p(\mathbb{R}^d)}^p.
\end{align*}
\end{lemma}

\begin{proof}
Let $\eta\in\mathcal S(\mathbb{R}^d)$ satisfy
$\widehat\eta\in C_{\rm c}^\infty(\mathbb{R}^d)$ and $\widehat\eta(\xi)=1$
for any $\xi\in B(0,C_0)$.
For any $x\in\mathbb{R}^d$ and $j\in\mathbb Z$, let $\eta_j(x)=2^{jd}\eta(2^jx)$.
From the Fourier support condition on $f$ and $\eta$, we infer that $f=\eta_j*f$.
Let $M\in(d,\infty)$. For any $x\in\mathbb{R}^d$ and $j\in\mathbb Z$, let
$K_j(x):=2^{jd}(1+2^j|x|)^{-M}$. Since $\eta\in\mathcal S(\mathbb{R}^d)$,
it follows that
$$|\eta_j(w)|\lesssim 2^{jd}\left(1+2^j|w|\right)^{-M}$$
for any $w\in\mathbb{R}^d$.
Using this and the fact that, for any $x\in Q$, $y\in C_1Q$, and $z\in\mathbb R^d$,
$1+2^j|x-z|\lesssim 1+2^j|y-z|$,
we find that,
for any $j\in\mathbb Z$, $Q\in\mathscr D$ having
$\ell(Q)=2^{-j}$, $x\in Q$, and $y\in C_1Q$,
\begin{align*}
|f(y)|\lesssim \int_{\mathbb{R}^d}
\frac{2^{jd}|f(z)|}{(1+2^j|x-z|)^M}\,dz
=(K_j*|f|)(x),
\end{align*}
Raising this estimate to the $p$-th
power and integrating over $x\in Q$, we find that
\begin{align*}
|Q|\sup_{y\in C_1Q}\left|f(y)\right|^p\lesssim
\int_Q\left(K_j*|f|\right)(x)^p\,dx.
\end{align*}
Since $M\in(d,\infty)$, it follows that
$$\|K_j\|_{L^1(\mathbb{R}^d)}=\int_{\mathbb{R}^d}(1+|z|)^{-M}\,dz<\infty$$
for any $j\in\mathbb Z$.
By the disjointness of the cubes of generation $j$ and Young's inequality,
we conclude that
\begin{align*}
\sum_{\genfrac{}{}{0pt}{}{Q\in\mathscr D}{\ell(Q)=2^{-j}}}
|Q|\sup_{y\in C_1Q}|f(y)|^p &\lesssim
\left\|K_j*|f|\right\|_{L^p(\mathbb{R}^d)}^p \le
\left\|K_j\right\|_{L^1(\mathbb{R}^d)}^p\|f\|_{L^p(\mathbb{R}^d)}^p \sim
\|f\|_{L^p(\mathbb{R}^d)}^p.
\end{align*}
This completes the proof of Lemma~\ref{lem:sampling}.
\end{proof}

By \cite[p.~2]{HytonenKairema}, we find that there exist a positive constant
$C_{\mathrm{cov}}$ and $3^d$ adjacent Euclidean dyadic grids $\{\mathscr D_t\}_{t=1}^{3^d}$
such that, for any $x\in\mathbb{R}^d$ and $r\in(0,\infty)$, there exist $t\in[1,3^d]\cap\mathbb Z$
and $Q\in\mathscr D_t$ satisfying
$\overline{B(x,r)}\subset Q^\circ$ and $2r\le\ell(Q)\le C_{\mathrm{cov}}r$,
here and thereafter, the notation $Q^\circ$ denotes the \emph{interior} of $Q$.
Let $\mathscr Q:=\bigcup_{t=1}^{3^d}\mathscr D_t$.

Let $\mathscr D(Q_0)$ denote the half-open dyadic descendants of a cube
$Q_0$.
Let $p\in(1,\infty)$ and $s\in(0,\infty)$.
For any cube $Q_0$, any compact set $E\subset Q_0$, and any
$N\in(j(Q_0),\infty)\cap\mathbb Z$, let
$$\mathscr D_N(Q_0):=\{Q\in\mathscr D(Q_0):j(Q)\le N\}$$
and $U_{\widehat{Q_0}}=0$ and $U_Q=1$ for any $Q\in\mathscr D(Q_0)$ satisfying
$\ell(Q)=2^{-N}$ and $Q\cap E\ne\varnothing$,
and let
\begin{align}\label{eq:finite-treecap}
\mathsf C_N(E;Q_0):=\inf\left[\sum_{Q\in\mathscr D_N(Q_0)}
\ell(Q)^{d-sp}\left|U_Q-U_{\widehat Q}\right|^p\right]
\end{align}
over all families $\{U_Q\}_{Q\in\mathscr D_N(Q_0)}\subset\mathbb{R}$
satisfying these constraints. Recall that $\widehat Q_0$ denotes the parent of $Q_0$.
By the definition of
$\mathsf C_N(E;Q_0)$, we obtain
$\mathsf C_{N+1}(E;Q_0)\le\mathsf C_N(E;Q_0)$
for any
$N\in(j(Q_0),\infty)\cap\mathbb Z$.
For any cube $Q_0$, any compact set $E\subset Q_0$,
and any $p\in(1,\infty)$, we define the \emph{tree capacity} by
\begin{align*}
\operatorname{cap}_{\mathscr D,p}(E;Q_0)
:=\lim_{N\to\infty}\mathsf C_N(E;Q_0).
\end{align*}
For studies of capacities on trees, see
\cite{BenjaminiPeres1992} and \cite[Section~5.1]{ArcozziRochbergSawyerWick2014}.

The following lemma relates the tree capacity to the classical capacity.

\begin{lemma}\label{lem:capcompare}
Let $s\in(0,2)\setminus\{1\}$ and $p\in(1,\frac{d}{s})$.
Then there exists a positive constant $C$ such that, for any dyadic
cube $Q_0$, any compact set $E\subset Q_0$, and any
$N\in(j(Q_0),\infty)\cap\mathbb Z$,
\begin{align}\label{eq:capcompare-sp}
\mathsf C_N(E;Q_0)\ge C^{-1}\operatorname{cap}_{s,p}(E)\quad\text{and}\quad
\operatorname{cap}_{\mathscr D,p}(E;Q_0)\ge C^{-1}\operatorname{cap}_{s,p}(E).
\end{align}
\end{lemma}

\begin{proof}
Let $p\in(1,\infty)$.
Let $\{U_Q\}_{Q\in\mathscr D_N(Q_0)}\subset\mathbb{R}$ satisfy
$U_{\widehat{Q_0}}=0$ and $U_Q=1$ for any $Q\in\mathscr D_N(Q_0)$
with $j(Q)=N$ and $Q\cap E\ne\varnothing$.
Let $\psi\in C_{\rm c}^\infty(\mathbb{R}^d)$ satisfy $0\le\psi\le1$,
$\operatorname{supp}\psi\subset[-2,2]^d$, and $\psi(x)=1$
for any $x\in[-1,1]^d$.
For any $Q\in\mathscr D_N(Q_0)$ and $x\in\mathbb{R}^d$, let
$$\psi_Q(x):=\psi\left(\frac{x-x_Q}{\ell(Q)}\right).$$
Using the smoothness and compact support of $\psi$, we find that
there exists a positive constant $A$ such that
$A^{-1}\ell(Q)^{s-\frac{d}{p}}\psi_Q$ is an $(s,p)$-molecule
in the sense of \cite[(3.1)--(3.2)]{FrazierJawerth1985}
for any $Q\in\mathscr D_N(Q_0)$.
For any $Q\in\mathscr D_N(Q_0)$ and $x\in\mathbb{R}^d$, let
$c_Q:=\max\{U_Q-U_{\widehat Q},0\}$ and
$$v(x):=\sum_{Q\in\mathscr D_N(Q_0)}c_Q\psi_Q(x).$$
From these definitions, we infer that
$0\le c_Q\le|U_Q-U_{\widehat Q}|$ for any $Q\in\mathscr D_N(Q_0)$.
Using $\psi_Q\ge0$, $\psi_Q=1$ on $\overline Q$, and the constraints
on $U_Q$, we find that, for any $x\in E$,
\begin{align*}
v(x)
=\sum_{Q\in\mathscr D_N(Q_0)}c_Q\psi_Q(x)\ge\sum_{\genfrac{}{}{0pt}{}{Q\in\mathscr D_N(Q_0)}{x\in Q}}c_Q
\ge\sum_{\genfrac{}{}{0pt}{}{Q\in\mathscr D_N(Q_0)}{x\in Q}}
\left(U_Q-U_{\widehat Q}\right)=1,
\end{align*}
which, together with $\operatorname{supp}\psi_Q\subset4\overline Q$
and $Q\subset Q_0$, yields
\begin{align}\label{eq:tree-test-admissibility}
v\in C_{\rm c}^\infty(5Q_0)
\quad\text{and}\quad
E\subset\left\{x\in\mathbb R^d:v(x)\ge1\right\}.
\end{align}
By this and $v\in C_{\rm c}^\infty(\mathbb{R}^d)$,
we have, for any $\varepsilon\in(0,\infty)$,
$(1+\varepsilon)v\ge1$ on an open neighborhood of $E$.
In the notation of \cite[Theorem~3.1(b)]{FrazierJawerth1985},
let $\alpha:=s$, $q:=p$,
$m_Q:=A^{-1}\ell(Q)^{s-\frac dp}\psi_Q$, and
$s_Q:=A\ell(Q)^{\frac dp-s}c_Q$ for any
$Q\in\mathscr D_N(Q_0)$.
Since $p\in(1,\infty)$ and $s\in(0,\infty)$, it follows that
the required moment order can be chosen as
$\max\{\lfloor d(\frac1p-1)_+-s\rfloor,-1\}=-1$,
and hence the moment condition is vacuous.
Therefore, from \eqref{eq:besovLP} and
\cite[Theorem~3.1(b)]{FrazierJawerth1985}, we infer that
\begin{align}\label{eq:tree-test-energy}
[v]_{s,p}^p &\sim\|v\|_{\dot B^s_{p,p}(\mathbb R^d)}^p \lesssim\sum_{Q\in\mathscr D_N(Q_0)}\ell(Q)^{d-sp}c_Q^p \le\sum_{Q\in\mathscr D_N(Q_0)} \ell(Q)^{d-sp}\left|U_Q-U_{\widehat Q}\right|^p.
\end{align}
For $p\in(1,\frac ds)$, by the definition of
$\operatorname{cap}_{s,p}$, we have, for any $\varepsilon\in(0,\infty)$,
$$\operatorname{cap}_{s,p}(E)\le(1+\varepsilon)^p[v]_{s,p}^p.$$
Using \eqref{eq:tree-test-energy}, taking the infimum over all
families satisfying the stated constraints, and letting
$\varepsilon\to0^+$, we obtain
$\operatorname{cap}_{s,p}(E)\lesssim\mathsf C_N(E;Q_0)$.
From this, the definition of tree capacity, and letting $N\to\infty$,
it follows that \eqref{eq:capcompare-sp} holds.
This completes the proof of Lemma~\ref{lem:capcompare}.
\end{proof}

In the following lemma, we use the minimizer of
\eqref{eq:finite-treecap} to construct a nonnegative measure.

\begin{lemma}\label{lem:tree-equilibrium}
Let $p\in(1,\infty)$, $s\in(0,\infty)$, $Q_0$ be a dyadic cube,
$E\subset Q_0$ be compact, and
$N\in(j(Q_0),\infty)\cap\mathbb Z$.
Then there exist a real family $\{U_Q\}_{Q\in\mathscr D_N(Q_0)}$
attaining the infimum
of $\mathsf C_N(E;Q_0)$ and satisfying
$0\le U_{\widehat Q}\le U_Q\le1$ for any $Q\in\mathscr D_N(Q_0)$.
Moreover, there exists a nonnegative and finitely supported measure $\mu_N$ on $E$
satisfying, for any $Q\in\mathscr D_N(Q_0)$,
\begin{align}\label{eq:finite-current}
\mu_N(Q)=\ell(Q)^{d-sp}\left(U_Q-U_{\widehat Q}\right)^{p-1}.
\end{align}
\end{lemma}

\begin{proof}
Let $$\mathscr L:=\left\{R\in\mathscr D_N(Q_0):j(R)=N,\ R\cap E\ne\varnothing\right\}$$
and $T(t):=\min\{1,\max\{0,t\}\}$ for any $t\in\mathbb{R}$.
From the definition of $T$, we infer that, for any $t,r\in\mathbb{R}$,
\begin{align}\label{eq:tree-truncation}
T(0)=0,\quad T(1)=1,\quad T(t)\in[0,1],
\quad\text{and}\quad |T(t)-T(r)|\le|t-r|.
\end{align}
Using \eqref{eq:tree-truncation}, the compactness of
$[0,1]^{\#\mathscr D_N(Q_0)}$ with $\#\mathscr D_N(Q_0)$ denoting the cardinality
of $\mathscr D_N(Q_0)$,
and the continuity of the sum in
\eqref{eq:finite-treecap}, we find that there exists a family
$\{V_Q\}_{Q\in\mathscr D_N(Q_0)}\subset[0,1]$ with
$V_{\widehat{Q_0}}=0$ and $V_R=1$ for any $R\in\mathscr L$ such that
\begin{align}\label{eq:tree-initial-minimum}
\sum_{Q\in\mathscr D_N(Q_0)}\ell(Q)^{d-sp}\left|V_Q-V_{\widehat Q}\right|^p
=\mathsf C_N(E;Q_0).
\end{align}
Let $U_{\widehat{Q_0}}:=0$ and
$U_Q:=\max_{\genfrac{}{}{0pt}{}{R\in\mathscr D_N(Q_0)}{Q\subset R}}V_R$
for any $Q\in\mathscr D_N(Q_0)$.
From these definitions, we deduce that $U_R=1$ for any $R\in\mathscr L$
and, for any $Q\in\mathscr D_N(Q_0)$,
$0\le U_{\widehat Q}\le U_Q\le1$,
$V_{\widehat Q}\le U_{\widehat Q}$,
and $$0\le U_Q-U_{\widehat Q}
=\max\left\{V_Q-U_{\widehat Q},0\right\}
\le\left|V_Q-V_{\widehat Q}\right|.$$
By this, \eqref{eq:tree-initial-minimum},
and the definition of $\mathsf C_N(E;Q_0)$, we have
\begin{align}\label{eq:tree-minimum}
\mathsf C_N(E;Q_0) &\le\sum_{Q\in\mathscr D_N(Q_0)}\ell(Q)^{d-sp}\left|U_Q-U_{\widehat Q}\right|^p \le\sum_{Q\in\mathscr D_N(Q_0)}\ell(Q)^{d-sp}\left|V_Q-V_{\widehat Q}\right|^p =\mathsf C_N(E;Q_0).
\end{align}
This shows that $\{U_Q\}_{Q\in\mathscr D_N(Q_0)}$
attain the infimum
of $\mathsf C_N(E;Q_0)$.

For any $Q\in\mathscr D_N(Q_0)$, let $\delta_Q:=U_Q-U_{\widehat Q}$.
By replacing only $U_Q$ with $t\in\mathbb{R}$ and using
\eqref{eq:tree-minimum}, we obtain, for any
$Q\in\mathscr D_N(Q_0)$ with $j(Q)<N$,
\begin{align}
0&=\left.\frac{d}{dt}\left[
\ell(Q)^{d-sp}\left|t-U_{\widehat Q}\right|^p+
\sum_{\genfrac{}{}{0pt}{}{R\in\mathscr D_N(Q_0)}{\widehat R=Q}}
\ell(R)^{d-sp}|U_R-t|^p\right]\right|_{t=U_Q}\notag\\
&=p \ell(Q)^{d-sp}\delta_Q^{p-1}
-p\sum_{\genfrac{}{}{0pt}{}{R\in\mathscr D_N(Q_0)}{\widehat R=Q}}
\ell(R)^{d-sp}\delta_R^{p-1}.
\label{eq:tree-first-order}
\end{align}
From \eqref{eq:tree-minimum} and a similar argument, we infer that
$p \ell(R)^{d-sp}\delta_R^{p-1}=0$ and hence $\delta_R=0$ for any
$R\in\mathscr D_N(Q_0)$ with $j(R)=N$ and $R\cap E=\varnothing$.
For any $R\in\mathscr L$, let $z_R\in E\cap R$.
For any Borel set $A\subset\mathbb{R}^d$, let
\begin{align}\label{eq:tree-measure}
\mu_N(A):=\sum_{R\in\mathscr L}\ell(R)^{d-sp}\delta_R^{p-1}\mathbf1_{A}(z_R),
\end{align}
which depends on the choice of
$\{z_R\}_{R\in\mathscr L}$.
Since $\mathscr L$ is finite, $\ell(R)^{d-sp}\delta_R^{p-1}\ge0$, and $z_R\in E$ exists
for any $R\in\mathscr L$, it follows that $\mu_N\ge0$ and
$\operatorname{supp}\mu_N\subset\{z_R:R\in\mathscr L\}\subset E$.
Using \eqref{eq:tree-measure}, \eqref{eq:tree-first-order}, and
$\delta_R=0$ for any $R\in\mathscr D_N(Q_0)$ satisfying
$j(R)=N$ and $R\cap E=\varnothing$,
we obtain, for any $Q\in\mathscr D_N(Q_0)$,
\begin{align*}
\mu_N(Q) &=\sum_{R\in\mathscr L}\ell(R)^{d-sp}\delta_R^{p-1}\mathbf1_{Q}(z_R) =\sum_{\genfrac{}{}{0pt}{}{R\in\mathscr L}{R\subset Q}}\ell(R)^{d-sp}\delta_R^{p-1}\\ &=\ell(Q)^{d-sp}\delta_Q^{p-1}=\ell(Q)^{d-sp}\left(U_Q-U_{\widehat Q}\right)^{p-1}.
\end{align*}
This completes the proof of Lemma~\ref{lem:tree-equilibrium}.
\end{proof}

The following lemma is a weighted Hardy inequality.

\begin{lemma}\label{lem:rayhardy}
Let $p\in(1,\infty)$, $\beta\in(0,1]$, $M\in\mathbb N$, and
$\{h_n\}_{n\in[0,M]\cap\mathbb Z}$ and
$\{z_n\}_{n\in[0,M]\cap\mathbb Z}$ be real sequences satisfying
$h_0>h_1>\cdots>h_M=0$ and $z_M=0$.
For any $n\in[1,M]\cap\mathbb Z$, let $\delta_n:=h_{n-1}-h_n$,
and let $G\subset[1,M]\cap\mathbb Z$ satisfy,
for any $n\in G$,
\begin{align}\label{eq:densitydrop}
\delta_n\ge\beta h_{n-1}.
\end{align}
Then there exists a positive constant $C$, independent of $\beta$,
$M$, $G$, $\{h_n\}_{n\in[0,M]\cap\mathbb Z}$, and
$\{z_n\}_{n\in[0,M]\cap\mathbb Z}$, such that
\begin{align*}
\sum_{n\in G}\frac{|z_n|^p}{\delta_n^{p-1}}
\le C\beta^{-p}\sum_{n=1}^M\frac{|z_n-z_{n-1}|^p}{\delta_n^{p-1}}.
\end{align*}
\end{lemma}

\begin{proof}
From \eqref{eq:densitydrop}, we infer that, for any $n\in G$,
\begin{align}\label{eq:ray-density-comparison}
\delta_n^{1-p}\le\beta^{1-p}h_{n-1}^{1-p}
\le\beta^{-p}\delta_n h_{n-1}^{-p}.
\end{align}
For any $t\in[0,h_0]$, let
\begin{align*}
f(t):=\sum_{n=1}^M\frac{|z_n-z_{n-1}|}{\delta_n}
\mathbf1_{(h_n,h_{n-1})}(t)
\quad\text{and}\quad Z(t):=\int_0^t f(r)\,dr.
\end{align*}
Using these definitions and the fundamental theorem of calculus,
we obtain $f\ge0$, $Z\in W^{1,p}(0,h_0)$, $Z(0)=0$,
and $Z'=f$ almost everywhere on $(0,h_0)$.
From $z_M=0$ and the triangle inequality, we deduce that, for any
$n\in[1,M]\cap\mathbb Z$ and $t\in[h_n,h_{n-1}]$,
\begin{align}\label{eq:ray-endpoint}
|z_n|\le\sum_{k=n+1}^M|z_k-z_{k-1}|=Z(h_n)\le Z(t).
\end{align}
By the one-dimensional Hardy inequality
\cite[Theorem 327]{HardyLittlewoodPolya}
and $Z(0)=0$, we have
\begin{align}\label{eq:ray-core-bound}
\int_0^{h_0}\frac{|Z(t)|^p}{t^p}\,dt
\lesssim\int_0^{h_0}|f(t)|^p\,dt
=\sum_{n=1}^M\frac{|z_n-z_{n-1}|^p}{\delta_n^{p-1}}.
\end{align}
Using \eqref{eq:ray-density-comparison}, \eqref{eq:ray-endpoint}, and
\eqref{eq:ray-core-bound}, we conclude that
\begin{align*}
\sum_{n\in G}\frac{|z_n|^p}{\delta_n^{p-1}} &\le\beta^{-p}\sum_{n\in G}\int_{h_n}^{h_{n-1}} \frac{|Z(t)|^p}{t^p}\,dt \le\beta^{-p}\int_0^{h_0}\frac{|Z(t)|^p}{t^p}\,dt \lesssim\beta^{-p}\sum_{n=1}^M \frac{|z_n-z_{n-1}|^p}{\delta_n^{p-1}}.
\end{align*}
This completes the proof of Lemma~\ref{lem:rayhardy}.
\end{proof}

\subsection{An Affine Estimate}\label{s3.2}

For any dyadic grid $\mathscr D$, any $Q\in\mathscr D$,
any $u\in\mathcal S(\mathbb{R}^d)$, and any $x\in\mathbb{R}^d$, let
$a_Q:=S_{j(Q)}u(x_Q)$, $b_Q:=\nabla S_{j(Q)}u(x_Q)$, and
$P_Q(x):=a_Q+b_Q\cdot(x-x_Q)$.

\begin{lemma}\label{lem:affine-gradient}
Let $p\in(1,\infty)$,
$s\in(0,1+\frac1p)\setminus\{1\}$,
and $\mathscr D$ be a dyadic grid.
Then there exists a positive constant $C$ such that, for any real-valued
$u\in\mathcal S(\mathbb{R}^d)$,
\begin{align}\label{eq:dangerous-final}
\sum_{\genfrac{}{}{0pt}{}{Q\in\mathscr D}{|a_Q|\le 4\sqrt d\,\ell(Q)|b_Q|}}
\ell(Q)^{d+p-sp}|b_Q|^p\le C[u]_{s,p}^p.
\end{align}
\end{lemma}

\begin{proof}
Let $$\mathscr B:=\left\{Q\in\mathscr D:|a_Q|\le 4\sqrt d\,\ell(Q)|b_Q|\right\}.$$
By $S_ju=\sum_{k\in\mathbb Z:\,k<j}\Delta_ku$ and Bernstein's
inequality (see, for instance,
\cite[Lemma~2.1]{BahouriCheminDanchin}), we have, for any $j\in\mathbb Z$,
\begin{align*}
2^{-j(2-s)}\left\|\nabla^2S_ju\right\|_{L^p(\mathbb{R}^d)}
\lesssim\sum_{\genfrac{}{}{0pt}{}{k\in\mathbb Z}{k<j}}
2^{-(j-k)(2-s)}2^{ks}\|\Delta_ku\|_{L^p(\mathbb{R}^d)}.
\end{align*}
From this, $2-s>0$, Young's inequality on $\ell^p(\mathbb Z)$,
and \eqref{eq:besovLP}, we deduce that
\begin{align}\label{eq:hessian-sum}
\sum_{j\in\mathbb Z}2^{-jp(2-s)}
\left\|\nabla^2S_ju\right\|_{L^p(\mathbb{R}^d)}^p
\lesssim\sum_{k\in\mathbb Z}2^{ksp}
\|\Delta_ku\|_{L^p(\mathbb{R}^d)}^p\sim[u]_{s,p}^p.
\end{align}
For any $R\in\mathscr D$, let
\begin{align}\label{eq:jetdefects}
B_{\widehat R,R}:=b_R-b_{\widehat R}
\quad\text{and}\quad
A_{\widehat R,R}:=a_R-a_{\widehat R}
-b_{\widehat R}\cdot(x_R-x_{\widehat R}).
\end{align}
Using $S_{j(R)}-S_{j(R)-1}=\Delta_{j(R)-1}$ and Taylor's theorem,
we obtain, for any $R\in\mathscr D$,
\begin{align}\label{1941}
\left|A_{\widehat R,R}\right|+\ell\left(\widehat R\right)\left|B_{\widehat R,R}\right|
&\lesssim\sup_{x\in\widehat R}|\Delta_{j(R)-1}u(x)|
+\ell\left(\widehat R\right)\sup_{x\in\widehat R}|\nabla\Delta_{j(R)-1}u(x)|\nonumber\\*
&\quad+\ell\left(\widehat R\right)^2
\sup_{x\in\widehat R}\left|\nabla^2S_{j(R)-1}u(x)\right|.
\end{align}
From the definitions of $\Delta_j$ and $S_j$, we infer that the Fourier supports of
$\Delta_ju$, $\nabla\Delta_ju$, and $\nabla^2S_ju$ are contained in
$B(0,2^{j+2})$ for any $j\in\mathbb Z$.
By this, \eqref{1941}, Lemma~\ref{lem:sampling}, Bernstein's inequality,
\eqref{eq:besovLP}, and \eqref{eq:hessian-sum}, we obtain
\begin{align}\label{eq:jetsquare}
&\sum_{R\in\mathscr D}
\left[\ell\left(\widehat R\right)^{d-sp}\left|A_{\widehat R,R}\right|^p
+\ell\left(\widehat R\right)^{d+p-sp}\left|B_{\widehat R,R}\right|^p\right]\notag\\
&\quad\lesssim\sum_{j\in\mathbb Z}2^{jsp}
\left[\|\Delta_ju\|_{L^p(\mathbb{R}^d)}^p
+2^{-jp}\|\nabla\Delta_ju\|_{L^p(\mathbb{R}^d)}^p
+2^{-2jp}\left\|\nabla^2S_ju\right\|_{L^p(\mathbb{R}^d)}^p\right]\notag\\
&\quad\lesssim\sum_{j\in\mathbb Z}2^{jsp}
\left[\|\Delta_ju\|_{L^p(\mathbb{R}^d)}^p
+2^{-2jp}\left\|\nabla^2S_ju\right\|_{L^p(\mathbb{R}^d)}^p\right]
\lesssim[u]_{s,p}^p.
\end{align}

For any $\varepsilon\in(0,\frac14]$, let $\mathscr G_\varepsilon$ consist of the pairs
$(\widehat R,R)\in\mathscr B\times\mathscr B$ satisfying
\begin{align}\label{eq:goodedge}
\left|B_{\widehat R,R}\right|\le\varepsilon|b_{\widehat R}|
\quad\text{and}\quad
\left|A_{\widehat R,R}\right|\le\varepsilon\ell\left(\widehat R\right)|b_{\widehat R}|.
\end{align}
For any $T\in\mathscr B$, let $\mathscr C(T)$ consist of $T$ and all
$Q\in\mathscr B$ satisfying $Q\subsetneqq T$ and
$(\widehat R,R)\in\mathscr G_\varepsilon$ for any $R\in\mathscr D$
satisfying $Q\subset R\subsetneqq T$.
Let $L\in\mathbb N$ be sufficiently large
which will be determined below, let
$\varepsilon\in(0,\min\{\frac1{4L},2^{-L}\}]$, and let
$T\in\mathscr B$ satisfy $b_T\ne0$.
Using \eqref{eq:goodedge} successively, we find that, for any
$k\in[0,L]\cap\mathbb Z$ and $Q\in\mathscr C(T)$ satisfying
$j(Q)=j(T)+k$,
\begin{align}\label{eq:corona-slope-comparison}
\frac34|b_T|\le(1-\varepsilon)^k|b_T|
\le|b_Q|\le(1+\varepsilon)^k|b_T|\le2|b_T|.
\end{align}
For any $S\in\mathscr C(T)$ satisfying $j(S)=j(T)+L$ and
$k\in[0,L]\cap\mathbb Z$, let $R_k\in\mathscr D$ satisfy
$S\subset R_k\subset T$ and $j(R_k)=j(T)+k$.
Let $y:=x_S-\frac{a_S}{|b_S|^2}b_S$ and
$H_T:=\{x\in\mathbb{R}^d:P_T(x)=0\}$.
By $S\in\mathscr B$ and the definition of $P_S$, we have
$P_S(y)=0$ and $|y-x_S|\le 4\sqrt d\,\ell(S)$.
For any
$k\in[1,L]\cap\mathbb Z$, from $S\subset R_k$, we infer that
$|y-x_{R_k}|\le5\sqrt d\,\ell(\widehat{R_k})$.
Using this, \eqref{eq:jetdefects}, \eqref{eq:goodedge},
and \eqref{eq:corona-slope-comparison}, we obtain
\begin{align*}
|P_T(y)| &\le\sum_{k=1}^{L}\left(\left|A_{\widehat{R_k},R_k}\right|
+\left|B_{\widehat{R_k},R_k}\right|\,|y-x_{R_k}|\right) \lesssim\varepsilon|b_T|\sum_{k=1}^{L}\ell\left(\widehat{R_k}\right) \le2\varepsilon\ell(T)|b_T|.
\end{align*}
By this, the distance formula for $H_T$, and $\varepsilon\le2^{-L}$,
we conclude that
\begin{align*}
\operatorname{dist}(x_S,H_T)
\le|x_S-y|+\frac{|P_T(y)|}{|b_T|}
\lesssim\left(2^{-L}+\varepsilon\right)\ell(T)\lesssim2^{-L}\ell(T).
\end{align*}
Thus, the centers of cubes $\{S\in\mathscr C(T):j(S)=j(T)+L\}$
lie in a slab around $H_T$
whose thickness is comparable to the descendant edge length
$2^{-L}\ell(T)$. From this and the disjointness of the
cubes $\{S: S\in\mathscr C(T),j(S)=j(T)+L\}$,
we deduce that
$$
\#\left\{S\in\mathscr C(T):j(S)=j(T)+L\right\}
\lesssim2^{L(d-1)}.
$$
Using this, \eqref{eq:corona-slope-comparison}, and
$\ell(S)=2^{-L}\ell(T)$, we obtain
\begin{align}\label{eq:corona-block-bound}
\sum_{\genfrac{}{}{0pt}{}{S\in\mathscr C(T)}{j(S)=j(T)+L}}
\ell(S)^{d+p-sp}|b_S|^p
&\lesssim2^{L(d-1)}2^{-L(d+p-sp)}
\ell(T)^{d+p-sp}|b_T|^p\notag\\
&=2^{-L[1-p(s-1)]}\ell(T)^{d+p-sp}|b_T|^p
\end{align}
uniformly in $L$ and $\varepsilon$. This calculation isolates the
origin of the restriction $s<1+\frac1p$, which makes the corona
estimate a strict contraction
for sufficiently large $L$. By \eqref{eq:goodedge}, we conclude that, for any $Q\in\mathscr C(T)$,
if $b_T=0$, then $b_Q=0$.
From \eqref{eq:corona-block-bound}, with an implicit constant independent
of $L$ and $\varepsilon$, and $1-p(s-1)>0$, we infer that there exists
a sufficiently large $L\in\mathbb N$ such that, for any
$\varepsilon\in(0,\min\{\frac1{4L},2^{-L}\}]$ and $T\in\mathscr B$,
\begin{align*}
\sum_{\genfrac{}{}{0pt}{}{S\in\mathscr C(T)}{j(S)=j(T)+L}}
\ell(S)^{d+p-sp}|b_S|^p
\le\frac12\ell(T)^{d+p-sp}|b_T|^p.
\end{align*}
Using this successively, we obtain, for any $T\in\mathscr B$ and
$k\in\mathbb Z_+$,
\begin{align}\label{eq:blockcontract}
\sum_{\genfrac{}{}{0pt}{}{Q\in\mathscr C(T)}{j(Q)=j(T)+kL}}
\ell(Q)^{d+p-sp}|b_Q|^p
\le2^{-k}\ell(T)^{d+p-sp}|b_T|^p.
\end{align}
From \eqref{eq:corona-slope-comparison} and the fact that
$\#\{Q\in\mathscr C(T):j(Q)=j(T)+r\}\le2^{rd}$,
it follows that, for any $T\in\mathscr B$ and
$r\in[0,L-1]\cap\mathbb Z$,
\begin{align*}
\sum_{\genfrac{}{}{0pt}{}{Q\in\mathscr C(T)}{j(Q)=j(T)+r}}
\ell(Q)^{d+p-sp}|b_Q|^p
\le2^p2^{rp(s-1)}\ell(T)^{d+p-sp}|b_T|^p.
\end{align*}
By this, \eqref{eq:blockcontract}, and the
summation over $r\in[0,L-1]\cap\mathbb Z$ and
$k\in\mathbb Z_+$ with $L$ fixed, we conclude that,
for any $T\in\mathscr B$,
\begin{align}\label{eq:coronapacking}
\sum_{Q\in\mathscr C(T)}\ell(Q)^{d+p-sp}|b_Q|^p
\lesssim\sum_{k=0}^{\infty}2^{-k}\ell(T)^{d+p-sp}|b_T|^p
=2\ell(T)^{d+p-sp}|b_T|^p.
\end{align}

Let $R\in\mathscr B$.
Using \eqref{eq:jetdefects}, the definition of $\mathscr B$, and
$|x_R-x_{\widehat R}|\le\frac{\sqrt d}{4}\ell(\widehat R)$,
we find that, if $\widehat R\notin\mathscr B$ and
$|B_{\widehat R,R}|<\frac{|b_R|}{4}$, then
\begin{align*}
\frac34|b_R|\le|b_{\widehat R}|\le\frac54|b_R|
\quad\text{and}\quad
\left|A_{\widehat R,R}\right|\ge\frac{\sqrt d}{2}\ell\left(\widehat R\right)|b_R|.
\end{align*}
From this and the case $|B_{\widehat R,R}|\ge\frac{|b_R|}{4}$,
we infer that, if $\widehat R\notin\mathscr B$, then
\begin{align}\label{eq:birthpoint}
\ell(R)^{d+p-sp}|b_R|^p
\lesssim\ell\left(\widehat R\right)^{d-sp}\left|A_{\widehat R,R}\right|^p
+\ell\left(\widehat R\right)^{d+p-sp}\left|B_{\widehat R,R}\right|^p.
\end{align}
By $|b_R|\le|b_{\widehat R}|+|B_{\widehat R,R}|$ and \eqref{eq:goodedge},
we find that, if $\widehat R\in\mathscr B$ and $(\widehat R,R)\notin\mathscr G_\varepsilon$,
then either
$|b_R|\le(1+\varepsilon^{-1})|B_{\widehat R,R}|$ or
$|b_R|\le(1+\varepsilon)\varepsilon^{-1}\ell(\widehat R)^{-1}|A_{\widehat R,R}|$.
Using this, we conclude that \eqref{eq:birthpoint} also holds in this case.
Thus, if $\widehat R\in\mathscr B$ and $(\widehat R,R)\notin\mathscr G_\varepsilon$, or $\widehat R\notin\mathscr B$,
then \eqref{eq:birthpoint} holds.

For any $m,n\in\mathbb Z$ satisfying $m\le n$, let
$$\mathscr B_{m,n}:=\left\{Q\in\mathscr B:j(Q)\in[m,n]\right\}\quad\text{and}
\quad\mathscr T_{m,n}:=\left\{T\in\mathscr B_{m,n}:j(T)=m\text{ or }
\left(\widehat T,T\right)\notin\mathscr G_\varepsilon\right\}.$$
From the definition of $\mathscr C(T)$, we infer that
$$\bigcup_{T\in\mathscr T_{m,n}}[\mathscr C(T)
\cap\mathscr B_{m,n}]=\mathscr B_{m,n}$$
and $\mathscr C(T)\cap\mathscr C(T')\cap\mathscr B_{m,n}=\varnothing$
for any distinct $T,T'\in\mathscr T_{m,n}$.
Using this, \eqref{eq:coronapacking}, \eqref{eq:birthpoint},
and \eqref{eq:jetsquare}, we obtain
\begin{align}\label{eq:affine-truncation}
\sum_{Q\in\mathscr B_{m,n}}\ell(Q)^{d+p-sp}|b_Q|^p
\lesssim[u]_{s,p}^p
+\sum_{\genfrac{}{}{0pt}{}{Q\in\mathscr D}{j(Q)=m}}
\ell(Q)^{d+p-sp}|b_Q|^p.
\end{align}
Let $K\in\mathcal S(\mathbb{R}^d)$ be the convolution kernel of $S_0$,
i.e., for any $m\in\mathbb Z$,
$S_mu=[2^{md}K(2^m\cdot)]*u$.
From Young's inequality and a change of variables, we deduce that,
for any $m\in\mathbb Z$,
\begin{align*}
\|\nabla S_mu\|_{L^p(\mathbb{R}^d)} &\le\left\|2^{m(d+1)}(\nabla K)(2^m\cdot)\right\|_{L^p(\mathbb{R}^d)} \|u\|_{L^1(\mathbb{R}^d)} =2^{m[1+d(1-\frac1p)]}\|\nabla K\|_{L^p(\mathbb{R}^d)} \|u\|_{L^1(\mathbb{R}^d)}.
\end{align*}
By this and Lemma~\ref{lem:sampling}, we obtain,
for any $m\in\mathbb Z$,
\begin{align*}
\sum_{\genfrac{}{}{0pt}{}{Q\in\mathscr D}{j(Q)=m}}
\ell(Q)^{d+p-sp}|b_Q|^p
\lesssim2^{mp(s-1)}\|\nabla S_mu\|_{L^p(\mathbb{R}^d)}^p\lesssim2^{m[sp+d(p-1)]}\|u\|_{L^1(\mathbb{R}^d)}^p.
\end{align*}
Using this, $sp+d(p-1)>0$, \eqref{eq:affine-truncation}, and the
monotone convergence theorem, and letting $n\to\infty$ and then
$m\to-\infty$, we conclude that \eqref{eq:dangerous-final} holds.
This completes the proof of Lemma~\ref{lem:affine-gradient}.
\end{proof}

\subsection{Proof of Theorem~\ref{thm:generalp}(i)}\label{s3.3}

The following estimate depends on the finite-tree energy condition,
and will be applied to the cubes selected by the capacitary distance.

\begin{proposition}\label{prop:dyadic-hardy-energy}
Let $p\in(1,\infty)$, $s\in(0,1+\frac1p)\setminus\{1\}$,
$\alpha\in(0,1]$, $\Omega\subset\mathbb R^d$ be open,
$F:=\mathbb R^d\setminus\Omega$, and $\mathscr D$ be a dyadic grid.
For any $\mathscr S\subset\mathscr D$, let
$\{A_Q\}_{Q\in\mathscr S}$ be pairwise disjoint measurable sets
with $A_Q\subset\Omega\cap Q$.
Assume that, for any $Q\in\mathscr S$, there exists a compact set
$E_Q\subset F\cap Q$ such that, for any
$N\in(j(Q),\infty)\cap\mathbb Z$,
\begin{align}\label{eq:dyadic-tree-thickness}
\mathsf C_N(E_Q;Q)\ge\alpha\ell(Q)^{d-sp}.
\end{align}
Then there exists a positive constant $C$, depending only on
$d$, $s$, and $p$, such that, for any $u\in C_{\rm c}^\infty(\Omega)$,
\begin{align}\label{eq:dyadic-hardy-energy}
\sum_{Q\in\mathscr S}\ell(Q)^{-sp}\int_{A_Q}|u(x)|^p\,dx
\le C\alpha^{-p'}[u]_{s,p}^p.
\end{align}
\end{proposition}

\begin{proof}
Let $u\in C_{\rm c}^\infty(\Omega)$.
From $I=H_j+S_j$ for any $j\in\mathbb Z$, we infer that,
for any $Q\in\mathscr D$ and $x\in Q$,
\begin{align}\label{eq:higherdecomp}
u(x)=H_{j(Q)}u(x)+S_{j(Q)}u(x).
\end{align}
By this and the triangle inequality, we conclude that it suffices
to prove that
\begin{align*}
\sum_{Q\in\mathscr S}\ell(Q)^{-sp}
\int_{A_Q}\left[|H_{j(Q)}u(x)|^p+|S_{j(Q)}u(x)|^p\right]\,dx
\lesssim\alpha^{-p'}[u]_{s,p}^p.
\end{align*}

We first estimate $H_{j(Q)}u$.
Using $H_ju=\sum_{k\in\mathbb Z:\,k\ge j}\Delta_ku$,
Young's inequality on $\ell^p(\mathbb Z)$, and
\eqref{eq:besovLP}, we obtain
\begin{align*}
\sum_{j\in\mathbb Z}2^{jsp}\|H_ju\|_{L^p(\mathbb{R}^d)}^p &\le\sum_{j\in\mathbb Z}\left[ \sum_{\genfrac{}{}{0pt}{}{k\in\mathbb Z}{k\ge j}} 2^{-(k-j)s}2^{ks}\|\Delta_ku\|_{L^p(\mathbb{R}^d)}\right]^p \lesssim\sum_{k\in\mathbb Z}2^{ksp} \|\Delta_ku\|_{L^p(\mathbb{R}^d)}^p\sim[u]_{s,p}^p.
\end{align*}
From this and the disjointness of the sets $A_Q$, we deduce that
\begin{align}\label{eq:higher-high-selected}
\sum_{Q\in\mathscr S}\ell(Q)^{-sp}
\int_{A_Q}|H_{j(Q)}u(x)|^p\,dx
&=\sum_{j\in\mathbb Z}2^{jsp}
\sum_{\genfrac{}{}{0pt}{}{Q\in\mathscr S}{j(Q)=j}}
\int_{A_Q}|H_ju(x)|^p\,dx\notag\\
&\le\sum_{j\in\mathbb Z}2^{jsp}\|H_ju\|_{L^p(\mathbb{R}^d)}^p
\lesssim[u]_{s,p}^p.
\end{align}

Next, we estimate $S_{j(Q)}u$.
For any $Q\in\mathscr D$ and $x\in Q$, let
$R_Q(x):=S_{j(Q)}u(x)-P_Q(x)$ and hence
\begin{align}\label{eq:lowpass-affine-decomp}
S_{j(Q)}u(x)=R_Q(x)+P_Q(x).
\end{align}
From the definition of $S_j$, we infer that
$\operatorname{supp}\widehat{\nabla^2S_ju}\subset B(0,2^{j+1})$
for any $j\in\mathbb Z$.
By Taylor's theorem, we have, for any $Q\in\mathscr D$ and $x\in Q$,
$$|R_Q(x)|\lesssim\ell(Q)^2\sup_{y\in Q}\left|\nabla^2S_{j(Q)}u(y)\right|.$$
Using this, $A_Q\subset Q$, Lemma~\ref{lem:sampling},
and \eqref{eq:hessian-sum}, we obtain
\begin{align}\label{eq:curvature}
\sum_{Q\in\mathscr S}\ell(Q)^{-sp}\int_{A_Q}|R_Q(x)|^p\,dx
&\le\sum_{Q\in\mathscr D}\ell(Q)^{d-sp}
\sup_{x\in Q}|R_Q(x)|^p\notag\\
&\lesssim\sum_{j\in\mathbb Z}2^{-jp(2-s)}
\left\|\nabla^2S_ju\right\|_{L^p(\mathbb{R}^d)}^p\lesssim[u]_{s,p}^p.
\end{align}
By this and \eqref{eq:lowpass-affine-decomp}, we conclude that
it suffices to estimate the contribution of $P_Q$.
Let $$\mathscr S_a:=\left\{Q\in\mathscr S:|a_Q|>4\sqrt d\,\ell(Q)|b_Q|\right\}$$ and
$\mathscr S_b:=\mathscr S\setminus\mathscr S_a$.
For any $Q\in\mathscr S_b$ and $x\in Q$, from the definitions of
$P_Q$ and $\mathscr S_b$ and $|x-x_Q|\le\frac{\sqrt d}{2}\ell(Q)$,
we infer that $|P_Q(x)|\lesssim\ell(Q)|b_Q|$.
By this, $A_Q\subset Q$, and Lemma~\ref{lem:affine-gradient}, we obtain
\begin{align}\label{eq:dangerPfinal}
\sum_{Q\in\mathscr S_b}\ell(Q)^{-sp}
\int_{A_Q}|P_Q(x)|^p\,dx
\lesssim\sum_{Q\in\mathscr S_b}\ell(Q)^{d+p-sp}|b_Q|^p
\lesssim[u]_{s,p}^p.
\end{align}

Now, we estimate the sum over $\mathscr S_a$.
Using the definitions of $P_Q$ and $\mathscr S_a$, we find that,
for any $Q\in\mathscr S_a$ and $x\in Q$,
\begin{align}\label{eq:safeP}
|P_Q(x)|\sim|a_Q|.
\end{align}
From this and $P_Q=S_{j(Q)}u-R_Q$, we deduce that, for any
$Q\in\mathscr S_a$ and $z\in F\cap Q$,
\begin{align}\label{eq:safe-constant-pointwise}
|a_Q|^p\lesssim|S_{j(Q)}u(z)|^p+|R_Q(z)|^p.
\end{align}
Let $Q_0\in\mathscr D$ and
$\mathscr S_0\subset\mathscr S_a\cap\mathscr D(Q_0)$ be finite and nonempty.
Let $E:=\bigcup_{Q\in\mathscr S_0}E_Q$ and
$N\in(\max_{Q\in\mathscr S_0}j(Q),\infty)\cap\mathbb Z$.
By Lemma~\ref{lem:tree-equilibrium},
we obtain a minimizing family $\{U_Q\}_{Q\in\mathscr D_N(Q_0)}$
and a nonnegative and finitely supported measure $\mu_N$ on $E$.
For any $Q\in\mathscr D_N(Q_0)$, let
$\delta_Q:=U_Q-U_{\widehat Q}$.
Assume that $Q\in\mathscr D_N(Q_0)$ satisfies $j(Q)<N$ and $U_Q=1$.
From Lemma~\ref{lem:tree-equilibrium}, we deduce that
$1=U_Q\le U_R\le1$ and $\delta_R=0$ for any
$R\in\mathscr D_N(Q_0)$ satisfying $\widehat R=Q$.
Using this and \eqref{eq:tree-first-order}, we obtain
\begin{align*}
\ell(Q)^{d-sp}\delta_Q^{p-1}
=\sum_{\genfrac{}{}{0pt}{}{R\in\mathscr D_N(Q_0)}{\widehat R=Q}}
\ell(R)^{d-sp}\delta_R^{p-1}=0.
\end{align*}
From this, $\ell(Q)^{d-sp}>0$, and $p\in(1,\infty)$, we deduce that
$\delta_Q=0$ and $U_{\widehat Q}=U_Q-\delta_Q=1$.
By iteration, we conclude that $U_{\widehat{Q_0}}=1$,
which is a contradiction. Then,
for any $Q\in\mathscr D_N(Q_0)$,
\begin{align}\label{eq:positive-parent-gap-H}
1-U_{\widehat Q}>0.
\end{align}
From $U_R=1$ for any $R\in\mathscr D_N(Q_0)$ with $j(R)=N$
and $R\cap E\ne\varnothing$, we infer that, for any
$Q\in\mathscr D_N(Q_0)$ and $z\in E\cap Q$,
$\sum_{R\in\mathscr D_N(Q)}\delta_R\mathbf1_R(z)=1-U_{\widehat Q}$,
which, together with Lemma~\ref{lem:tree-equilibrium}, further implies that
\begin{align}
\sum_{R\in\mathscr D_N(Q)}\ell(R)^{d-sp}\delta_R^p
&=\int_{E\cap Q}\sum_{R\in\mathscr D_N(Q)}
\delta_R\mathbf1_{R}(z)\,d\mu_N(z)\notag\\
&=\left(1-U_{\widehat Q}\right)\mu_N(Q)
=\ell(Q)^{d-sp}\delta_Q^{p-1}\left(1-U_{\widehat Q}\right).
\label{eq:tree-local-energy-H}
\end{align}
For any $Q\in\mathscr S_0$ and $R\in\mathscr D_N(Q)$,
let $V_{\widehat Q}:=0$ and
$V_R:=\frac{U_R-U_{\widehat Q}}{1-U_{\widehat Q}}$.
From the constraints on $U$ and $E_Q\subset E\cap Q$, we infer that,
for any $R\in\mathscr D_N(Q)$ satisfying $j(R)=N$ and
$R\cap E_Q\ne\varnothing$, $U_R=1$ and hence $V_R=1$.
Using this and $V_{\widehat Q}=0$, we conclude that
$\{V_R\}_{R\in\mathscr D_N(Q)\cup\{\widehat Q\}}$
satisfies the constraints in the definition of $\mathsf C_N(E_Q;Q)$.
From \eqref{eq:dyadic-tree-thickness}, the definition of
$\mathsf C_N(E_Q;Q)$, and \eqref{eq:tree-local-energy-H},
it follows that, for any $Q\in\mathscr S_0$,
\begin{align*}
\alpha\ell(Q)^{d-sp} &\le\mathsf C_N(E_Q;Q) \le \left(1-U_{\widehat Q}\right)^{-p}\sum_{R\in\mathscr D_N(Q)} \ell(R)^{d-sp}\delta_R^p =\ell(Q)^{d-sp}\left(\frac{\delta_Q}{1-U_{\widehat Q}}\right)^{p-1},
\end{align*}
which, combined with \eqref{eq:positive-parent-gap-H}, further implies
that, for any $Q\in\mathscr S_0$,
\begin{align}\label{eq:drop-H}
\delta_Q\gtrsim\alpha^{\frac1{p-1}}\left(1-U_{\widehat Q}\right)>0.
\end{align}
By \eqref{eq:safe-constant-pointwise}, the integration with respect to
$\frac{d\mu_N(z)}{\mu_N(Q)}$ over $E\cap Q$, \eqref{eq:finite-current},
and \eqref{eq:curvature}, we have
\begin{align}\label{eq:safe-average-H}
\sum_{Q\in\mathscr S_0}\ell(Q)^{d-sp}|a_Q|^p
\lesssim\sum_{Q\in\mathscr S_0}\delta_Q^{1-p}
\int_{E\cap Q}|S_{j(Q)}u(z)|^p\,d\mu_N(z)+[u]_{s,p}^p.
\end{align}
From the triangle inequality, we deduce that
\begin{align}\label{eq:lowpass-split-H}
&\sum_{Q\in\mathscr S_0}\delta_Q^{1-p}
\int_{E\cap Q}|S_{j(Q)}u(z)|^p\,d\mu_N(z)\notag\\
&\quad\lesssim\sum_{Q\in\mathscr S_0}\delta_Q^{1-p}
\int_{E\cap Q}|S_{j(Q)}u(z)-S_Nu(z)|^p\,d\mu_N(z)
+\sum_{Q\in\mathscr S_0}\delta_Q^{1-p}
\int_{E\cap Q}|S_Nu(z)|^p\,d\mu_N(z).
\end{align}
Let $z\in\operatorname{supp}\mu_N$, $M:=N-j(Q_0)+1$, and
$R_n\in\mathscr D$ satisfy $z\in R_n$ and
$j(R_n)=j(Q_0)+n-1$ for any $n\in[0,M]\cap\mathbb Z$.
Since $\mu_N$ is finitely supported and $z\in\operatorname{supp}\mu_N$,
it follows that $\mu_N(R_n)>0$ for any $n\in[1,M]\cap\mathbb Z$.
From this, \eqref{eq:finite-current}, $R_0=\widehat{Q_0}$,
and the constraints at $j(R_M)=N$, we infer that
$0=U_{R_0}<U_{R_1}<\cdots<U_{R_M}=1$.
By $j(R_M)=N$, we have $S_{j(R_M)}u(z)-S_Nu(z)=0$.
Using the fact that $j(R_n)-j(R_{n-1})=1$ and
$S_j-S_{j-1}=\Delta_{j-1}$ for any $j\in\mathbb Z$,
we obtain, for any $n\in[1,M]\cap\mathbb Z$,
$$[S_{j(R_n)}u(z)-S_Nu(z)]
-[S_{j(R_{n-1})}u(z)-S_Nu(z)]
=\Delta_{j(R_n)-1}u(z).$$
From \eqref{eq:drop-H}, it follows that there exists $\beta\in(0,1]$
with $\beta\sim\alpha^{\frac1{p-1}}$ such that
$\delta_{R_n}\ge\beta(1-U_{R_{n-1}})$ for any
$n\in[1,M]\cap\mathbb Z$ satisfying $R_n\in\mathscr S_0$.
Applying Lemma~\ref{lem:rayhardy} to
$\{1-U_{R_n}\}_{n\in[0,M]\cap\mathbb Z}$ and
$\{S_{j(R_n)}u(z)-S_Nu(z)\}_{n\in[0,M]\cap\mathbb Z}$
with $G:=\{n\in[1,M]\cap\mathbb Z:R_n\in\mathscr S_0\}$,
$\beta^{-p}\sim\alpha^{-p'}$,
and integration with respect to $\mu_N$, we obtain
\begin{align}\label{eq:integrated-ray-H}
&\sum_{Q\in\mathscr S_0}\delta_Q^{1-p}
\int_{E\cap Q}|S_{j(Q)}u(z)-S_Nu(z)|^p\,d\mu_N(z)\notag\\
&\quad\lesssim\alpha^{-p'}
\sum_{\genfrac{}{}{0pt}{}{Q\in\mathscr D_N(Q_0)}{\delta_Q>0}}\delta_Q^{1-p}
\int_{E\cap Q}|\Delta_{j(Q)-1}u(z)|^p\,d\mu_N(z).
\end{align}
From \eqref{eq:finite-current}, Lemma~\ref{lem:sampling}, and
\eqref{eq:besovLP}, we deduce that
\begin{align}\label{eq:finite-energy-charge-H}
&\sum_{\genfrac{}{}{0pt}{}{Q\in\mathscr D_N(Q_0)}{\delta_Q>0}}\delta_Q^{1-p}
\int_{E\cap Q}|\Delta_{j(Q)-1}u(z)|^p\,d\mu_N(z)\notag\\
&\quad\le\sum_{Q\in\mathscr D(Q_0)}\ell(Q)^{d-sp}
\sup_{x\in Q}|\Delta_{j(Q)-1}u(x)|^p
\lesssim\sum_{j\in\mathbb Z}2^{jsp}\|\Delta_{j-1}u\|_{L^p(\mathbb{R}^d)}^p
\sim[u]_{s,p}^p.
\end{align}
It remains to consider the second term of \eqref{eq:lowpass-split-H}.
By \eqref{eq:finite-current} and $\operatorname{supp}\mu_N\subset E\subset F$,
we have
\begin{align}\label{eq:terminal-error-H}
\sum_{Q\in\mathscr S_0}\delta_Q^{1-p}
\int_{E\cap Q}|S_Nu(z)|^p\,d\mu_N(z)
\le\left[\sup_{z\in F}|S_Nu(z)|\right]^p
\sum_{Q\in\mathscr S_0}\ell(Q)^{d-sp}.
\end{align}
From \eqref{eq:lowpass-split-H}, \eqref{eq:integrated-ray-H},
\eqref{eq:finite-energy-charge-H}, and \eqref{eq:terminal-error-H},
we infer that
\begin{align}\label{eq:lowpass-zero-set-H}
\sum_{Q\in\mathscr S_0}\delta_Q^{1-p}
\int_{E\cap Q}|S_{j(Q)}u(z)|^p\,d\mu_N(z)
\lesssim\alpha^{-p'}[u]_{s,p}^p
+\left[\sup_{z\in F}|S_Nu(z)|\right]^p
\sum_{Q\in\mathscr S_0}\ell(Q)^{d-sp}.
\end{align}
By the definition of $S_N$, $u\in C_{\rm c}^\infty(\mathbb{R}^d)$,
and $u|_F=0$, we conclude that
$\|S_Nu-u\|_{L^\infty(\mathbb{R}^d)}\to0$ and hence
$\sup_{z\in F}|S_Nu(z)|\to0$ as $N\to\infty$. Using this,
the finiteness of $\mathscr S_0$, \eqref{eq:safe-average-H}, \eqref{eq:lowpass-zero-set-H},
and $1\le\alpha^{-p'}$ and letting $N\to\infty$ with
$\mathscr S_0$ fixed, we obtain
\begin{align}\label{eq:safe-finite-H}
\sum_{Q\in\mathscr S_0}\ell(Q)^{d-sp}|a_Q|^p
\lesssim\alpha^{-p'}[u]_{s,p}^p.
\end{align}
Let $\mathscr S_1\subset\mathscr S_a$ be finite and nonempty.
Let $m\in\mathbb Z$ satisfy
$m<\min_{Q\in\mathscr S_1}j(Q)$ and
$\operatorname{diam}(\bigcup_{Q\in\mathscr S_1}Q)<2^{-m}$.
Since the projection of $\bigcup_{Q\in\mathscr S_1}Q$ onto each
coordinate axis has diameter less than $2^{-m}$, it follows that
this set meets at most $2^d$ cubes $R\in\mathscr D$ satisfying $j(R)=m$.
Let
$$\{R_i\}_{i=1}^K:=\left\{R\in\mathscr D:j(R)=m\ \text{and}\
\mathscr S_1\cap\mathscr D(R)\ne\varnothing\right\}$$
with $K\in\mathbb N$.
From this and the definition of a dyadic grid, we infer that
$K\in[1,2^d]\cap\mathbb Z$,
$R_i\cap R_j=\varnothing$ for any distinct
$i,j\in[1,K]\cap\mathbb Z$, and
$\mathscr S_1\subset\bigcup_{i=1}^K\mathscr D(R_i)$.
From \eqref{eq:safe-finite-H} with $Q_0=R_i$ and
$\mathscr S_0=\mathscr S_1\cap\mathscr D(R_i)$ for any
$i\in[1,K]\cap\mathbb Z$ satisfying
$\mathscr S_1\cap\mathscr D(R_i)\ne\varnothing$,
from summation in $i$, and from $K\le2^d$, we deduce that
\eqref{eq:safe-finite-H} with $\mathscr S_1$ in place of $\mathscr S_0$ holds.
By this and the monotone convergence theorem, we obtain
\begin{align*}
\sum_{Q\in\mathscr S_a}\ell(Q)^{d-sp}|a_Q|^p
\lesssim\alpha^{-p'}[u]_{s,p}^p.
\end{align*}
From this, \eqref{eq:safeP}, and $A_Q\subset Q$, we deduce that
\begin{align}\label{eq:safePfinal}
\sum_{Q\in\mathscr S_a}\ell(Q)^{-sp}
\int_{A_Q}|P_Q(x)|^p\,dx\lesssim\alpha^{-p'}[u]_{s,p}^p.
\end{align}
Using \eqref{eq:lowpass-affine-decomp}, the triangle inequality,
\eqref{eq:curvature}, $\mathscr S=\mathscr S_a\cup\mathscr S_b$,
\eqref{eq:dangerPfinal}, \eqref{eq:safePfinal}, and $1\le\alpha^{-p'}$,
we obtain
\begin{align*}
\sum_{Q\in\mathscr S}\ell(Q)^{-sp}
\int_{A_Q}|S_{j(Q)}u(x)|^p\,dx
\lesssim\alpha^{-p'}[u]_{s,p}^p.
\end{align*}
By this, \eqref{eq:higher-high-selected}, and \eqref{eq:higherdecomp},
we conclude that \eqref{eq:dyadic-hardy-energy} holds.
This completes the proof of Proposition~\ref{prop:dyadic-hardy-energy}.
\end{proof}

Next, we are able to prove Theorem \ref{thm:generalp}(i) in the case
$s\ne1$.

\begin{proof}[Proof of Theorem~\ref{thm:generalp}(i) for $s\ne1$]
Let $p\in(1,\infty)$,
$s\in(0,\min\{1+\frac1p,\frac dp\})\setminus\{1\}$,
$\alpha\in(0,1]$, and $u\in C_{\rm c}^\infty(\Omega)$.
Let $x\in\Omega$ satisfy $d_{\alpha,s,p}(x)\in(0,\infty)$.
From the definition of $d_{\alpha,s,p}(x)$, we infer that there exists
$r_x\in[d_{\alpha,s,p}(x),2d_{\alpha,s,p}(x))$ satisfying
$$\frac{\operatorname{cap}_{s,p}(\overline{F\cap B(x,r_x)})}
{\operatorname{cap}_{s,p}(\overline{B(x,r_x)})}\ge\alpha.$$
Using the adjacent-grid covering property, we find that there exist
$t\in[1,3^d]\cap\mathbb Z$ and $Q\in\mathscr D_t$ satisfying
$\overline{B(x,r_x)}\subset Q^\circ$ and
$2r_x\le\ell(Q)\le C_{\mathrm{cov}}r_x$.
By the choice of $r_x$, we have $\overline{F\cap B(x,r_x)}\subset F\cap Q$ and
\begin{align}\label{eq:selection-cap}
\operatorname{cap}_{s,p}\left(\overline{F\cap B(x,r_x)}\right)
\ge\alpha\kappa_{d,s,p}r_x^{d-sp}
\ge\alpha\kappa_{d,s,p}C_{\mathrm{cov}}^{-(d-sp)}\ell(Q)^{d-sp}.
\end{align}
Let $c:=\min\{1,\kappa_{d,s,p}C_{\mathrm{cov}}^{-(d-sp)}\}$ and
$\mathscr A$ be the family of all cubes $Q\in\mathscr Q$
for which there exists a compact set $E_Q\subset F\cap Q$ such that
$\operatorname{cap}_{s,p}(E_Q)\ge c\alpha\ell(Q)^{d-sp}$.
Using \eqref{eq:selection-cap} and the choice of $Q$, we find that
there exists a positive constant $C$ such that each $x\in\Omega$
with $d_{\alpha,s,p}(x)<\infty$ admits a cube $Q\in\mathscr A$
satisfying $x\in Q$ and
$C^{-1}\ell(Q)\le d_{\alpha,s,p}(x)\le C\ell(Q)$.
Let $\{Q_m\}_{m\in\mathbb N}$ enumerate $\mathscr Q$ without repetitions.
For any $m\in\mathbb N$, let
$$G_m:=\left\{x\in Q_m\cap\Omega:C^{-1}\ell(Q_m)\le d_{\alpha,s,p}(x)\le C\ell(Q_m)\right\}$$
if $Q_m\in\mathscr A$, and let $G_m:=\varnothing$ otherwise. Let
$A_{Q_m}:=G_m\setminus\bigcup_{k=1}^{m-1}G_k$.
Since $d_{\alpha,s,p}$ is measurable by Lemma~\ref{lem:distance-measurable},
it follows that the sets $A_{Q_m}$ are measurable, pairwise disjoint,
and satisfy $A_Q=\varnothing$ for any $Q\notin\mathscr A$.
Using the definition of $G_m$ and $A_{Q_m}\subset G_m$, we obtain,
for any $Q\in\mathscr Q$ and $x\in A_Q$,
\begin{align}\label{eq:selection-scale}
C^{-1}\ell(Q)\le d_{\alpha,s,p}(x)\le C\ell(Q).
\end{align}
From the existence of an admissible cube for each $x$, we deduce that
$$\bigcup_{m\in\mathbb N}A_{Q_m}
=\left\{x\in\Omega:d_{\alpha,s,p}(x)\in(0,\infty)\right\}.$$
Using \eqref{eq:selection-scale}, the disjointness of the sets $A_Q$,
and $\infty^{-sp}=0$, we obtain
\begin{align}\label{eq:selection-integral}
\int_\Omega\frac{|u(x)|^p}{d_{\alpha,s,p}(x)^{sp}}\,dx
=\sum_{Q\in\mathscr A}\int_{A_Q}
\frac{|u(x)|^p}{d_{\alpha,s,p}(x)^{sp}}\,dx
\lesssim\sum_{Q\in\mathscr A}\ell(Q)^{-sp}
\int_{A_Q}|u(x)|^p\,dx.
\end{align}
For every $Q\in\mathscr A$, let $E_Q$ be such a compact set.
By Lemma~\ref{lem:capcompare}, we have, for every
$N\in(j(Q),\infty)\cap\mathbb Z$,
\begin{align*}
\operatorname{cap}_{s,p}(E_Q)\lesssim\mathsf C_N(E_Q;Q).
\end{align*}
Thus, \eqref{eq:dyadic-tree-thickness} holds with $c_1\alpha$ in place of
$\alpha$, for a fixed $c_1\in(0,c]$.
This, combined with Proposition~\ref{prop:dyadic-hardy-energy} and
\eqref{eq:selection-integral}, then completes the proof of
Theorem~\ref{thm:generalp}(i) in the case $s\ne1$.
\end{proof}

\section{Sharpness of the Dependence on $\alpha$}
\label{sec:alpha-sharpness}

For any $n\in[1,\infty)\cap\mathbb N$ and any nonempty
$U\subset\mathbb{R}^n$, let $\operatorname{diam}(U):=\sup_{x,y\in U}|x-y|$,
and let $\operatorname{diam}(\varnothing):=0$.
For any $q\in(0,n]$, $\delta\in(0,\infty)$, and
$A\subset\mathbb{R}^n$, we define the \emph{$q$-dimensional Hausdorff measure}
$\mathcal H^q$ by
\begin{align*}
\mathcal H^q(A)&:=\lim_{\delta\to0^+}\mathcal H^q_\delta(A)
:=\lim_{\delta\to0^+}\inf_{\genfrac{}{}{0pt}{}{A\subset\bigcup_{i\in\mathbb N}U_i}{U_i\subset\mathbb{R}^n,\ \operatorname{diam}(U_i)
\le\delta,\,\forall\, i\in\mathbb N}}
\sum_{i\in\mathbb N}\operatorname{diam}(U_i)^q.
\end{align*}

In this section, let $p\in(1,\infty)$, $s\in(0,\min\{1+\frac1p,\frac dp\})$,
and $q:=d-sp$.
Let $\tau\in(0,d]$. A compact set $E\subset\mathbb R^d$ with
$\operatorname{diam}(E)\in(0,\infty)$ is said to be \emph{$\tau$-regular}
if there exists a positive constant $C$ such that, for any
$x\in E$ and $r\in(0,\operatorname{diam}(E)]$,
\begin{align*}
C^{-1}r^\tau
\le\mathcal H^\tau(E\cap B(x,r))
\le Cr^\tau.
\end{align*}
Let $Q:=[0,\frac1{8\sqrt d}]^d$.
Using the $d$-regularity of $Q$ and
\cite[Theorem~3.1]{MattilaSaaranen2009}, we find that there exists
a nonempty compact $q$-regular set $E\subset Q\subset B(0,\frac14)$.
Let $\mu:=\mathcal H^q|_E$, $\rho(x):=\operatorname{dist}(x,E)$ for any
$x\in\mathbb{R}^d$, and $E_r:=\{x\in\mathbb{R}^d:\rho(x)<r\}$ for any
$r\in(0,\infty)$.
By the $q$-regularity of $E$, we have $\mu(E)\sim1$ and, for any
$y\in E$ and $r\in(0,1]$,
\begin{align}\label{eq:compact-sharp-regularity}
\mu(B(y,r))\sim r^q.
\end{align}
For any $r\in(0,1]$, let $Y_r\subset E$ be maximal among families
satisfying $|a-b|\ge r$ for any distinct $a,b\in Y_r$.
Using this maximality, the disjointness of $\{B(a,\frac r3)\}_{a\in Y_r}$,
and \eqref{eq:compact-sharp-regularity}, we obtain
$\#Y_r\sim r^{-q}$ and
$$\bigcup_{a\in Y_r}B\left(a,\frac r3\right)\subset E_r
\subset\bigcup_{a\in Y_r}B(a,2r),$$ which further implies that
\begin{align}\label{eq:compact-sharp-volume}
|E_r|\sim r^{d-q}=r^{sp}.
\end{align}
For any $\varepsilon\in(0,\frac1{16}]$, let
\begin{align*}
F_\varepsilon:=\bigcup_{a\in Y_\varepsilon}
\overline{B\left(a,\frac\varepsilon8\right)}
\quad\text{and}\quad
\Omega_\varepsilon:=B(0,1)\setminus F_\varepsilon.
\end{align*}
Then we have the following conclusion.

\begin{lemma}\label{lem:compact-sharp-capacity}
There exists a positive constant $C$ such that, for any
$\varepsilon\in(0,\frac1{64}]$, $y\in E$, and
$r\in[8\varepsilon,\frac18]$,
\begin{align}\label{eq:compact-sharp-capacity}
\operatorname{cap}_{s,p}\left(F_\varepsilon\cap\overline{B(y,r)}\right)
\ge Cr^q\left(\log\frac r\varepsilon\right)^{1-p}.
\end{align}
\end{lemma}

\begin{proof}
Let $\varepsilon\in(0,\frac1{64}]$, $y\in E$,
$r\in[8\varepsilon,\frac18]$, $\mathcal Y:=Y_\varepsilon\cap B(y,\frac r2)$, and
$M:=\#\mathcal Y$.
By the maximality of $Y_\varepsilon$ and $\varepsilon\le\frac r8$,
we have $E\cap B(y,\frac r4)\subset
\bigcup_{a\in\mathcal Y}B(a,\varepsilon)$.
Using this, \eqref{eq:compact-sharp-regularity}, and the disjointness
of $\{B(a,\frac\varepsilon3)\}_{a\in\mathcal Y}$, we obtain
\begin{align*}
r^q
&\sim\mu\left(B\left(y,\frac r4\right)\right)
\le\sum_{a\in\mathcal Y}\mu(B(a,\varepsilon))
\sim M\varepsilon^q
\end{align*}
and
\begin{align*}
M\varepsilon^q
&\sim\sum_{a\in\mathcal Y}\mu\left(B\left(a,\frac\varepsilon3\right)\right)
\le\mu(B(y,r))\sim r^q.
\end{align*}
Let $b\in C_{\rm c}^\infty(B(0,\frac1{16}))$ satisfy $b\ge0$ and
$\int_{\mathbb{R}^d}b(x)\,dx=1$.
For any $x\in\mathbb{R}^d$, let
$$f(x):=M^{-1}\sum_{a\in\mathcal Y}\frac1{\varepsilon^{d}}
b\left(\frac{x-a}{\varepsilon}\right)$$
and $d\nu(x):=f(x)\,dx$.
From these definitions and $\varepsilon\le\frac r8$, it follows that
\begin{align}\label{eq:compact-sharp-measure}
\nu(\mathbb{R}^d)=1
\quad\text{and}\quad
\operatorname{supp}\nu\subset F_\varepsilon\cap B\left(y,\frac{3r}4\right).
\end{align}
By the definition of $f$, $M\sim(\frac r\varepsilon)^q$,
and $q=d-sp$, we have,
for any $x\in\mathbb{R}^d$ and $t\in(0,\varepsilon]$,
\begin{align}\label{eq:compact-sharp-growth-small}
\nu(B(x,t)) &=\int_{B(x,t)}f(z)\,dz \le\|f\|_{L^\infty(\mathbb{R}^d)}|B(x,t)|
=\frac1{M\varepsilon^{d}}\|b\|_{L^\infty(\mathbb{R}^d)}|B(x,t)|
\sim \frac{t^d}{r^{q}\varepsilon^{sp}}.
\end{align}
From \eqref{eq:compact-sharp-regularity} and $\mu(E)\sim1$,
we infer that $\mu(B(x,t))\lesssim t^q$ for any $x\in\mathbb{R}^d$
and $t\in(0,\infty)$.
Using this and the disjointness of
$\{B(a,\frac\varepsilon3)\}_{a\in Y_\varepsilon}$, we obtain,
for any $x\in\mathbb{R}^d$ and $t\in[\varepsilon,\infty)$,
\begin{align*}
\varepsilon^q\#(Y_\varepsilon\cap B(x,t))
\lesssim\sum_{a\in Y_\varepsilon\cap B(x,t)}
\mu\left(B\left(a,\frac\varepsilon3\right)\right)\le\mu\left(B\left(x,t+\frac\varepsilon3\right)\right)\lesssim t^q.
\end{align*}
From this,
$\operatorname{supp}[b(\frac{\cdot-a}{\varepsilon})]\subset
B(a,\frac\varepsilon{16})$ for any $a\in\mathcal Y$,
the definition of $f$, and \eqref{eq:compact-sharp-measure},
it follows that, for any $x\in\mathbb{R}^d$ and $t\in[\varepsilon,\infty)$,
\begin{align}\label{eq:compact-sharp-growth-large}
\nu(B(x,t))
\le\min\left\{1,
\frac{\#(\mathcal Y\cap B(x,t+\frac\varepsilon{16}))}{M}\right\}
\lesssim\min\left\{1,\frac{t^q}{r^{q}}\right\}.
\end{align}

We first assume $s=1$, and hence $p\in(1,d)$ and $q=d-p$.
For any $x\in\mathbb{R}^d$, let
$$I_1\nu(x):=\int_{\mathbb{R}^d}\frac 1{|x-z|^{d-1}}\,d\nu(z).$$
By Wolff's inequality (see, for instance,
\cite[pp.~164--165]{HedbergWolff1983}),
\eqref{eq:compact-sharp-growth-small},
\eqref{eq:compact-sharp-growth-large}, and
\eqref{eq:compact-sharp-measure}, we have
\begin{align*}
\|I_1\nu\|_{L^{p'}(\mathbb{R}^d)}^{p'}
&\sim\int_{\mathbb{R}^d}\int_0^\infty
\left[\frac{\nu(B(x,t))}{t^{d-p}}\right]^{\frac1{p-1}}
\frac{dt}{t}\,d\nu(x)\\
&\lesssim \frac1{r^{\frac q{p-1}}}
\left[\int_0^\varepsilon\left(\frac t\varepsilon\right)^{p'}\frac{dt}{t}
+\int_\varepsilon^r\frac{dt}{t}
+\int_r^\infty\left(\frac rt\right)^{\frac q{p-1}}\frac{dt}{t}\right]
\sim \frac1{r^{\frac q{p-1}}}\log\frac r\varepsilon.
\end{align*}
Using this, the fundamental theorem of calculus, averaging over
the unit sphere, Tonelli's theorem, and H\"older's inequality,
we obtain, for any $\varphi\in C_{\rm c}^\infty(\mathbb{R}^d)$,
\begin{align}\label{eq:compact-sharp-pairing}
\left|\int_{\mathbb{R}^d}\varphi(x)\,d\nu(x)\right|
&\lesssim\int_{\mathbb{R}^d}|\nabla\varphi(x)|I_1\nu(x)\,dx
\le\|\nabla\varphi\|_{L^p(\mathbb{R}^d)}
\|I_1\nu\|_{L^{p'}(\mathbb{R}^d)} \lesssim[\varphi]_{s,p}
\,\frac1{r^{\frac qp}} \left(\log\frac r\varepsilon\right)^{\frac1{p'}}.
\end{align}

Now, we assume that $s\ne1$.
For any $j\in\mathbb Z$, let
$\widetilde\Delta_j:=\sum_{k=j-2}^{j+2}\Delta_k$.
From the fact that $\widetilde\Delta_j\Delta_j=\Delta_j$
for any $j\in\mathbb Z$, the uniform convergence of $\sum_{j\in\mathbb Z}\Delta_j\varphi$
to $\varphi$, H\"older's inequalities, and \eqref{eq:besovLP}, we infer that,
for any $\varphi\in C_{\rm c}^\infty(\mathbb{R}^d)$,
\begin{align}\label{920}
\left|\int_{\mathbb{R}^d}\varphi(x)f(x)\,dx\right| &=\left|\sum_{j\in\mathbb Z}\int_{\mathbb{R}^d} \Delta_j\varphi(x)\widetilde\Delta_jf(x)\,dx\right| \lesssim[\varphi]_{s,p} \left[\sum_{j\in\mathbb Z}2^{-jsp'} \left\|\widetilde\Delta_jf\right\|_{L^{p'}(\mathbb{R}^d)}^{p'}\right]^{\frac1{p'}}.
\end{align}
Let $K\in\mathcal S(\mathbb{R}^d)$ be the convolution kernel of
$\widetilde\Delta_0$, and let $K_t(x):=\frac1{t^{d}}K(\frac xt)$ for any
$t\in(0,\infty)$ and $x\in\mathbb{R}^d$.
By Young's inequality, \eqref{eq:compact-sharp-measure},
the definition of $f$,
and $M\sim(\frac r\varepsilon)^q$, we have, for any $t\in(0,\infty)$,
\begin{align}\label{eq:compact-sharp-kernel-large}
\|K_t*f\|_{L^{p'}(\mathbb{R}^d)}
&\le\|K_t\|_{L^{p'}(\mathbb{R}^d)}\|f\|_{L^1(\mathbb{R}^d)}
\sim \frac1{t^{\frac dp}}
\end{align}
and
\begin{align}\label{eq:compact-sharp-kernel-small}
\|K_t*f\|_{L^{p'}(\mathbb{R}^d)}
&\le \frac1{M^{\frac1p}\varepsilon^{\frac dp}}
\|K_t\|_{L^1(\mathbb{R}^d)}\|b\|_{L^{p'}(\mathbb{R}^d)}
\sim \frac1{r^{\frac qp}\varepsilon^{s}}.
\end{align}
Using the Schwartz decay of $K$ and
\eqref{eq:compact-sharp-growth-large}, we obtain, for any
$t\in[\varepsilon,r]$ and $x\in\mathbb{R}^d$,
\begin{align*}
|K_t*f(x)| &\lesssim \frac1{t^{d}}\sum_{m\in\mathbb Z_+}
\frac1{2^{m(d+1)}} \nu\left(B\left(x,2^{m+1}t\right)\right)
\lesssim \frac1{r^{q}t^{sp}} \sum_{m\in\mathbb Z_+}\frac1{2^{m(d+1-q)}}
\lesssim \frac1{r^{q}t^{sp}}.
\end{align*}
From this, Young's inequality,
and H\"older's inequality, it follows that, for any $t\in[\varepsilon,r]$,
\begin{align}\label{eq:compact-sharp-kernel-middle}
\|K_t*f\|_{L^{p'}(\mathbb{R}^d)}
&\le\|K_t*f\|_{L^\infty(\mathbb{R}^d)}^{\frac1p}
\|K_t*f\|_{L^1(\mathbb{R}^d)}^{\frac1{p'}}
\lesssim \frac1{r^{\frac qp}t^{s}}
\|K_t\|_{L^1(\mathbb{R}^d)}^{\frac1{p'}}
\|f\|_{L^1(\mathbb{R}^d)}^{\frac1{p'}}\sim \frac1{r^{\frac qp}t^{s}}.
\end{align}
Using $\widetilde\Delta_jf=K_{2^{-j}}*f$,
\eqref{eq:compact-sharp-kernel-small},
\eqref{eq:compact-sharp-kernel-middle},
\eqref{eq:compact-sharp-kernel-large}, $q=d-sp$, and $r\ge8\varepsilon$,
we obtain
\begin{align*}
\sum_{j\in\mathbb Z}2^{-jsp'} \left\|\widetilde\Delta_jf\right\|_{L^{p'}(\mathbb{R}^d)}^{p'}
&\lesssim \frac1{r^{\frac q{p-1}}} \left[ \sum_{\genfrac{}{}{0pt}{}{j\in\mathbb Z}{2^{-j}\le\varepsilon}} \left(\frac{2^{-j}}\varepsilon\right)^{sp'} +\sum_{\genfrac{}{}{0pt}{}{j\in\mathbb Z}{\varepsilon<2^{-j}<r}}1 +\sum_{\genfrac{}{}{0pt}{}{j\in\mathbb Z}{2^{-j}\ge r}} \left(\frac r{2^{-j}}\right)^{\frac q{p-1}} \right]\\
&\sim \frac1{r^{\frac q{p-1}}}
\left(1+\#\left\{j\in\mathbb Z:\varepsilon<\frac1{2^{j}}
<r\right\}\right) \sim \frac1{r^{\frac q{p-1}}}\log\frac r\varepsilon,
\end{align*}
which, together with \eqref{920}, further implies
that \eqref{eq:compact-sharp-pairing} when $s\ne1$ holds.
By this and \eqref{eq:compact-sharp-measure}, we have, for any
$\varphi\in C_{\rm c}^\infty(\mathbb{R}^d)$ satisfying $\varphi\ge1$ on an open
neighbourhood of $F_\varepsilon\cap\overline{B(y,r)}$,
\begin{align*}
1\le\int_{\mathbb{R}^d}\varphi(x)\,d\nu(x)
\lesssim[\varphi]_{s,p}r^{-\frac qp}
\left(\log\frac r\varepsilon\right)^{\frac1{p'}}.
\end{align*}
Taking the infimum and using the definition of $\operatorname{cap}_{s,p}$,
we obtain \eqref{eq:compact-sharp-capacity}.
This completes the proof of Lemma~\ref{lem:compact-sharp-capacity}.
\end{proof}

\begin{proof}[Proof of Theorem~\ref{thm:generalp}(ii)]
For any $N\in[2,\infty)\cap\mathbb N$, let
$\varepsilon_N:=4^{-(3N+2)}$ and
$K_N:=E_{4^{-(N+1)}}\setminus E_{4^{-2N}}$.
Let $b\in C_{\rm c}^\infty(B(0,\frac1{16}))$ satisfy $b\ge0$ and
$\int_{\mathbb{R}^d}b(x)\,dx=1$.
For any $t\in(0,\infty)$ and $x\in\mathbb{R}^d$,
let $b_t(x):=\frac1{t^{d}}b(\frac xt)$.
For any $j\in\mathbb N$, let $r_j:=\frac1{4^{j}}$ and
$\eta_j:=b_{\frac{r_j}4}*\mathbf1_{E_{\frac{r_j}2}}$.
By these definitions, Young's inequality, and
\eqref{eq:compact-sharp-volume}, we have
$0\le\eta_j\le1$, $\eta_j=1$ on $E_{\frac{r_j}4}$,
$\operatorname{supp}\eta_j\subset E_{\frac{3r_j}4}$, and, for any $m\in\{0,1,2\}$,
\begin{align}\label{eq:compact-sharp-cutoffs}
\|D^m\eta_j\|_{L^p(\mathbb{R}^d)}
\le\left\|D^mb_{\frac{r_j}4}\right\|_{L^1(\mathbb{R}^d)}
\left|E_{\frac{r_j}2}\right|^{\frac1p}\lesssim r_j^{s-m}.
\end{align}
Let
\begin{align}\label{eq:compact-sharp-test}
u_N:=\frac1N\sum_{j=1}^N(\eta_j-\eta_{2N+j}).
\end{align}
Using $\operatorname{supp}\eta_{j+1}\subset E_{\frac{3r_j}{16}}
\subset E_{\frac{r_j}4}$ and $\eta_j=1$ on $E_{\frac{r_j}4}$,
we obtain $\eta_j\ge\eta_{j+1}$, and hence $0\le u_N\le1$.
By $\frac{r_{3N}}4=4\varepsilon_N$, $\eta_j=1$ on $E_{\frac{r_j}4}$,
and $\operatorname{supp}\eta_j\subset E_{\frac{3r_j}4}$, we have
$u_N=0$ on $E_{4\varepsilon_N}$ and
$\operatorname{supp} u_N\subset E_{\frac3{16}}\subset B(0,1)$.
From this and $F_{\varepsilon_N}\subset E_{\frac{\varepsilon_N}4}$,
we infer that $u_N\in C_{\rm c}^\infty(\Omega_{\varepsilon_N})$.
By the definition of $K_N$, we have, for any $x\in K_N$ and
 $j\in[1,N]\cap\mathbb N$,
$\rho(x)<4^{-(N+1)}=\frac{r_N}4\le\frac{r_j}4$ and
$\rho(x)\ge4^{-2N}>\frac{3r_{2N+j}}4$.
Using this and the definition of $\eta_j$,
we obtain $\eta_j(x)=1$ and $\eta_{2N+j}(x)=0$, which imply
that $u_N(x)=1$ for any $x\in K_N$.

For $s=1$, we have $\operatorname{supp}\nabla\eta_j\subset
\{x\in\mathbb{R}^d:\rho(x)\in[\frac{r_j}4,\frac{3r_j}4]\}$ for any
$j\in\mathbb N$.
From this and $r_{j+1}=\frac{r_j}4$, we deduce that
$\operatorname{supp}\nabla\eta_j\cap\operatorname{supp}\nabla\eta_k=\varnothing$ for any
distinct $j,k\in\mathbb N$.
Using this, \eqref{eq:compact-sharp-test}, and
\eqref{eq:compact-sharp-cutoffs}, we obtain
\begin{align*}
[u_N]_{1,p}^p
=\frac1{N^p}\sum_{j=1}^N
\left[\|\nabla\eta_j\|_{L^p(\mathbb{R}^d)}^p
+\|\nabla\eta_{2N+j}\|_{L^p(\mathbb{R}^d)}^p\right]
\lesssim N^{1-p}.
\end{align*}
For $s\ne1$, by the $L^p$ boundedness of $\Delta_k$, Bernstein's inequality,
and \eqref{eq:compact-sharp-cutoffs}, we have, for any
$k\in\mathbb Z$ and $j\in\mathbb N$,
\begin{align*}
2^{ks}\|\Delta_k\eta_j\|_{L^p(\mathbb{R}^d)}
&\lesssim2^{ks}\min\left\{\|\eta_j\|_{L^p(\mathbb{R}^d)},
\frac1{2^{2k}}\left\|D^2\eta_j\right\|_{L^p(\mathbb{R}^d)}\right\}\\
&\lesssim\min\left\{2^{(k-2j)s},\frac1{2^{(k-2j)(2-s)}}\right\}
\le\frac1{2^{\min\{s,2-s\}|k-2j|}}.
\end{align*}
From this, \eqref{eq:compact-sharp-test}, and the triangle
inequality, it follows that, for any $k\in\mathbb Z$,
\begin{align*}
2^{ks}\|\Delta_ku_N\|_{L^p(\mathbb{R}^d)}
\lesssim\frac1N\sum_{j=1}^N
\left(\frac1{2^{\min\{s,2-s\}|k-2j|}}+\frac1{2^{\min\{s,2-s\}|k-4N-2j|}}\right).
\end{align*}
By this, Young's inequality on $\ell^p(\mathbb Z)$, and
\eqref{eq:besovLP}, we have
\begin{align*}
[u_N]_{s,p}^p
\sim\sum_{k\in\mathbb Z}2^{ksp}
\|\Delta_ku_N\|_{L^p(\mathbb{R}^d)}^p\lesssim\frac{2N}{N^p}
\left(\sum_{m\in\mathbb Z}\frac1{2^{\min\{s,2-s\}|m|}}\right)^p
\lesssim \frac1{N^{p-1}}.
\end{align*}
From \eqref{eq:compact-sharp-volume}, it follows that
$\int_{r_{2N}}^{r_{N+1}}|E_t|t^{-sp}\frac{dt}{t}
\sim\log\frac{r_{N+1}}{r_{2N}}=(N-1)\log4$
and $$\left|E_{r_{2N}}\right|\int_{r_{2N}}^{r_{N+1}}\frac1{t^{sp}}\frac{dt}{t}
=\frac{|E_{r_{2N}}|}{sp}\left(\frac1{r_{2N}^{sp}}-\frac1{r_{N+1}^{sp}}\right)
\lesssim1.$$
Using this and Tonelli's theorem, we find that there exists
$N_0\in[4,\infty)\cap\mathbb N$ such that, for any
$N\in[N_0,\infty)\cap\mathbb N$,
\begin{align*}
\int_{K_N}\frac{1}{\rho(x)^{sp}}\,dx
&\ge sp\int_{r_{2N}}^{r_{N+1}} \frac1{t^{sp+1}}|E_t\setminus E_{r_{2N}}|
\,dt =sp\int_{r_{2N}}^{r_{N+1}} \frac1{t^{sp}}(|E_t|-|E_{r_{2N}}|)
\,\frac{dt}{t}\gtrsim N.
\end{align*}
By these estimates, we conclude that, for any
$N\in[N_0,\infty)\cap\mathbb N$,
\begin{align}\label{eq:compact-sharp-profile}
[u_N]_{s,p}^p\lesssim \frac1{N^{p-1}}
\quad\text{and}\quad
\int_{K_N}\frac{1}{\rho(x)^{sp}}\,dx\gtrsim N.
\end{align}

Let $C$ be the constant in \eqref{eq:compact-sharp-capacity}.
Let $\gamma:=\min\{1,\frac{C(3\log4)^{1-p}}{3^q\kappa_{d,s,p}}\}$
and $\alpha_*:=\gamma(N_0+1)^{1-p}$.
For any $\alpha\in(0,\alpha_*]$, let
$N:=\lfloor(\frac\gamma\alpha)^{\frac1{p-1}}\rfloor$.
From these definitions, we infer that $N\in[N_0,\infty)\cap\mathbb N$,
$\alpha\le\gamma N^{1-p}$, and $N\sim\alpha^{-\frac1{p-1}}$.
Moreover, $\Omega_{\varepsilon_N}$ is bounded, has
smooth boundary, and is connected whenever
$d\in[2,\infty)\cap\mathbb N$.
For any $x\in K_N$, let $y\in E$ satisfy
$|x-y|=\rho(x)$.
By the monotonicity of capacity and
\eqref{eq:compact-sharp-capacity},
we conclude that
\begin{align*}
\operatorname{cap}_{s,p}\left(\overline{\left(\mathbb{R}^d\setminus\Omega_{\varepsilon_N}\right)\cap B(x,3\rho(x))}\right)
&\ge\operatorname{cap}_{s,p}\left(F_{\varepsilon_N}\cap\overline{B(y,\rho(x))}\right)
\ge C\rho(x)^q\left[\log\frac {\rho(x)}{\varepsilon_N}\right]^{1-p}\\
&\ge C(3\log4)^{1-p}N^{1-p}\rho(x)^q
\ge\alpha\kappa_{d,s,p}3^q\rho(x)^q.
\end{align*}
From this and the definition of $d_{\alpha,s,p}$, we infer that
$d_{\alpha,s,p}(x)\le3\rho(x)$ for any $x\in K_N$.
Using this, $u_N=1$ on $K_N$, \eqref{eq:compact-sharp-profile},
and the definition of $\mathcal H_{d,s,p}$, we obtain
\begin{align*}
\mathcal H_{d,s,p}(\alpha) &\ge\frac{1}{[u_N]_{s,p}^p}\int_{\Omega_{\varepsilon_N}} \frac{|u_N(x)|^p}{d_{\alpha,s,p}(x)^{sp}}\,dx \ge\frac{1}{[u_N]_{s,p}^p 3^{sp}}
\int_{K_N}\frac{1}{\rho(x)^{sp}}\,dx \gtrsim\frac{N}{N^{1-p}}=N^p\sim\alpha^{-p'}.
\end{align*}
This completes the proof of Theorem~\ref{thm:generalp}(ii).
\end{proof}

\section{Optimality of the Range of $s$}\label{sec:counterexamples}

\begin{lemma}\label{lem:thickE}
Let $p\in(1,d-1)$, $s\in[1+\frac1p,\frac dp)$,
$q:=d-p(s-1)$, and
$E\subset\mathbb{R}^d$ be a compact $q$-regular set.
Then there exists a positive constant $C$ such that, for any
$y\in E$ and $r\in(0,\operatorname{diam}(E)]$,
\begin{align}\label{eq:thickE}
\operatorname{cap}_{s,p}(E\cap B(y,r))\ge Cr^{d-sp}.
\end{align}
\end{lemma}

\begin{proof}
Let $m:=\lfloor s\rfloor$,
$\sigma:=s-m\in[0,1)$, $y\in E$, $r\in(0,\operatorname{diam}(E)]$, and
$\nu:=\mathcal H^q|_{E\cap\overline{B(y,\frac r2)}}$.
By the $q$-regularity of $E$, we have, for any $z\in\mathbb{R}^d$ and $t\in(0,\infty)$,
\begin{align}\label{eq:nu-mass}
\nu(\mathbb{R}^d)\sim r^q
\quad\text{and}\quad
\nu(B(z,t))\lesssim\min\{t^q,r^q\}.
\end{align}
Let $\mathscr D$ be a dyadic grid and
$j_0\in\mathbb Z$ satisfy $2^{-j_0-1}<r\le2^{-j_0}$.
From H\"older's inequality, Lemma~\ref{lem:sampling}, and
\eqref{eq:nu-mass}, it follows that, for any
$j\in\mathbb Z$ and $\varphi\in C_{\rm c}^\infty(\mathbb{R}^d)$,
\begin{align}\label{eq:nu-frequency-bound}
\left|\int_{\mathbb{R}^d}\Delta_j\varphi(z)\,d\nu(z)\right|
&\le\left[\sum_{\genfrac{}{}{0pt}{}{Q\in\mathscr D}{j(Q)=j}}
|Q|\sup_{z\in Q}|\Delta_j\varphi(z)|^p\right]^{\frac1p}
\left[2^{\frac{jd}{p-1}}
\sum_{\genfrac{}{}{0pt}{}{Q\in\mathscr D}{j(Q)=j}}\nu(Q)^{p'}\right]^{\frac1{p'}}\notag\\
&\lesssim\|\Delta_j\varphi\|_{L^p(\mathbb{R}^d)}
\left[2^{\frac{jd}{p-1}}
\sum_{\genfrac{}{}{0pt}{}{Q\in\mathscr D}{j(Q)=j}}\nu(Q)^{p'}\right]^{\frac1{p'}}\notag\\
&\lesssim\|\Delta_j\varphi\|_{L^p(\mathbb{R}^d)}
\begin{cases}
r^q2^{\frac{jd}{p}},&j\le j_0,\\
r^{\frac q{p'}}2^{\frac{j(d-q)}p},&j>j_0.
\end{cases}
\end{align}
Let $j\in\mathbb Z$ and $\varphi\in C_{\rm c}^\infty(\mathbb{R}^d)$.
Using Bernstein's inequality,
we obtain
\begin{align}\label{eq:counter-bernstein}
\|\Delta_j\varphi\|_{L^p(\mathbb{R}^d)}
\sim2^{-jm}\|D^m\Delta_j\varphi\|_{L^p(\mathbb{R}^d)}
=2^{-jm}\|\Delta_jD^m\varphi\|_{L^p(\mathbb{R}^d)}.
\end{align}
From the $L^p$ boundedness of $\Delta_j$ and the definition of
$[\,\cdot\,]_{s,p}$, we infer that, if $\sigma=0$,
\begin{align}\label{eq:counter-integer-frequency}
\left\|\Delta_jD^m\varphi\right\|_{L^p(\mathbb{R}^d)}
\lesssim\left\|D^m\varphi\right\|_{L^p(\mathbb{R}^d)}=[\varphi]_{s,p}.
\end{align}
By the definition of $[\,\cdot\,]_{s,p}$ with $s=m+\sigma$
and by \eqref{eq:besovLP}, we conclude that, if $\sigma\in(0,1)$,
\begin{align*}
2^{j\sigma}\left\|\Delta_jD^m\varphi\right\|_{L^p(\mathbb{R}^d)}
\le\left[\sum_{k\in\mathbb Z}2^{k\sigma p}
\left\|\Delta_kD^m\varphi\right\|_{L^p(\mathbb{R}^d)}^p\right]^{\frac1p}
\sim\left[D^m\varphi\right]_{\sigma,p}=[\varphi]_{s,p}.
\end{align*}
From this, \eqref{eq:counter-bernstein},
\eqref{eq:counter-integer-frequency}, and $s=m+\sigma$, we infer
that, for any $j\in\mathbb Z$ and
$\varphi\in C_{\rm c}^\infty(\mathbb{R}^d)$,
\begin{align*}
\|\Delta_j\varphi\|_{L^p(\mathbb{R}^d)}
\sim\frac1{2^{jm}}\|\Delta_jD^m\varphi\|_{L^p(\mathbb{R}^d)}
\lesssim\frac1{2^{js}}[\varphi]_{s,p}.
\end{align*}
By this, $\varphi=\sum_{j\in\mathbb Z}\Delta_j\varphi$,
\eqref{eq:nu-frequency-bound},
$d-sp\in(0,\infty)$, and $s-\frac{d-q}{p}=1$, we conclude that
\begin{align*}
\left|\int_{\mathbb{R}^d}\varphi(z)\,d\nu(z)\right|
\lesssim[\varphi]_{s,p}
\left[r^q\sum_{\genfrac{}{}{0pt}{}{j\in\mathbb Z}{j\le j_0}}2^{\frac{j(d-sp)}p}
+r^{\frac q{p'}}\sum_{\genfrac{}{}{0pt}{}{j\in\mathbb Z}{j>j_0}}\frac1{2^{j}}\right]
\lesssim r^{q-\frac{d-sp}p}[\varphi]_{s,p}.
\end{align*}
Using this, \eqref{eq:nu-mass}, and the fact that
$\nu(\mathbb{R}^d)\le\int_{\mathbb{R}^d}\varphi\,d\nu$ for any
$\varphi\in C_{\rm c}^\infty(\mathbb{R}^d)$ satisfying $\varphi\ge1$ near $E\cap\overline{B(y,\frac r2)}$,
we find that $[\varphi]_{s,p}^p\gtrsim r^{d-sp}$.
By the definition and monotonicity of $\operatorname{cap}_{s,p}$, we conclude that
\eqref{eq:thickE} holds. This completes the proof of
Lemma~\ref{lem:thickE}.
\end{proof}

Next, we show the range of $s$ in Theorem
\ref{thm:generalp} is optimal, namely Theorem~\ref{thm:counterexample}.

\begin{proof}[Proof of Theorem~\ref{thm:counterexample}]
Let $p\in(1,d-1)$, $s\in[1+\frac1p,\frac dp)$,
$m:=\lfloor s\rfloor$, and $\sigma:=s-m\in[0,1)$. Let
$q:=d-p(s-1)\in(p,d-1].$
Let $Q:=[0,\frac{1}{8\sqrt{d-1}}]^{d-1}\times\{0\}$.
Using the $(d-1)$-regularity of $Q$ and
\cite[Theorem~3.1]{MattilaSaaranen2009}, we find that, for
$q\in(p,d-1)$, there exists a nonempty compact $q$-regular set
$E\subset Q$.
For $q=d-1$, let $E:=Q$.
Let $r_0\in(0,\operatorname{diam}(E))$.
By Lemma~\ref{lem:thickE}, we obtain \eqref{eq:thickE} for
any $y\in E$ and $r\in(0,r_0)$.
From $E\subset Q$, we infer that
\begin{align}\label{eq:Ehyperplane}
E\subset B_{\mathbb{R}^{d-1}}\left(0,\frac14\right)\times\{0\},
\end{align}
where $B_{\mathbb{R}^{d-1}}$ denotes the $(d-1)$-dimension ball.
For any $x\in\mathbb{R}^d$ and $r\in(0,\infty)$, let
$\rho_E(x):=\operatorname{dist}(x,E)$ and
$N_r(E):=\{x\in\mathbb{R}^d:\rho_E(x)<r\}$.
Let $\Omega:=B(0,3)\setminus E$ and
$$\alpha_0:=\min\left\{1,\frac{C}{3^{d-sp}\kappa_{d,s,p}}\right\},$$
where $C$ is the same as in Lemma~\ref{lem:thickE}.
For any $x\in B(0,1)$ satisfying $\rho_E(x)\in(0,r_0)$,
let $y\in E$ satisfy $|x-y|=\rho_E(x)$.
By $E\subset\mathbb{R}^d\setminus\Omega$, $B(y,\rho_E(x))\subset B(x,3\rho_E(x))$,
and \eqref{eq:thickE}, we have
\begin{align*}
\operatorname{cap}_{s,p}\left(\overline{\left(\mathbb{R}^d\setminus\Omega\right)\cap B(x,3\rho_E(x))}\right)
\ge C\rho_E(x)^{d-sp}
\ge\alpha_0\kappa_{d,s,p}\left[3\rho_E(x)\right]^{d-sp}.
\end{align*}
Using this and the definition of $d_{\alpha_0,s,p}$, we obtain, for any
$x\in B(0,1)$ satisfying $\rho_E(x)\in(0,r_0)$,
\begin{align}\label{eq:dlessrho}
\rho_E(x)\le d_{\alpha_0,s,p}(x)\le3\rho_E(x).
\end{align}

Let $\chi\in C_{\rm c}^\infty(B(0,2))$ satisfy $0\le\chi\le1$ and
$\chi(x)=1$ for any $x\in B(0,1)$.

For any $r\in(0,r_0)$, let $\{y_i\}_{i=1}^{M_r}\subset E$ be
maximal among families satisfying $|y_i-y_j|\ge r$ for any distinct
$i,j\in[1,M_r]\cap\mathbb Z$.
From the maximality of $\{y_i\}_{i=1}^{M_r}$ and the disjointness of
the balls $B(y_i,\frac r3)$, it follows that
\begin{align}\label{eq:counter-covering}
\bigcup_{i=1}^{M_r}B\left(y_i,\frac r3\right)
\subset N_r(E)\subset\bigcup_{i=1}^{M_r}B(y_i,2r),
\end{align}
which, together with the $q$-regularity of $E$, further implies that
\begin{align*}
M_r r^q
&\lesssim\sum_{i=1}^{M_r}\mathcal H^q\left(E\cap B\left(y_i,\frac r3\right)\right)
\le\mathcal H^q(E)
\le\sum_{i=1}^{M_r}\mathcal H^q(E\cap B(y_i,r))
\lesssim M_r r^q,
\end{align*}
and hence $M_r\sim r^{-q}$.
Using this, \eqref{eq:counter-covering}, the disjointness
of the balls
$\{B(y_i,\frac r3)\}_{i=1}^{M_r}$, and \eqref{eq:Ehyperplane}, we obtain, for any $r\in(0,r_0)$,
\begin{align}\label{eq:tubevolume}
|N_r(E)|\sim r^{d-q}
\end{align}
and
\begin{align}\label{eq:lower-tube-H}
\int_{\mathbb{R}^{d-1}}\mathbf1_{\{\rho_E(x',0)<r\}}\,dx'
\sim r^{d-1-q}.
\end{align}

Let $\omega\in C_{\rm c}^\infty(B(0,1))$ satisfy $\omega\ge0$ and
$\int_{\mathbb{R}^d}\omega(x)\,dx=1$.
For any $\varepsilon\in(0,\frac{r_0}{4})$ and $x\in\mathbb{R}^d$, let
\begin{align}\label{2138}
\eta_\varepsilon(x):=\int_{N_{\frac{3\varepsilon}{2}}(E)}
\left(\frac4\varepsilon\right)^d
\omega\left(\frac{4(x-y)}\varepsilon\right)\,dy.
\end{align}
By this definition, we have $\eta_\varepsilon\in C_{\rm c}^\infty(\mathbb{R}^d)$,
$0\le\eta_\varepsilon\le1$, $\operatorname{supp}\eta_\varepsilon\subset N_{2\varepsilon}(E)$,
and, for any $k\in\mathbb Z_+$ and $x\in N_\varepsilon(E)$,
\begin{align}\label{eq:cutoffprops}
\eta_\varepsilon(x)=1\quad\text{and}\quad
\sup_{y\in\mathbb{R}^d}\left|D^k\eta_\varepsilon(y)\right|\lesssim\varepsilon^{-k}.
\end{align}
From \eqref{eq:Ehyperplane}, $\operatorname{supp}\eta_\varepsilon\subset N_{2\varepsilon}(E)$,
and the definition of $\chi$, it follows that, for any
$\varepsilon\in(0,\frac{r_0}{4})$,
$\operatorname{supp}\eta_\varepsilon\subset B(0,1)$ and
$\chi\eta_\varepsilon=\eta_\varepsilon$.
For any $\varepsilon\in(0,\frac{r_0}{4})$ and $x:=(x_1,\ldots,x_d)\in\mathbb{R}^d$, let
\begin{align*}
u_\varepsilon(x):=x_d\chi(x)\left[1-\eta_\varepsilon(x)\right]
=x_d\chi(x)-x_d\eta_\varepsilon(x),
\end{align*}
and we immediately have $u_\varepsilon\in C_{\rm c}^\infty(\Omega)$.

We first estimate $[u_\varepsilon]_{s,p}$.
Using \eqref{eq:Ehyperplane},
$\operatorname{supp}\eta_\varepsilon\subset N_{2\varepsilon}(E)\subset B(0,1)$,
Leibniz's rule, \eqref{eq:cutoffprops}, \eqref{eq:tubevolume}, and
$d-q=p(s-1)$, we obtain, for any $\varepsilon\in(0,\frac{r_0}{4})$ and
$k\in\mathbb Z_+$,
\begin{align}\label{eq:dercutoff}
\left\|D^k(x_d\eta_\varepsilon)\right\|_{L^p(\mathbb{R}^d)}
&\lesssim\varepsilon^{1-k}|N_{2\varepsilon}(E)|^{\frac1p}
\lesssim\varepsilon^{1-k+\frac{d-q}p}=\varepsilon^{s-k}.
\end{align}
By this with $k=m$ and by the definition of $[\,\cdot\,]_{s,p}$,
we have, if $\sigma=0$, $[x_d\eta_\varepsilon]_{s,p}\lesssim1$.
If $\sigma\in(0,1)$, using \eqref{eq:dercutoff} with $k=m$
and $k=m+1$,
we obtain, for any $h\in\mathbb{R}^d$,
\begin{align*}
\|D^m(x_d\eta_\varepsilon)(\,\cdot+h)-D^m(x_d\eta_\varepsilon)\|_{L^p(\mathbb{R}^d)}
\lesssim\min\left\{\varepsilon^\sigma,|h|\varepsilon^{\sigma-1}\right\}.
\end{align*}
From this and the definition of $[\,\cdot\,]_{s,p}$, we deduce that,
if $\sigma\in(0,1)$,
\begin{align*}
[x_d\eta_\varepsilon]_{s,p}^p
&\lesssim\varepsilon^{(\sigma-1)p}
\int_{\{h\in\mathbb{R}^d:\,|h|<\varepsilon\}}|h|^{p-d-\sigma p}\,dh
+\varepsilon^{\sigma p}\int_{\{h\in\mathbb{R}^d:\,|h|\ge\varepsilon\}}|h|^{-d-\sigma p}\,dh
\sim1.
\end{align*}
By this and the case $\sigma=0$, we conclude that
\begin{align*}
\sup_{\varepsilon\in(0,\frac{r_0}{4})}[x_d\eta_\varepsilon]_{s,p}<\infty.
\end{align*}
Using this, the definition of $u_\varepsilon$,
the triangle inequality for $[\cdot]_{s,p}$, and
$x_d\chi\in C_{\rm c}^\infty(\mathbb{R}^d)$, we obtain
\begin{align}\label{eq:uniformueps}
\sup_{\varepsilon\in(0,\frac{r_0}{4})}[u_\varepsilon]_{s,p}
\le[x_d\chi]_{s,p}
+\sup_{\varepsilon\in(0,\frac{r_0}{4})}[x_d\eta_\varepsilon]_{s,p}<\infty.
\end{align}

Now, we estimate the Hardy integrals of $u_\varepsilon$.
Let $k_0\in\mathbb N$ satisfy
$4^{-k_0}<\frac{r_0}{\sqrt5}$.
For any $k\in[k_0,\infty)\cap\mathbb Z$, let
\begin{align*}
G_k:=\left\{(x',x_d)\in\mathbb{R}^{d-1}\times\mathbb{R}:
\rho_E(x',0)<4^{-k},\ x_d\in(4^{-k},2\cdot4^{-k})\right\}.
\end{align*}
From \eqref{eq:Ehyperplane} and the disjointness of the intervals
$\{(4^{-k},2\cdot4^{-k})\}_{k=k_0}^\infty$, we infer that $G_k\subset\Omega\cap B(0,1)$ and
$G_k\cap G_l=\varnothing$ for any distinct
$k,l\in[k_0,\infty)\cap\mathbb Z$.
Using \eqref{eq:Ehyperplane}, \eqref{eq:cutoffprops}, and
\eqref{eq:lower-tube-H}, we find that, for any
$\varepsilon\in(0,\frac{r_0}{4})$,
$k\in[k_0,\infty)\cap\mathbb Z$ satisfying $4^{-k}\ge2\varepsilon$, and
$x\in G_k$,
$4^{-k}<\rho_E(x)<\sqrt54^{-k}$,
$|u_\varepsilon(x)|=|x_d|\sim 4^{-k}$, and
$|G_k|\sim 4^{-k(d-q)}.$
By this, \eqref{eq:dlessrho}, and $d-q+p-sp=0$, we obtain
\begin{align*}
\int_{G_k}\frac{|u_\varepsilon(x)|^p}{d_{\alpha_0,s,p}(x)^{sp}}\,dx
\sim 4^{-k(p-sp)}|G_k|\sim 4^{-k(d-q+p-sp)}=1.
\end{align*}
From this, the disjointness of the sets $\{G_k\}_{k=k_0}^\infty$,
and the fact that
$\#\{k\in[k_0,\infty)\cap\mathbb Z:4^{-k}\ge2\varepsilon\}
\sim\log\frac1\varepsilon$, we deduce that
\begin{align*}
\int_\Omega\frac{|u_\varepsilon(x)|^p}{d_{\alpha_0,s,p}(x)^{sp}}\,dx
\gtrsim\log\frac1\varepsilon\to\infty,
\end{align*}
as $\varepsilon\to0^+$.
This, together with \eqref{eq:uniformueps},
completes the proof of Theorem~\ref{thm:counterexample}.
\end{proof}

\section{Extension to Relative Capacities}\label{sec:relative-capacity}

To enlarge the range of $s$, we replace absolute capacity by relative
capacity and consider the corresponding capacitary distance.
In Subsection \ref{6.1}, we prove
Theorem~\ref{thm:relative-capacity}, which gives the Hardy inequality
for the admissible order $s$ and the sharp dependence on $\alpha$.
In Subsection \ref{6.2}, we construct counterexamples for the
remaining order and prove
Theorem~\ref{thm:relative-classification}, thereby proving that the
admissible range of $s$ is optimal.

\subsection{Proof of Theorem~\ref{thm:relative-capacity}}\label{6.1}

We first show that the relative capacitary distance $\widetilde d_{\alpha,s,p}$
is also Lebesgue measurable.
Indeed, for any $(x,r)\in\mathbb R^d\times(0,\infty)$,
the change of variables $\xi:=\frac{y-x}{r}$ transforms
the outer ball $B(x,2r)$ into $B(0,2)$.
Under this change of variables, the image of
$\overline{F\cap B(x,r)}$ is
$\overline{\{\frac{y-x}{r}:y\in F\cap B(x,r)\}}$.
The map from $(x,r)$ to this normalized compact set is
Borel measurable.
Moreover, the map
$K\longmapsto\operatorname{cap}_{s,p}(K;B(0,2))$
is upper semicontinuous on the space of compact subsets
of $B(0,2)$. By the scaling and composition of Borel measurable maps,
we conclude that the admissibility set defined by \eqref{eq:relative-distance}
is Borel in $\mathbb R^d\times(0,\infty)$.
Applying \cite[Lemma~2.6]{CXYZ} to its restriction to
$\Omega\times(0,t)$ for any $t\in(0,\infty)$, we obtain the
Lebesgue measurability of
$\{x\in\Omega:\widetilde d_{\alpha,s,p}(x)<t\}$,
and hence of $\widetilde d_{\alpha,s,p}$.

Let $\widetilde\kappa_{d,s,p}
:=\operatorname{cap}_{s,p}(\overline{B(0,1)};B(0,2)).$
By translation invariance and scaling, we have
$$\operatorname{cap}_{s,p}\left(\overline{B(x,r)};B(x,2r)\right)
=\widetilde\kappa_{d,s,p}r^{d-sp}$$
for every $x\in\mathbb R^d$ and $r>0$.
The following two lemmas give some basic properties of $\widetilde d_{\alpha,s,p}$

\begin{lemma}\label{lem:relative-localization}
Let $p\in(1,\infty)$, $s\in(0,2)$, and $L\in[2,\infty)$.
Then there exists a positive constant $C$ such that, for any
$x\in\mathbb R^d$, any $r\in(0,\infty)$, and any compact set
$K\subset\overline{B(x,r)}$,
\begin{align}\label{eq:relative-localization}
\operatorname{cap}_{s,p}(K;B(x,Lr))
\le\operatorname{cap}_{s,p}(K;B(x,2r))
\le C\operatorname{cap}_{s,p}(K;B(x,Lr))
\end{align}
and, if $sp<d$, then
\begin{align}\label{eq:absolute-relative-comparison}
\operatorname{cap}_{s,p}(K)
\le\operatorname{cap}_{s,p}(K;B(x,2r))
\le C\operatorname{cap}_{s,p}(K).
\end{align}
\end{lemma}

\begin{proof}
By translation invariance and scaling, it suffices to prove
\eqref{eq:relative-localization} and
\eqref{eq:absolute-relative-comparison} for $x=0$ and $r=1$.
We first show that \eqref{eq:relative-localization} holds.
Let $\chi\in C_{\rm c}^\infty(B(0,2))$ satisfy $\chi=1$ on an open
neighbourhood of $\overline{B(0,1)}$ and $0\le\chi\le1$ on
$\mathbb R^d$.
Let $\sigma\in(0,1)$ and any $g\in C_{\rm c}^\infty(\mathbb R^d)$.
Using $\chi=0$ on $\mathbb R^d\setminus B(0,2)$,
$|\chi(y)-\chi(z)|\lesssim|y-z|$ for any $y,z\in\mathbb R^d$,
and convexity, we obtain
\begin{align}\label{eq:relative-cutoff-difference}
[\chi g]_{\sigma,p}^p
&\le2^{p-1}\iint_{B(0,4)\times B(0,4)}
\frac{|\chi(y)|^p|g(y)-g(z)|^p}{|y-z|^{d+\sigma p}}\,dy\,dz\notag\\
&\quad+2^{p-1}\iint_{B(0,4)\times B(0,4)}
\frac{|g(z)|^p|\chi(y)-\chi(z)|^p}{|y-z|^{d+\sigma p}}\,dy\,dz\notag\\
&\quad+2\int_{B(0,2)}\int_{\mathbb R^d\setminus B(0,4)}
\frac{|\chi(y)g(y)|^p}{|y-z|^{d+\sigma p}}\,dz\,dy\notag\\
&\lesssim[g]_{\sigma,p}^p+\|g\|_{L^p(B(0,4))}^p
\left[\int_{B(0,8)}|h|^{p-d-\sigma p}\,dh
+\int_{\mathbb R^d\setminus B(0,2)}|h|^{-d-\sigma p}\,dh\right]\notag\\
&\lesssim[g]_{\sigma,p}^p+\|g\|_{L^p(B(0,4))}^p.
\end{align}
By the fundamental theorem of calculus and the triangle inequality
for $\|\cdot\|_{L^p(\mathbb R^d)}$,
we have, for any $h\in\mathbb R^d$,
$$\left\|g(\cdot+h)-g\right\|_{L^p(\mathbb R^d)}
\le\min\left\{|h|\|\nabla g\|_{L^p(\mathbb R^d)},
2\|g\|_{L^p(\mathbb R^d)}\right\},$$
which further implies that
\begin{align}\label{eq:relative-lower-order-difference}
[g]_{\sigma,p}^p
&\le\|\nabla g\|_{L^p(\mathbb R^d)}^p
\int_{B(0,1)}|h|^{p-d-\sigma p}\,dh
+2^p\|g\|_{L^p(\mathbb R^d)}^p
\int_{\mathbb R^d\setminus B(0,1)}|h|^{-d-\sigma p}\,dh\notag\\
&\lesssim\|\nabla g\|_{L^p(\mathbb R^d)}^p
+\|g\|_{L^p(\mathbb R^d)}^p.
\end{align}
For $s\in(1,2)$, let $\sigma:=s-1$.
Using the product rule, \eqref{eq:relative-cutoff-difference}, and
\eqref{eq:relative-lower-order-difference}, we find that, for any
$v\in C_{\rm c}^\infty(\mathbb R^d)$,
\begin{align*}
[\chi v]_{s,p} &\le[\chi\nabla v]_{\sigma,p}+[v\nabla\chi]_{\sigma,p} \lesssim[v]_{s,p}+\|\nabla v\|_{L^p(B(0,4))} +\|v\nabla\chi\|_{L^p(\mathbb R^d)} +\|\nabla(v\nabla\chi)\|_{L^p(\mathbb R^d)}\\
&\lesssim[v]_{s,p}+\|v\|_{L^p(B(0,4))} +\|\nabla v\|_{L^p(B(0,4))}.
\end{align*}
From this, \eqref{eq:relative-cutoff-difference} with $\sigma:=s$
when $s\in(0,1)$, and the product rule when $s=1$, it follows that,
for any $v\in C_{\rm c}^\infty(\mathbb R^d)$,
\begin{align}\label{eq:relative-common-cutoff}
[\chi v]_{s,p}\lesssim[v]_{s,p}+
\begin{cases}
\|v\|_{L^p(B(0,4))},&s\in(0,1],\\
\|v\|_{L^p(B(0,4))}+\|\nabla v\|_{L^p(B(0,4))},&s\in(1,2).
\end{cases}
\end{align}
Let $v\in C_{\rm c}^\infty(B(0,L))$.
By $v=0$ on $B(0,3L)\setminus B(0,2L)$, we conclude that,
for any $\sigma\in(0,1)$,
\begin{align}\label{919}
[v]_{\sigma,p}^p
\ge2\int_{B(0,L)}\int_{B(0,3L)\setminus B(0,2L)}
\frac{|v(y)|^p}{|y-z|^{d+\sigma p}}\,dz\,dy
\gtrsim\|v\|_{L^p(\mathbb R^d)}^p.
\end{align}
For $s=1$, the Poincar\'e inequality gives
$\|v\|_{L^p(\mathbb R^d)}\lesssim[v]_{1,p}.$
For $s\in(1,2)$, using the Poincar\'e inequality and \eqref{919}
with $\sigma:=s-1$, we obtain
\begin{align}\label{922}
\|v\|_{L^p(\mathbb R^d)}+\|\nabla v\|_{L^p(\mathbb R^d)}
\lesssim[\nabla v]_{s-1,p}=[v]_{s,p}.
\end{align}
From these estimates and \eqref{eq:relative-common-cutoff}, we deduce that
$[\chi v]_{s,p}\lesssim[v]_{s,p}$.
Using this, the definition of relative capacity, and taking the infimum
over these $v$, we obtain
$$\operatorname{cap}_{s,p}(K;B(0,2))
\lesssim\operatorname{cap}_{s,p}(K;B(0,L)).$$
By the monotonicity of $\operatorname{cap}_{s,p}$ in
the outer open set, we conclude that
\eqref{eq:relative-localization} holds.

Next, we prove \eqref{eq:absolute-relative-comparison}.
Assume $sp<d$, and let $v\in C_{\rm c}^\infty(\mathbb R^d)$.
From H\"older's inequality on $B(0,4)$ and the homogeneous Sobolev
inequality, it follows that
\begin{align*}
\|v\|_{L^p(B(0,4))}
\lesssim\|v\|_{L^{\frac{dp}{d-sp}}(\mathbb R^d)}
\lesssim[v]_{s,p}.
\end{align*}
For $s\in(1,2)$, still using H\"older's inequality on $B(0,4)$ and the homogeneous Sobolev
inequality, we obtain
\begin{align*}
\|\nabla v\|_{L^p(B(0,4))}
\lesssim\|\nabla v\|_{L^{\frac{dp}{d-(s-1)p}}(\mathbb R^d)}
\lesssim[\nabla v]_{s-1,p}=[v]_{s,p}.
\end{align*}
From these estimates and \eqref{eq:relative-common-cutoff}, we infer
that $[\chi v]_{s,p}\lesssim[v]_{s,p}$.
Using the definition
of the capacity $\operatorname{cap}_{s,p}$ and taking the infimum over all such $v$ satisfying
$v\ge1$ on an open neighbourhood of $K$, we obtain
$\operatorname{cap}_{s,p}(K;B(0,2))\lesssim\operatorname{cap}_{s,p}(K)$.
By this and the definition
of the capacity $\operatorname{cap}_{s,p}$ again, we conclude that
\eqref{eq:absolute-relative-comparison} holds. This completes the proof of
Lemma~\ref{lem:relative-localization}.
\end{proof}

\begin{lemma}\label{lem:supercritical-distance}
Let $p\in(1,\infty)$ and $s\in(\frac dp,\infty)$.
Then there exists $\alpha_0\in(0,1)$ such that, for any
open set $\Omega\subset\mathbb R^d$, any $\alpha\in(0,\alpha_0]$,
and any $x\in\Omega$,
\begin{align*}
\widetilde d_{\alpha,s,p}(x)=\operatorname{dist}(x,F).
\end{align*}
\end{lemma}

\begin{proof}
Using Poincar\'e inequalities, the definition of
$[\,\cdot\,]_{\tau,p}$ when $\tau\in\mathbb N$, and \eqref{919}
when $\tau\notin\mathbb N$, we obtain,
for any $\tau\in(0,\infty)$ and $v\in C_{\rm c}^\infty(B(0,2))$,
\begin{align}\label{1109}
\sum_{j=0}^{\lfloor \tau\rfloor}\|D^jv\|_{L^p(\mathbb R^d)}
\lesssim\|D^{\lfloor \tau\rfloor}v\|_{L^p(\mathbb R^d)}
\lesssim[v]_{\tau,p}.
\end{align}
From this and the Sobolev--Morrey embedding,
we infer that, for any $v\in C_{\rm c}^\infty(B(0,2))$,
\begin{align*}
\|v\|_{L^\infty(\mathbb R^d)}
\lesssim
\sum_{j=0}^{\lfloor s\rfloor}\|D^jv\|_{L^p(\mathbb R^d)}+[v]_{s,p}
\lesssim[v]_{s,p}.
\end{align*}
Using this and scaling, we find that, for any $x\in\mathbb R^d$,
 $r\in(0,\infty)$, $y\in\overline{B(x,r)}$, and
$v\in C_{\rm c}^\infty(B(x,2r))$ satisfying $v\ge1$ on an open
neighbourhood of $\{y\}$,
\begin{align*}
1\le|v(y)|\le\|v\|_{L^\infty(\mathbb R^d)}
\lesssim r^{s-\frac dp}[v]_{s,p}.
\end{align*}
From this and the definition of relative capacity, we deduce that
there exists a positive constant $c$ such that, for any
$x\in\mathbb R^d$, $r\in(0,\infty)$, and
$y\in\overline{B(x,r)}$,
\begin{align}\label{eq:relative-point-capacity}
\operatorname{cap}_{s,p}(\{y\};B(x,2r))\ge cr^{d-sp}.
\end{align}
Let $\alpha_0:=\min\{\frac12,\frac{c}{2\widetilde\kappa_{d,s,p}}\}$.
For any open set $\Omega\subset\mathbb R^d$ with
$F:=\mathbb R^d\setminus\Omega\ne\varnothing$, any
$\alpha\in(0,\alpha_0]$, and any $x\in\Omega$,
using $F\cap B(x,r)=\varnothing$ when
$r\in(0,\operatorname{dist}(x,F)]$, we obtain
$\widetilde d_{\alpha,s,p}(x)\ge\operatorname{dist}(x,F)$.
For any $r\in(\operatorname{dist}(x,F),\infty)$, let
$y\in F\cap B(x,r)$.
By \eqref{eq:relative-point-capacity}, the monotonicity of
$\operatorname{cap}_{s,p}$, and the definition
of $\alpha_0$, we conclude that
\begin{align*}
\operatorname{cap}_{s,p}\left(\overline{F\cap B(x,r)};B(x,2r)\right)
&\ge\operatorname{cap}_{s,p}(\{y\};B(x,2r))
\ge cr^{d-sp}\ge\alpha\widetilde\kappa_{d,s,p}r^{d-sp}.
\end{align*}
Using this and the definition of $\widetilde d_{\alpha,s,p}$,
we obtain $\widetilde d_{\alpha,s,p}(x)\le r$.
Letting $r\to\operatorname{dist}(x,F)^+$, we find that
$\widetilde d_{\alpha,s,p}(x)\le\operatorname{dist}(x,F)$,
which, together with the reverse inequality, implies
the desired equality.
For $F=\varnothing$, using $\inf\varnothing=\infty$ and
$\operatorname{dist}(x,\varnothing)=\infty$, we obtain the same conclusion.
This completes the proof of Lemma~\ref{lem:supercritical-distance}.
\end{proof}

\begin{remark}
Let $p\in(1,\infty)$ and $s\in(0,1+\frac1p)$ satisfy $sp>d$,
and let $\alpha_0$ be the same as in
Lemma~\ref{lem:supercritical-distance}.
By Lemma~\ref{lem:supercritical-distance}, we conclude that, for any
$\alpha\in(0,\alpha_0]$, \eqref{eq:relative-main} is the Hardy inequality
with the Euclidean distance to the complement and
$\widetilde{\mathcal H}_{d,s,p}(\alpha)
=\widetilde{\mathcal H}_{d,s,p}(\alpha_0)$.
For $s=1$, this reduces to the classical Hardy inequality on
arbitrary open sets with $p>d$; see, for instance,
\cite[(3.2)]{KinnunenMartio1997}.
\end{remark}

By \eqref{eq:absolute-relative-comparison} and the definition of the
(relative) capacitary distance, we find that there exists $A\in[1,\infty)$ such that, if $sp<d$,
then, for any $x\in\Omega$ and $\alpha\in(0,\infty)$,
\begin{align*}
d_{\frac\alpha A,s,p}(x)\le\widetilde d_{\alpha,s,p}(x)
\le d_{A\alpha,s,p}(x).
\end{align*}
From this and the definitions of the best constants, we infer that
\begin{align}\label{eq:relative-best-comparison}
\mathcal H_{d,s,p}(A\alpha)
\le\widetilde{\mathcal H}_{d,s,p}(\alpha)
\le\mathcal H_{d,s,p}\left(\frac\alpha A\right).
\end{align}

Now, we prove Theorem~\ref{thm:relative-capacity}(i).

\begin{proof}[Proof of Theorem~\ref{thm:relative-capacity}(i)]
In the case $s\in(0,\min\{\frac dp,1+\frac1p\})$,
the desired conclusion follows from
\eqref{eq:relative-best-comparison} and Theorem~\ref{thm:generalp}(i).

Now, let $s\in[\frac dp,1+\frac1p)$ (if this interval is nonempty). We first consider the case $s=1$.
Let $f\in W^{1,p}_0(\Omega)$ satisfy $0\le f\le1$,
and let $v_{2^j}$ and $S_f$ be the same as in
Lemma~\ref{lem:resolvent-sum}.
For any $(j,z)\in\mathbb Z\times\mathbb Q^d$, let
$E_{j,z}:=F\cap\overline{B(z,2^{\frac jp})}$ and
$h_{j,z}$ be the minimizer in Lemma~\ref{lem:local-barrier}
with $E:=E_{j,z}$ and $R:=2^{\frac jp}$.
By Lemma~\ref{lem:local-comparison}, the definition of the essential
infimum, and the countability of $\mathbb Z\times\mathbb Q^d$,
we conclude that there exists a Lebesgue-null set $N_1\subset\Omega$
such that, for any $(j,z)\in\mathbb Z\times\mathbb Q^d$ and
$y\in[\Omega\cap B(z,4\cdot2^{\frac jp})]\setminus N_1$,
$$v_{2^j}(y)\le1-h_{j,z}(y)\quad \mathrm{and}\quad
h_{j,z}(y)\ge\operatorname*{ess\,inf}_{w\in
B(z,4\cdot2^{\frac jp})}h_{j,z}(w).$$

Let $x\in\Omega\setminus N_1$ satisfy $f(x)=1$ and
$\widetilde d_{\alpha,1,p}(x)<\infty$,
$r\in[\widetilde d_{\alpha,1,p}(x),2\widetilde d_{\alpha,1,p}(x))$
be admissible in \eqref{eq:relative-distance},
$j\in\mathbb Z$ satisfy
$4\widetilde d_{\alpha,1,p}(x)\le R:=2^{\frac jp}
<2^{2+\frac1p}\widetilde d_{\alpha,1,p}(x)$, and
$z\in\mathbb Q^d$ satisfy $|z-x|<\widetilde d_{\alpha,1,p}(x)$.
From $r+|z-x|<3\widetilde d_{\alpha,1,p}(x)<R$, we infer that
$\overline{F\cap B(x,r)}\subset E_{j,z}$.
By the admissibility of $r$,
$\overline{F\cap B(x,r)}\subset E_{j,z}$,
\eqref{eq:relative-localization},
the Sobolev characterization of relative capacity,
\eqref{eq:barrier-energy-oscillation}, and
Lemma~\ref{lem:local-comparison}, we have
\begin{align*}
\alpha r^{d-p}
&\lesssim\operatorname{cap}_{1,p}
\left(\overline{F\cap B(x,r)};B(x,2r)\right)
\lesssim\operatorname{cap}_{1,p}(E_{j,z};B(z,8R))\\
&\le[h_{j,z}]_{1,p}^p
\le\mathcal J_{2^j}(h_{j,z};0)
\lesssim R^{d-p}\left[1-v_{2^j}(x)\right]^{p-1}.
\end{align*}
From this, the definition of $S_f$, and
$r\sim R\sim\widetilde d_{\alpha,1,p}(x)$, we deduce that
\begin{align*}
S_f(x)\ge\frac{|1-v_{2^j}(x)|^p}{2^j}
\gtrsim\frac{\alpha^{p'}}{\widetilde d_{\alpha,1,p}(x)^p},
\end{align*}
which, together with \eqref{eq:resolvent-sum}, further implies that
\begin{align*}
\int_{\{x\in\Omega:f(x)=1\}}
\frac{1}{\widetilde d_{\alpha,1,p}(x)^p}\,dx
\lesssim\alpha^{-p'}\int_\Omega S_f(x)\,dx
\lesssim\alpha^{-p'}[f]_{1,p}^p.
\end{align*}
Using Lemma~\ref{lem:plateau}, we obtain the upper bound with
$\alpha^{-p'}$ for $s=1$.

We next consider the case $s\ne1$. Let $u\in C_{\rm c}^\infty(\Omega)$.
Let $x\in\Omega$ satisfy $\widetilde d_{\alpha,s,p}(x)\in(0,\infty)$.
From the definition of $\widetilde d_{\alpha,s,p}(x)$, we infer that there exists
$r_x\in[\widetilde d_{\alpha,s,p}(x),2\widetilde d_{\alpha,s,p}(x))$ satisfying
$$\frac{\operatorname{cap}_{s,p}(\overline{F\cap B(x,r_x)};B(x,2r_x))}
{\operatorname{cap}_{s,p}(\overline{B(x,r_x)};B(x,2r_x))}\ge\alpha.$$
Using the adjacent-grid covering property, we find that there exist
$t\in[1,3^d]\cap\mathbb Z$ and $Q\in\mathscr D_t$ satisfying
$\overline{B(x,r_x)}\subset Q^\circ$ and
$2r_x\le\ell(Q)\le C_{\mathrm{cov}}r_x$.
Let
$C_{\mathrm{loc}}\in[1,\infty)$ be the constant in
\eqref{eq:relative-localization} with $L:=3\sqrt d\,C_{\mathrm{cov}}$.
By $(5Q)^\circ\subset B(x,Lr_x)$, monotonicity in the outer domain,
\eqref{eq:relative-localization}, and the choice of $r_x$,
we have $\overline{F\cap B(x,r_x)}\subset F\cap Q$ and
\begin{align}\label{eq:relative-selection-cap}
\operatorname{cap}_{s,p}\left(\overline{F\cap B(x,r_x)};(5Q)^\circ\right)
&\ge\operatorname{cap}_{s,p}\left(\overline{F\cap B(x,r_x)};B(x,Lr_x)\right)\notag\\
&\ge C_{\mathrm{loc}}^{-1}
\operatorname{cap}_{s,p}\left(\overline{F\cap B(x,r_x)};B(x,2r_x)\right)\notag\\
&\ge\alpha C_{\mathrm{loc}}^{-1}\widetilde\kappa_{d,s,p}r_x^{d-sp}
\ge\alpha C_{\mathrm{loc}}^{-1}\widetilde\kappa_{d,s,p}\ell(Q)^{d-sp}.
\end{align}
Here the last inequality uses $d-sp\le0$ and $\ell(Q)\ge2r_x$.
Let $c:=\min\{1,C_{\mathrm{loc}}^{-1}\widetilde\kappa_{d,s,p}\}$ and
$\mathscr A$ be the family of all cubes $Q\in\mathscr Q$
for which there exists a compact set $E_Q\subset F\cap Q$ such that
$\operatorname{cap}_{s,p}(E_Q;(5Q)^\circ)\ge c\alpha\ell(Q)^{d-sp}$.
Using \eqref{eq:relative-selection-cap} and the choice of $Q$, we find that
there exists a positive constant $C$ such that each $x\in\Omega$
with $\widetilde d_{\alpha,s,p}(x)<\infty$ admits a cube $Q\in\mathscr A$
satisfying $x\in Q$ and
$C^{-1}\ell(Q)\le\widetilde d_{\alpha,s,p}(x)\le C\ell(Q)$.
Let $\{Q_m\}_{m\in\mathbb N}$ enumerate $\mathscr Q$ without repetitions.
For any $m\in\mathbb N$, let
$$G_m:=\left\{x\in Q_m\cap\Omega:C^{-1}\ell(Q_m)\le\widetilde d_{\alpha,s,p}(x)\le C\ell(Q_m)\right\}$$
if $Q_m\in\mathscr A$, and let $G_m:=\varnothing$ otherwise. Let
$A_{Q_m}:=G_m\setminus\bigcup_{k=1}^{m-1}G_k$.
Since $\widetilde d_{\alpha,s,p}$ is measurable, it follows that the sets
$\{A_{Q_m}\}_{m\in\mathbb N}$ are measurable, pairwise disjoint,
and satisfy $A_Q=\varnothing$ for any $Q\notin\mathscr A$.
Using the definition of $G_m$ and $A_{Q_m}\subset G_m$, we obtain,
for any $Q\in\mathscr Q$ and $x\in A_Q$,
\begin{align}\label{eq:relative-selection-scale}
C^{-1}\ell(Q)\le\widetilde d_{\alpha,s,p}(x)\le C\ell(Q).
\end{align}
From the existence of an admissible cube for each $x$, we deduce that
$$\bigcup_{m\in\mathbb N}A_{Q_m}
=\left\{x\in\Omega:\widetilde d_{\alpha,s,p}(x)\in(0,\infty)\right\}.$$
Using \eqref{eq:relative-selection-scale}, the disjointness of the sets $\{A_Q\}_{Q\in\mathscr A}$,
and $\infty^{-sp}=0$, we obtain
\begin{align}\label{eq:relative-selection-integral}
\int_\Omega\frac{|u(x)|^p}{\widetilde d_{\alpha,s,p}(x)^{sp}}\,dx
=\sum_{Q\in\mathscr A}\int_{A_Q}
\frac{|u(x)|^p}{\widetilde d_{\alpha,s,p}(x)^{sp}}\,dx
\lesssim\sum_{Q\in\mathscr A}\ell(Q)^{-sp}
\int_{A_Q}|u(x)|^p\,dx.
\end{align}
By \eqref{eq:tree-test-admissibility}, \eqref{eq:tree-test-energy}
with $Q_0:=Q$ and $E:=E_Q$, and the definition of relative capacity,
we have, for every $N\in(j(Q),\infty)\cap\mathbb Z$,
\begin{align*}
\operatorname{cap}_{s,p}(E_Q;(5Q)^\circ)\lesssim\mathsf C_N(E_Q;Q).
\end{align*}
Thus, \eqref{eq:dyadic-tree-thickness} holds with $c_1\alpha$ in place of
$\alpha$, for a fixed $c_1\in(0,c]$.
Using Proposition~\ref{prop:dyadic-hardy-energy} and
\eqref{eq:relative-selection-integral}, we obtain the upper bound
with $\alpha^{-p'}$ for $s\ne1$.

Combining the above two cases, we find that, for $s\in[\frac dp,1+\frac1p)$ and any $\alpha\in(0,1]$,
\begin{align}\label{eq:relative-preliminary-upper}
\int_\Omega\frac{|u(x)|^p}{\widetilde d_{\alpha,s,p}(x)^{sp}}\,dx
\lesssim\alpha^{-p'}[u]_{s,p}^p.
\end{align}
This gives the desired upper bound when $s=\frac dp$.
For $s\in(\frac dp,1+\frac1p)$, we now show that the constant
can be chosen independently of $\alpha$. Let
$\alpha_0$ be the same as in Lemma~\ref{lem:supercritical-distance}.
By Lemma~\ref{lem:supercritical-distance} and
\eqref{eq:relative-preliminary-upper} with $\alpha:=\alpha_0$,
we have, for any $u\in C_{\rm c}^\infty(\Omega)$,
\begin{align}\label{eq:relative-low-order-euclidean}
\int_\Omega\frac{|u(x)|^p}{\operatorname{dist}(x,F)^{sp}}\,dx
=\int_\Omega\frac{|u(x)|^p}{\widetilde d_{\alpha_0,s,p}(x)^{sp}}\,dx
\lesssim\alpha_0^{-p'}[u]_{s,p}^p\lesssim[u]_{s,p}^p,
\end{align}
which, together with
$\widetilde d_{\alpha,s,p}(x)\ge\operatorname{dist}(x,F)$,
further implies that, for every $\alpha\in(0,1]$,
\begin{align*}
\int_\Omega\frac{|u(x)|^p}{\widetilde d_{\alpha,s,p}(x)^{sp}}\,dx
\le\int_\Omega\frac{|u(x)|^p}{\operatorname{dist}(x,F)^{sp}}\,dx
\lesssim[u]_{s,p}^p.
\end{align*}

Finally, we consider the case $s\in\mathcal A_{s,p}\setminus
(0,\max\{\frac dp,1+\frac1p\})$. Let $k\in\mathbb N$ satisfy
$s\in(k+\frac dp,k+1+\frac1p)$.
Applying \cite[Theorem~5]{Avkhadiev2006} successively to the derivatives
of $u$ of order $j\in[0,k-1]\cap\mathbb Z$, with the parameters
$n:=d$ and $s:=(s-j)p>d$ therein, we obtain
\begin{align}\label{eq:relative-higher-order-euclidean}
\int_\Omega\frac{|u|^p}{\operatorname{dist}(x,F)^{sp}}\,dx
\lesssim\int_\Omega\frac{|D^ku|^p}{\operatorname{dist}(x,F)^{(s-k)p}}\,dx.
\end{align}
Using \eqref{eq:relative-low-order-euclidean} with $u$ and $s$ replaced by $D^ku$ and $s-k$,
we find that
\begin{align*}
\int_\Omega\frac{|D^ku|^p}{\operatorname{dist}(x,F)^{(s-k)p}}\,dx
\lesssim[D^ku]_{s-k,p}^p\sim[u]_{s,p}^p.
\end{align*}
From this, \eqref{eq:relative-higher-order-euclidean}, and
$\widetilde d_{\alpha,s,p}(x)\ge\operatorname{dist}(x,F)$ for any $x\in\Omega$, we infer that
\eqref{eq:relative-main} holds with a constant independent of $\alpha$
for $s\in(k+\frac dp,k+1+\frac1p)$.

Combining the above three ranges, we complete
the proof of Theorem~\ref{thm:relative-capacity}(i).
\end{proof}

\begin{lemma}\label{lem:relative-critical-ball}
Let $p\in(1,\infty)$ and $s\in(0,2)$ satisfy $sp=d$.
Then, for any $r\in(0,\infty)$ and $\varepsilon\in(0,\frac r2]$,
\begin{align}\label{eq:relative-critical-ball}
\operatorname{cap}_{s,p}
\left(\overline{B(0,\varepsilon)};B(0,r)\right)
\gtrsim\left(\log\frac r\varepsilon\right)^{1-p}.
\end{align}
\end{lemma}

\begin{proof}
Let $\varphi\in C_{\rm c}^\infty(B(0,r))$ satisfy $\varphi\ge1$ on an open
neighbourhood of $\overline{B(0,\varepsilon)}$, $\tau:=\frac\varepsilon r$,
and $\psi(x):=\varphi(rx)$ for any $x\in\mathbb R^d$.
By scaling and $sp=d$, we have $\psi\in C_{\rm c}^\infty(B(0,1))$,
$[\psi]_{s,p}=[\varphi]_{s,p}$, $\tau\in(0,\frac12]$, and
$\psi\ge1$ on $\overline{B(0,\tau)}$.
For $s=1$, from $p=d$ and \cite[(2.2.14)]{MazyaSobolev}, it follows that
\begin{align*}
[\psi]_{1,p}^p
\ge\operatorname{cap}_{1,p}\left(\overline{B(0,\tau)};B(0,1)\right)
\sim\left(\log\frac1\tau\right)^{1-p}.
\end{align*}
For $s\ne1$, by Bernstein's inequality, H\"older's inequality, and $sp=d$,
we obtain, for any $j\in\mathbb Z_+$,
\begin{align}\label{eq:critical-frequency-average}
\left|\frac1{|B(0,\tau)|}\int_{B(0,\tau)}\Delta_j\psi(x)\,dx\right|
\lesssim\min\left\{1,\frac1{(2^j\tau)^{s}}\right\}
2^{js}\|\Delta_j\psi\|_{L^p(\mathbb R^d)}.
\end{align}
Using \eqref{eq:critical-frequency-average} and
$\psi=S_0\psi+\sum_{j\in\mathbb Z_+}\Delta_j\psi$, applying
Young's inequality to $S_0$, and using H\"older's inequality on
$\ell^p(\mathbb Z_+)$, we obtain
\begin{align*}
1&\le\frac1{|B(0,\tau)|}\int_{B(0,\tau)}\psi(x)\,dx
\le\|S_0\psi\|_{L^\infty(\mathbb R^d)}
+\sum_{j\in\mathbb Z_+}
\left|\frac1{|B(0,\tau)|}\int_{B(0,\tau)}\Delta_j\psi(x)\,dx\right|\\
&\quad\lesssim\|\psi\|_{L^p(\mathbb R^d)}
+\left[\sum_{j\in\mathbb Z_+}\min\left\{1,\left(2^j\tau\right)^{-sp'}\right\}\right]^{\frac1{p'}}
\left[\sum_{j\in\mathbb Z_+}2^{jsp}
\|\Delta_j\psi\|_{L^p(\mathbb R^d)}^p\right]^{\frac1p}\\
&\quad\lesssim\|\psi\|_{L^p(\mathbb R^d)}
+\left(\log\frac1\tau\right)^{\frac1{p'}}
\|\psi\|_{\dot B^s_{p,p}(\mathbb R^d)}.
\end{align*}
From this, \eqref{eq:besovLP}, and the Poincar\'e inequalities
\eqref{919} and \eqref{922}, it follows that
$[\psi]_{s,p}^p\gtrsim(\log\frac1\tau)^{1-p}$.
Using this, $[\psi]_{s,p}=[\varphi]_{s,p}$,
$\log\frac1\tau=\log\frac r\varepsilon$, and taking the infimum over
all such $\varphi$, we obtain \eqref{eq:relative-critical-ball} in the case $s\ne1$.
This completes the proof of Lemma~\ref{lem:relative-critical-ball}.
\end{proof}

\begin{proof}[Proof of Theorem~\ref{thm:relative-capacity}(ii)]
Using \eqref{eq:relative-best-comparison} and
Theorem~\ref{thm:generalp}(ii), we obtain the desired conclusion for $s\in(0,\frac dp)$.

For $s=\frac dp$, let $E:=\{0\}$, and define $\varepsilon_N$, $r_j$,
$\eta_j$, $u_N$, $K_N$, and $\Omega_{\varepsilon_N}$
as in the proof of Theorem~\ref{thm:generalp}(ii) with
this choice of $E$.
Since $E_r=B(0,r)$ for any $r\in(0,\infty)$ and $sp=d$,
it follows that $|E_r|\sim r^d=r^{sp}$.
By Young's inequality and the definition of $\eta_j$, we have,
for any $j\in\mathbb N$ and $m\in\{0,1,2\}$,
\begin{align*}
\|D^m\eta_j\|_{L^p(\mathbb R^d)}
\lesssim r_j^{-m}|B(0,r_j)|^{\frac1p}
\sim r_j^{\frac dp-m}=r_j^{s-m}.
\end{align*}
Using this and the same estimates as in
\eqref{eq:compact-sharp-cutoffs} and
\eqref{eq:compact-sharp-profile}, we obtain
$u_N\in C_{\rm c}^\infty(\Omega_{\varepsilon_N})$,
$u_N=1$ on $K_N$, and
\begin{align}\label{aaa}
[u_N]_{s,p}^p\lesssim N^{1-p}.
\end{align}
By polar coordinates and the definition of $K_N$, we have,
for any $N\in[2,\infty)\cap\mathbb N$,
\begin{align}\label{bbb}
\int_{K_N}\frac{1}{|x|^d}\,dx
&=d|B(0,1)|\int_{\frac1{4^{2N}}}^{\frac1{4^{N+1}}}\frac1r\,dr
=d|B(0,1)|(N-1)\log4\sim N.
\end{align}
From $B(x,6|x|)\subset B(0,8|x|)$ and
Lemma~\ref{lem:relative-critical-ball} with $r:=8|x|$ and
$\varepsilon:=\frac{\varepsilon_N}8$, it follows that, for any $x\in K_N$,
\begin{align*}
&\operatorname{cap}_{s,p}
\left(\overline{\left(\mathbb R^d\setminus\Omega_{\varepsilon_N}\right)\cap B(x,3|x|)};
B(x,6|x|)\right)\\
&\quad\ge\operatorname{cap}_{s,p}
\left(\overline{B\left(0,\frac{\varepsilon_N}8\right)};B(0,8|x|)\right)
\gtrsim\left(\log\frac{64|x|}{\varepsilon_N}\right)^{1-p}
\gtrsim N^{1-p},
\end{align*}
which further implies that there exists $\gamma\in(0,1)$ such that,
for any $\alpha\le\gamma N^{1-p}$ and $x\in K_N$,
$\widetilde d_{\alpha,s,p}(x)\le3|x|$.
Let $\alpha_*:=\gamma3^{1-p}$ and, for
$\alpha\in(0,\alpha_*]$, let
$N:=\lfloor(\frac\gamma\alpha)^{\frac1{p-1}}\rfloor$.
Using \eqref{aaa} and \eqref{bbb},
we obtain, for $\alpha\in(0,\alpha_*]$,
\begin{align*}
\widetilde{\mathcal H}_{d,s,p}(\alpha)
\ge\frac{1}{[u_N]_{s,p}^p}3^{-d}\int_{K_N}\frac{1}{|x|^d}\,dx
\gtrsim N^p\sim\alpha^{-p'}.
\end{align*}
For $s\in(\frac dp,\infty)$, let $\Omega:=B(0,2)$ and
$u\in C_{\rm c}^\infty(B(0,1))$ satisfy $u\ne0$.
By Lemma~\ref{lem:supercritical-distance}, for every sufficiently small $\alpha$,
\begin{align*}
\widetilde{\mathcal H}_{d,s,p}(\alpha)
\ge\frac{1}{[u]_{s,p}^p}\int_{B(0,2)}\frac{|u(x)|^p}{(2-|x|)^{sp}}\,dx>0.
\end{align*}
This completes the proof of Theorem~\ref{thm:relative-capacity}(ii).
\end{proof}

\subsection{Proof of Theorem~\ref{thm:relative-classification}}\label{6.2}

The following relative capacity estimate will allow the regular-set
construction in Section~\ref{sec:counterexamples} at arbitrary orders.

\begin{lemma}\label{lem:classified-capacity}
Let $p\in(1,\infty)$, $s\in(0,\infty)$, $L\in(1,\infty)$, and
$q\in[0,d)$ satisfy $q>d-sp$.
Let $E\subset\mathbb R^d$ be a nonempty compact $q$-regular set
if $q>0$, and a singleton if $q=0$.
Then there exists a positive constant $C$ such that, for any
$y\in E$ and $r\in(0,\operatorname{diam}(E))$ when $q>0$,
and for any $y\in E$ and $r\in(0,\infty)$ when $q=0$,
\begin{align*}
\operatorname{cap}_{s,p}(E\cap\overline{B(y,r)};B(y,Lr))
\ge Cr^{d-sp}.
\end{align*}
\end{lemma}

\begin{proof}
For $q=0$, let $y\in E:=\{y\}$ and $r\in(0,\infty)$.
By $q>d-sp$, we have $sp>d$.
Using \eqref{eq:relative-point-capacity} with $x:=y$ and
$r:=\frac{Lr}{2}$ therein, we obtain
\begin{align*}
\operatorname{cap}_{s,p}(E\cap\overline{B(y,r)};B(y,Lr))
=\operatorname{cap}_{s,p}(\{y\};B(y,Lr))\gtrsim r^{d-sp}.
\end{align*}

Next, we consider the case $q\in(0,d)$.
Let $y\in E$, $r\in(0,\operatorname{diam}(E))$, and
$\nu:=\mathcal H^q|_{E\cap\overline{B(y,\frac r2)}}$.
From the $q$-regularity of $E$, we infer that \eqref{eq:nu-mass} holds
for any $z\in\mathbb R^d$ and $t\in(0,\infty)$.
Let $j_0\in\mathbb Z$ satisfy $2^{-j_0-1}<r\le2^{-j_0}$ and
$\varphi\in C_{\rm c}^\infty(B(y,Lr))$.
Using \eqref{1109} and scaling, we find that
\begin{align}\label{eq:classified-poincare}
\|\varphi\|_{L^p(\mathbb R^d)}\lesssim r^s[\varphi]_{s,p}.
\end{align}
Using \eqref{eq:counter-bernstein} and
\eqref{eq:counter-integer-frequency} when $s\in\mathbb N$, or
\eqref{eq:counter-bernstein} and \eqref{eq:besovLP} when $s\notin\mathbb N$, we obtain,
for any $j\in\mathbb Z$,
\begin{align}\label{eq:classified-frequency}
\|\Delta_j\varphi\|_{L^p(\mathbb R^d)}
\sim\frac1{2^{j\lfloor s\rfloor}}\left\|
\Delta_jD^{\lfloor s\rfloor}\varphi\right\|_{L^p(\mathbb R^d)}
\lesssim\frac1{2^{js}}[\varphi]_{s,p}.
\end{align}
From \eqref{eq:nu-mass}, Bernstein's inequality, the
$L^p$ boundedness of $\Delta_j$, and \eqref{eq:classified-poincare},
it follows that, for any $j\in(-\infty,j_0]\cap\mathbb Z$,
\begin{align}\label{eq:classified-low-frequency}
\int_{\mathbb R^d}|\Delta_j\varphi(z)|\,d\nu(z)
&\le\nu(\mathbb R^d)\|\Delta_j\varphi\|_{L^\infty(\mathbb R^d)}
\lesssim r^q\|\Delta_j\varphi\|_{L^\infty(\mathbb R^d)}\notag\\
&\lesssim r^q2^{\frac{jd}p}\|\Delta_j\varphi\|_{L^p(\mathbb R^d)}
\lesssim r^q2^{\frac{jd}p}\|\varphi\|_{L^p(\mathbb R^d)}
\lesssim r^{q+s}2^{\frac{jd}p}[\varphi]_{s,p}.
\end{align}
Using \eqref{eq:nu-frequency-bound} and \eqref{eq:classified-frequency}, we obtain,
for any $j\in(j_0,\infty)\cap\mathbb Z$,
\begin{align}\label{eq:classified-high-frequency}
\left|\int_{\mathbb R^d}\Delta_j\varphi(z)\,d\nu(z)\right|
\lesssim r^{\frac q{p'}}2^{\frac{j(d-q)}p}
\|\Delta_j\varphi\|_{L^p(\mathbb R^d)}
\lesssim\frac{r^{\frac q{p'}}}{2^{j(s-\frac{d-q}p)}}[\varphi]_{s,p}.
\end{align}
By $\varphi\in\mathcal S(\mathbb R^d)$ and the definition of $\Delta_j$,
we have $\sum_{j=-n}^{n}\Delta_j\varphi\to\varphi$ in
$L^\infty(\mathbb R^d)$ as $n\to\infty$.
From this, \eqref{eq:classified-low-frequency},
\eqref{eq:classified-high-frequency}, $s-\frac{d-q}p>0$,
and $2^{-j_0-1}<r\le2^{-j_0}$, we infer that
\begin{align}\label{eq:classified-trace}
\left|\int_{\mathbb R^d}\varphi(z)\,d\nu(z)\right|
&\le\sum_{j\in\mathbb Z}
\left|\int_{\mathbb R^d}\Delta_j\varphi(z)\,d\nu(z)\right|\notag\\
&\lesssim[\varphi]_{s,p}
\left[r^{q+s}\sum_{j\in(-\infty,j_0]\cap\mathbb Z}2^{\frac{jd}p}
+r^{\frac q{p'}}\sum_{j\in(j_0,\infty)\cap\mathbb Z}
\frac1{2^{j[s-\frac{d-q}p]}}\right]
\sim r^{q+s-\frac dp}[\varphi]_{s,p}.
\end{align}
For any $\varphi\in C_{\rm c}^\infty(B(y,Lr))$ satisfying
$\varphi\ge1$ on an open neighbourhood of $E\cap\overline{B(y,r)}$,
using \eqref{eq:nu-mass}, the definition of $\nu$, and
\eqref{eq:classified-trace}, we obtain
\begin{align*}
r^q\sim\nu(\mathbb R^d)
\le\int_{\mathbb R^d}\varphi(z)\,d\nu(z)
\lesssim r^{q+s-\frac dp}[\varphi]_{s,p}.
\end{align*}
By this and the definition of relative capacity, we
complete the proof of
Lemma~\ref{lem:classified-capacity}.
\end{proof}

Now, we prove Theorem~\ref{thm:relative-classification}.

\begin{proof}[Proof of Theorem~\ref{thm:relative-classification}]
By the definition of $\mathcal A_{d,p}$, we find that there exists
$k\in\mathbb N$ such that $s\in[k+\frac1p,k+\frac dp]$.
Let $q:=d-p(s-k)\in[0,d-1]$.
If $q=0$, let $E:=\{0\}$.
If $q\in(0,d-1]$, let
$Q:=[0,\frac{1}{8\sqrt{d-1}}]^{d-1}\times\{0\}$.
In the case $q=d-1$, let $E:=Q$. If $q\in(0,d-1)$, let
$E\subset Q$ be a nonempty compact $q$-regular set,
whose existence follows from the $(d-1)$-regularity of $Q$ and
\cite[Theorem~3.1]{MattilaSaaranen2009}.
Let $r_0:=\operatorname{diam}(E)$ when $q\in(0,d-1]$
and $r_0:=\frac18$ when $q=0$.
Let $\Omega:=B(0,3)\setminus E$, and
$\rho_E(x):=\operatorname{dist}(x,E)$ for any $x\in\mathbb R^d$.
By these definitions, we have $r_0\in(0,\frac18]$ and
\begin{align}\label{2220}
E\subset B\left(0,\frac14\right)\cap\{x\in\mathbb R^d:x_d=0\}.
\end{align}

For any $x\in B(0,1)$ satisfying $\rho_E(x)\in(0,r_0)$, let
$y\in E$ satisfy $|x-y|=\rho_E(x)$.
By monotonicity of relative capacity and
Lemma~\ref{lem:classified-capacity},
we obtain
\begin{align*}
\operatorname{cap}_{s,p}
\left(\overline{(\mathbb R^d\setminus\Omega)\cap
B(x,3\rho_E(x))};B(x,6\rho_E(x))\right)
&\ge\operatorname{cap}_{s,p}
\left(E\cap\overline{B(y,\rho_E(x))};B(y,7\rho_E(x))\right)\\
&\ge C\rho_E(x)^{d-sp},
\end{align*}
where $C$ is the constant in Lemma~\ref{lem:classified-capacity}
with $L:=7$.
Let $\alpha_0:=\min\{\frac12,
\frac{C}{2\widetilde\kappa_{d,s,p}3^{d-sp}}\}$.
Using this, scaling of the ball capacity, and the definition of
$\widetilde d_{\alpha_0,s,p}$, we obtain, for any
$x\in B(0,1)$ satisfying $\rho_E(x)\in(0,r_0)$,
\begin{align}\label{eq:classified-distance}
\rho_E(x)\le\widetilde d_{\alpha_0,s,p}(x)\le3\rho_E(x).
\end{align}
For any $r\in(0,\infty)$, let
$N_r(E):=\{x\in\mathbb R^d:\rho_E(x)<r\}$.
For $q\in(0,d-1]$, applying the covering argument based on
\eqref{eq:counter-covering} to the set $E$ chosen above, we obtain
\eqref{eq:tubevolume} and \eqref{eq:lower-tube-H} for any
$r\in(0,r_0)$. For $q=0$, we obtain \eqref{eq:tubevolume} from
$N_r(E)=B(0,r)$ and, for $q\in(0,d-1]$,
\eqref{eq:lower-tube-H} from
$\rho_E(x',0)=|x'|$ for any $x'\in\mathbb R^{d-1}$.
Let $\chi$ and $\omega$ be as in the proof of
Theorem~\ref{thm:counterexample}, and, for any
$\varepsilon\in(0,\frac{r_0}{4})$, let $\eta_\varepsilon$ be defined
by \eqref{2138} using the present set $E$. Then the same construction
gives \eqref{eq:cutoffprops} and
$\operatorname{supp}\eta_\varepsilon\subset N_{2\varepsilon}(E)$.
From $E\subset B(0,\frac14)$, $r_0\le\frac18$, and the definition
of $\chi$, it follows that $\operatorname{supp}\eta_\varepsilon\subset B(0,1)$
and $\chi\eta_\varepsilon=\eta_\varepsilon$.
For any $\varepsilon\in(0,\frac{r_0}4)$ and
$x:=(x_1,\ldots,x_d)\in\mathbb R^d$, let
\begin{align*}
u_\varepsilon(x):=x_d^k\chi(x)[1-\eta_\varepsilon(x)]
=x_d^k\chi(x)-x_d^k\eta_\varepsilon(x).
\end{align*}
By \eqref{eq:cutoffprops} and $\chi\in C_{\rm c}^\infty(B(0,2))$,
we have $u_\varepsilon\in C_{\rm c}^\infty(\Omega)$.

We first estimate $[u_\varepsilon]_{s,p}$.
Using $|x_d|\le2\varepsilon$ for any $x\in N_{2\varepsilon}(E)$,
Leibniz's rule, \eqref{eq:cutoffprops}, \eqref{eq:tubevolume}, and
$d-q=p(s-k)$, we obtain, for any
$\varepsilon\in(0,\frac{r_0}4)$ and $j\in\mathbb Z_+$,
\begin{align*}
\left\|D^j(x_d^k\eta_\varepsilon)\right\|_{L^p(\mathbb R^d)}
\lesssim\sum_{\ell=0}^{\min\{j,k\}}
\varepsilon^{k-\ell}\varepsilon^{-(j-\ell)}
|N_{2\varepsilon}(E)|^{\frac1p}
\lesssim\varepsilon^{k-j+\frac{d-q}p}=\varepsilon^{s-j}.
\end{align*}
Using this formula in place of
\eqref{eq:dercutoff} in the argument leading to
\eqref{eq:uniformueps}, we obtain
\begin{align}\label{eq:classified-uniform-energy}
\sup_{\varepsilon\in(0,\frac{r_0}4)}[u_\varepsilon]_{s,p}
\le[x_d^k\chi]_{s,p}
+\sup_{\varepsilon\in(0,\frac{r_0}4)}[x_d^k\eta_\varepsilon]_{s,p}<\infty.
\end{align}

Now, we estimate the Hardy integrals of $u_\varepsilon$.
Let $n_0\in\mathbb N$ satisfy $\frac1{4^{n_0}}<\frac{r_0}{\sqrt5}$.
For any $n\in[n_0,\infty)\cap\mathbb N$, let
\begin{align*}
G_n:=\left\{(x',x_d)\in\mathbb R^{d-1}\times\mathbb R:
\rho_E(x',0)<\frac1{4^n}\ \text{and}\ x_d\in\left(\frac1{4^n},\frac2{4^n}\right)\right\}
\end{align*}
when $d\ge2$, and let $G_n:=(\frac1{4^n},\frac2{4^n})$ when $d=1$.
From \eqref{2220} and
the disjointness of the intervals
$\{(\frac1{4^n},\frac2{4^n})\}_{n\in[n_0,\infty)\cap\mathbb N}$,
we infer that $G_n\subset\Omega\cap B(0,1)$ and
$G_n\cap G_\ell=\varnothing$ for any distinct
$n,\ell\in[n_0,\infty)\cap\mathbb N$.
Using \eqref{2220}, \eqref{eq:cutoffprops},
and \eqref{eq:lower-tube-H} when $d\ge2$, we find that, for any
$\varepsilon\in(0,\frac{r_0}4)$,
$n\in[n_0,\infty)\cap\mathbb N$ satisfying $\frac1{4^n}\ge2\varepsilon$,
and $x\in G_n$,
$\frac1{4^n}<\rho_E(x)<\frac{\sqrt5}{4^n}<r_0$,
$|u_\varepsilon(x)|=|x_d|^k\sim\frac1{4^{nk}}$ and
$|G_n|\sim\frac1{4^{n(d-q)}}.$
By this, \eqref{eq:classified-distance}, and $kp-sp+d-q=0$, we obtain
\begin{align*}
\int_{G_n}\frac{|u_\varepsilon(x)|^p}
{\widetilde d_{\alpha_0,s,p}(x)^{sp}}\,dx
\gtrsim\frac1{4^{n(kp-sp)}}|G_n|
\sim\frac1{4^{n(kp-sp+d-q)}}=1.
\end{align*}
From this and the disjointness of
$\{G_n\}_{n\in[n_0,\infty)\cap\mathbb N}$, we deduce that
\begin{align}\label{eq:classified-hardy-divergence}
\int_\Omega\frac{|u_\varepsilon(x)|^p}
{\widetilde d_{\alpha_0,s,p}(x)^{sp}}\,dx
\gtrsim\sum_{n=n_0}^{\left\lfloor\log_4\frac1{2\varepsilon}\right\rfloor}1
=\left\lfloor\log_4\frac1{2\varepsilon}\right\rfloor-n_0+1
\sim\log\frac1\varepsilon\to\infty
\end{align}
as $\varepsilon\to0^+$.
Using \eqref{eq:classified-uniform-energy} and
\eqref{eq:classified-hardy-divergence},
we complete the proof of Theorem~\ref{thm:relative-classification}.
\end{proof}

\smallskip

\noindent\textbf{Acknowledgements}\quad
The main mathematical ideas and approaches in this article originated from
the authors. Generative AI tools were used as supplementary aids during the exploratory and
manuscript-preparation stages, particularly in refining the presentation and examining technical
details. All mathematical arguments and proofs in the final manuscript were independently verified
and finalized by the authors, who take full responsibility for their correctness.

Finally, we would like to mention the preprint `Capacitary Hardy inequalities, 
Dirichlet boundary conditions, and the limit as $p\downarrow 1$' 
(arXiv:2611.20397v1) by Takanobu Hara. This article appeared on the same 
day when we submitted our article to the arXiv (November 17, 2026). 
Hara’s preprint has some nontrivial overlap with our results and contains, 
in particular, our Hardy inequalities for $s=1$ and $p\in (1, \infty)$. 
It also has results in the limit $p\downarrow 1$, that are not in our preprint, 
but, in contrast to our article, it does not have results for orders $s \neq 1$ 
and does not have our counterexamples. We believe that both independent 
articles complement each in a way that is valuable for the community.

\bigskip

\noindent Yiqun Chen, Dachun Yang, Wen Yuan, and Yangyang Zhang
%(Corresponding author)

\smallskip

\noindent  Laboratory of Mathematics and Complex Systems
(Ministry of Education of China),
School of Mathematical Sciences, Institute for Advanced Study,
Beijing Normal University,
Beijing 100875, The People's Republic of China

\smallskip

\noindent {\it E-mails}: \texttt{yiqunchen@mail.bnu.edu.cn} (Y. Chen)

\noindent\phantom{{\it E-mails:}} \texttt{dcyang@bnu.edu.cn} (D. Yang)

\noindent\phantom{{\it E-mails:}} \texttt{wenyuan@bnu.edu.cn} (W. Yuan)

\noindent\phantom{{\it E-mails:}} \texttt{yangyzhang@bnu.edu.cn} (Y. Zhang)

\bigskip

\noindent Rupert L. Frank

\smallskip

\noindent Mathematisches Institut, Ludwig-Maximilians
Universit\"at M\"unchen, Theresienstr. 39, 80333 M\"unchen,
Germany; Munich Center for Quantum Science and Technology,
Schellingstr. 4, 80799
M\"unchen, Germany; Mathematics 253-37, Caltech, Pasadena, CA 91125, USA

\smallskip

\noindent {\it E-mail}: \texttt{r.frank@lmu.de}

\bigskip

\noindent Jie Xiao

\smallskip

\noindent Department of Mathematics and Statistics, Memorial University, St. John's, NL A1C5S7, Canada

\smallskip

\noindent {\it E-mail}: \texttt{jxiao@mun.ca}

\end{document}